\documentclass[final]{siamltex}         

\usepackage{graphicx}

\usepackage{subfigure} 

\graphicspath{{C:/NS/Pictures/}}

\title{Solution of the Cauchy problem for the Navier - Stokes and Euler equations}                                       
\author{ A. Tsionskiy, M. Tsionskiy  }               

\newcommand{\sz}[1]{\mbox{\Large #1 \normalsize}}
\newcommand{\szp}[1]{\mbox{\large #1 \normalsize}}

\begin{document}             

\maketitle                   

\begin{abstract}
Some known results regarding the Euler and Navier-Stokes equations were obtained by different authors. 
Existence and smoothness of the Navier-Stokes solutions in two dimensions have been known for a long time. 
Leray $\cite{jL34}$ showed that the Navier-Stokes equations in three space dimensions have a weak solution. 
Scheffer $\cite{vS76}, \cite{vS93}$ and Shnirelman $\cite{aS97}$ obtained weak solution of the Euler equations with compact support in spacetime. Caffarelli-Kohn-Nirenberg $\cite{CKN82}$ improved Scheffer's results , and F.-H. Lin $\cite{fL98}$ simplified the proof of the results of J. Leray. Many problems and conjectures about the behavior of solutions of the Euler and Navier-Stokes equations are described in the book of Bertozzi and Majda $\cite{BM02}$ or Constantin $\cite{pC01}$.

Solutions of the Navier-Stokes and Euler equations with initial conditions (Cauchy problem) for two and three dimensions  are obtained in the convergence series form by the iterative method using the Fourier and Laplace transforms in this paper. For several combinations of problem parameters numerical results were obtained and presented as graphs.

\end{abstract}



\pagestyle{myheadings}
\thispagestyle{plain}
\markboth{A. TSIONSKIY, M. TSIONSKIY}{SOLUTION OF THE NAVIER-STOKES AND EULER EQUATIONS}

\section{The mathematical setup}\ 

The Navier-Stokes equations describe the motion of a fluid in \(R^{N} ( N=2\;\rm{or}\; 3 ) \). We look for a viscous incompressible fluid filling all of \(R^{N}\)  here. The Navier-Stokes equations are then given by

\begin{equation}\label{eqn1}
\frac{\partial u_{k}}{\partial t} \; + \; \sum_{n=1}^{N} u_{n}\frac{\partial u_{k}}{\partial x_{n}}\; =\;\nu\Delta u_{k}\; - \; \frac{\partial p}{\partial x_{k}} \; + \;f_{k}(x,t)\;\;\;\;\;(x\in R^{N},\;\;t\geq 0,\;\;{1\leq k \leq N}) 
\end{equation}
\begin{equation}\label{eqn2}
\emph{div}\,\vec{u}\;=\; \sum_{n=1}^{N} \frac{\partial u_{n}}{\partial x_{n}}\; =\;0\;\;\;\;\;\;\;\;\;\;  (x\in R^{N},t\geq 0)
\end{equation}

with initial conditions

\begin{equation}\label{eqn3}
\vec{u}(x,0)\; = \; \vec{u}^{0}(x)\;\;\;\;\;\;\;\;\;\; (x\in R^{N})
\end{equation}

Here \(\vec{u}(x,t)=(u_{k}(x,t)) \in R^{N},\;\; ({1\leq k \leq N})  \;-\; \)is an unknown velocity vector  \(( N = 2\; \rm{or}\; 3 ),\; p\,(x,t)\;-\;\) is an unknown pressure, \(\vec{u}^{0}(x)\;\) is a given, \(C^{\infty}\) divergence-free vector field \(,\; f_{k}(x,t)\;\)are components of a given, externally applied force \(\vec{f}(x,t)\),  \(\nu\) is a positive coefficient of the viscosity (if \(\nu = 0\) then (1.1) - (1.3) are the Euler equations), and \( \Delta\;=\; \sum_{n=1}^{N} \frac{\partial^{2}}{\partial x_{n}^{2}}\;\) is the Laplacian in the space variables. Equation $(\ref{eqn1})$ is Newton's law for a fluid element subject. Equation $(\ref{eqn2})$ says that the fluid is incompressible. For physically reasonable solutions, we accept 

\begin{equation}\label{eqn4}
u_{k}(x,t) \rightarrow 0\;\;, \;\;
\frac{\partial u_{k}}{\partial x_{n}} \;\rightarrow\;0\;\; 
\rm {as} 
\;\;\mid x \mid \;\rightarrow\; \infty\;\;\;( {1\leq k \leq N} ,\;\; {1\leq n \leq N}) \;\;\; 
\end{equation}

Hence, we will restrict attention to initial conditions $\vec{u}^{0}$ and force $\vec{f}$ that satisfy

\begin{equation}\label{eqn5}
\mid\partial_{x}^{\alpha}\vec{u}^{0}(x)\mid\;\leq\;C_{\alpha K}(1+\mid x \mid)^{-K} \quad \rm{on }\;R^{N}\;\rm{ for\; any }\;\alpha \;\rm{ and }\;K>0.
\end{equation}
and
\begin{equation}\label{eqn6}
\mid\partial_{x}^{\alpha}\partial_{t}^{\beta}\vec{f}(x,t)\mid\;\leq\;C_{\alpha \beta K}(1+\mid x \mid +t)^{-K} \quad \rm{on }\;R^{N}\times[0,\infty)\; \rm{ for \;any }\;\alpha,\beta \;\rm{ and }\;K>0.
\end{equation}

We add (\( - \sum_{n=1}^{N} u_{n}\frac{\partial u_{k}}{\partial x_{n}}\;\)) to both sides of the equations (\ref{eqn1}). Then we have:

\begin{equation}\label{eqn7}
\frac{\partial u_{k}}{\partial t}\;=\;\nu\,\Delta\,u_{k}\;-\;\frac{\partial p}{\partial x_{k}}\;+\;f_{k}(x,t)- \; \sum_{n=1}^{N} u_{n}\frac{\partial u_{k}}{\partial x_{n}}\;\;\;\;\;\;\;\;\;\;\;\; (x\in R^{N},\;\;t\geq 0,\;\;{1\leq k \leq N})
\end{equation}
\begin{equation}\label{eqn8}
\emph{div}\,\vec{u}\;=\; \sum_{n=1}^{N} \frac{\partial u_{n}}{\partial x_{n}}\; =\;0\;\;\;\;\;\;\;\;\;\;  (x\in R^{N},t\geq 0)
\end{equation}
\begin{equation}\label{eqn9}
\vec{u}(x,0)\; = \; \vec{u}^{0}(x)\;\;\;\;\;\;\;\;\;\; (x\in R^{N})
\end{equation}

\begin{equation}\label{eqn10}
u_{k}(x,t) \rightarrow 0\;\;, \;\;\frac{\partial u_{k}}{\partial x_{n}}\;\rightarrow\;0\;\; \rm{as} \;\;\mid x \mid \;\rightarrow\; \infty\;\;\;( {1\leq k \leq N} ,\;\; {1\leq n \leq N}) \;\;\; 
\end{equation}

\begin{equation}\label{eqn11}
\mid\partial_{x}^{\alpha}\vec{u}^{0}(x)\mid\;\leq\;C_{\alpha K}(1+\mid x \mid)^{-K} \quad\rm{on }\;R^{N}\;\rm{ for \; any }\;\alpha\;\rm{ and }\;K>0.
\end{equation}

\begin{equation}\label{eqn12}
\mid\partial_{x}^{\alpha}\partial_{t}^{\beta}\vec{f}(x,t)\mid\;\leq\;C_{\alpha \beta K}(1+\mid x \mid +t)^{-K} \quad\rm{on }\;R^{N}\times[0,\infty)\;\rm{ for \; any }\;\alpha,\beta\;\rm{ and }\;K>0.
\end{equation}

We shall solve the system of equations (\ref{eqn7}) - (\ref{eqn12}) by the iterative method. To do so we write this system of equations in the following form:

\begin{equation}\label{eqn13}
\frac{\partial u_{jk}}{\partial t}\;=\;\nu\,\Delta\,u_{jk}\;-\;\frac{\partial p_{j}}{\partial x_{k}}\;+\;f_{jk}(x,t) \;\;\;\;\;\;\;\;\;\;\;\; (x\in R^{N},\;\;t\geq 0,\;\;{1\leq k \leq N})
\end{equation}
\begin{equation}\label{eqn14}
\emph{div}\,\vec{u}_{j}\;=\; \sum_{n=1}^{N} \frac{\partial u_{jn}}{\partial x_{n}}\; =\;0\;\;\;\;\;\;\;\;\;\;  (x\in R^{N},t\geq 0)
\end{equation}
\begin{equation}\label{eqn15}
\vec{u}_{j}(x,0)\; = \; \vec{u}^{0}(x)\;\;\;\;\;\;\;\;\;\; (x\in R^{N})
\end{equation}

\begin{equation}\label{eqn16}
u_{j k}(x,t) \rightarrow 0\;\;, \;\;\frac{\partial u_{j k}}{\partial x_{n}}\;\rightarrow\;0\;\; \rm{as} \;\;\mid x \mid \;\rightarrow\; \infty\;\;\;( {1\leq k \leq N} ,\;\; {1\leq n \leq N}) \;\;\; 
\end{equation}

\begin{equation}\label{eqn17}
\mid\partial_{x}^{\alpha}\vec{u}^{0}(x)\mid\;\leq\;C_{\alpha K}(1+\mid x \mid)^{-K} \quad\rm{on }\;R^{N}\;\rm{ for \; any }\;\alpha\;,\;K>0\;\;\;\rm{ and } \;\;C_{\alpha K} > 0.
\end{equation}

\begin{equation}\label{eqn18}
\mid\partial_{x}^{\alpha}\partial_{t}^{\beta}\vec{f}(x,t)\mid\;\leq\;C_{\alpha \beta K}(1+\mid x \mid +t)^{-K} \quad\rm{on }\;R^{N}\times[0,\infty)\;\rm{ for \; any }\;\alpha,\beta\;,\;K>0\;\;\rm{ and } \;\;C_{\alpha \beta K} > 0.
\end{equation}

Here j is the number of the iterative process step (j = 1,2,3,...).

\begin{equation}\label{eqn19}
f_{jk}(x,t)\; = \; f_{k}(x,t) \; - \; \sum_{n=1}^{N} u_{j-1,n}\frac{\partial u_{j-1,k}}{\partial x_{n}}\;\;\;\;\;\;({1\leq k \leq N})
\end{equation}

or the vector form

\begin{equation}\label{eqn20}
\vec{f}_{j}(x,t)\; = \; \vec{f}(x,t) \; - \;(\;\vec{u}_{j-1}\;\cdot\;\nabla\;)\;\vec{u}_{j-1}\;
\end{equation}
\\
For the first step of the iterative process (j = 1) we have:  

\[\\(\vec{u}_{0}\;\cdot\;\nabla\;)\;\vec{u}_{0}\;=\;0\]
and
\[\\\vec{f}_{1}(x,t)\; = \; \vec{f}(x,t) \]

\section{Solution. Case N = 2}\ 

We use Fourier transform (\ref{A2}) for equations $(\ref{eqn13})\; - \;(\ref{eqn20})$ and get:

\[\\U_{jk}( \gamma_{1},\gamma_{2},t)\;=\;F [u_{jk} ( x_{1}, x_{2},t)]\]

\[\\ F \sz{[} \frac{\partial^{2}u_{jk}( x_{1}, x_{2},t)}{\partial x^{2}_{s}}\mbox{\Large ]\normalsize}\;=\;-\gamma^{2}_{s}U_{jk}( \gamma_{1},\gamma_{2},t) \;\;\;\;\; \rm{[use (\ref{eqn16})]}\]

\[\\U_{k}^{0}( \gamma_{1},\gamma_{2})\;=\;F[u_{k}^{0} ( x_{1}, x_{2})]\]

\[\\P_{j}( \gamma_{1},\gamma_{2},t)\;=\;F[p_{j}\, ( x_{1}, x_{2},t)]\]

\[\\F_{jk}( \gamma_{1},\gamma_{2},t)\;=\;F[f_{jk} ( x_{1}, x_{2},t)]\]

\[ \\  k,s\;=\;1,2 \]

and then:

\begin{equation}\label{eqn21}
\frac{\partial U_{j1}( \gamma_{1},\gamma_{2},t )}{ \partial t}  \;=\;-\nu
( \gamma_{1}^{2}+\gamma_{2}^{2}) U_{j1}( \gamma_{1},\gamma_{2},t )\;+\;i\gamma_{1} P_{j}( \gamma_{1},\gamma_{2},t )\;+\; F_{j1}( \gamma_{1},\gamma_{2},t )
\end{equation}
\begin{equation}\label{eqn22}
\frac{\partial U_{j2}( \gamma_{1},\gamma_{2},t )}{ \partial t}  \;=\;-\nu
( \gamma_{1}^{2}+\gamma_{2}^{2}) U_{j2}( \gamma_{1},\gamma_{2},t )\;+\;i\gamma_{2} P_{j}( \gamma_{1},\gamma_{2},t )\;+\; F_{j2}( \gamma_{1},\gamma_{2},t )
\end{equation}

\begin{equation}\label{eqn23}
\gamma_{1} U_{j1}( \gamma_{1},\gamma_{2},t )\;+\; \gamma_{2} U_{j2}( \gamma_{1},\gamma_{2},t )\;=\;0
\end{equation}

\begin{equation}\label{eqn24}
U_{j1}(\gamma_{1},\gamma_{2},0)\;=\; U_{1}^{0}(\gamma_{1},\gamma_{2})
\end{equation}

\begin{equation}\label{eqn25}
U_{j2}(\gamma_{1},\gamma_{2},0)\;=\; U_{2}^{0}(\gamma_{1},\gamma_{2})
\end{equation}

Hence eliminate $P_{j}(\gamma_{1},\gamma_{2},t)$ from equations $(\ref{eqn21})$, $(\ref{eqn22})$ and find:

\begin{eqnarray} \label{eqn26}
\frac{\partial }{\partial t} \szp{[}  U_{j2}( \gamma_{1},\gamma_{2},t )\; -\;\frac{\gamma_{2}}{\gamma_{1}} U_{j1}( \gamma_{1},\gamma_{2},t) \szp{]} \;= 
\quad\quad\quad\quad\quad\quad \quad\quad \quad 
\nonumber \\ 
\nonumber \\ 
-\nu( \gamma_{1}^{2}+\gamma_{2}^{2}) \szp{[} U_{j2}( \gamma_{1},\gamma_{2},t )\; -\;\frac{\gamma_{2}}{\gamma_{1}} U_{j1}( \gamma_{1},\gamma_{2},t) \szp{]} \;+ \;\szp{[} F_{j2}( \gamma_{1},\gamma_{2},t )\; -\;\frac{\gamma_{2}}{\gamma_{1}} F_{j1}( \gamma_{1},\gamma_{2},t) \szp{]} 
\end{eqnarray}

\begin{equation}\label{eqn27}
\gamma_{1} U_{j1}( \gamma_{1},\gamma_{2},t )\;+\; \gamma_{2} U_{j2}( \gamma_{1},\gamma_{2},t )\;=\;0
\end{equation}

\begin{equation}\label{eqn28}
U_{j1}(\gamma_{1},\gamma_{2},0)\;=\; U_{1}^{0}(\gamma_{1},\gamma_{2})
\end{equation}

\begin{equation}\label{eqn29}
U_{j2}(\gamma_{1},\gamma_{2},0)\;=\; U_{2}^{0}(\gamma_{1},\gamma_{2})
\end{equation}

We use Laplace transform ($\ref{A4}$), ($\ref{A5}$)  for equations ($\ref{eqn26}$), ($\ref{eqn27}$)  and have:

\[\\U_{jk}^{\otimes} (\gamma_{1},\gamma_{2},\eta) \;=\;L[\,U_{jk}(\gamma_{1},\gamma_{2},t)\,] \;\;\;\;\;\;\;       \rm{k=1,2}\]

\begin{eqnarray}\label{eqn30}
\eta \,\szp{[}  U_{j2}^{\otimes}( \gamma_{1},\gamma_{2},\eta )\; -\;\frac{\gamma_{2}}{\gamma_{1}}\, U_{j1}^{\otimes}( \gamma_{1},\gamma_{2},\eta) \szp{]} \;-\; \szp{[}  U_{j2}( \gamma_{1},\gamma_{2},0 )\; -\;\frac{\gamma_{2}}{\gamma_{1}}\, U_{j1}( \gamma_{1},\gamma_{2},0) \szp{]} \;=  
\nonumber\\
\nonumber\\
-\nu\,( \gamma_{1}^{2}+\gamma_{2}^{2})\szp{[}  U_{j2}^{\otimes}                          ( \gamma_{1},\gamma_{2},\eta )\; -\;\frac{\gamma_{2}}{\gamma_{1}}\, U_{j1}^{\otimes}( \gamma_{1},\gamma_{2},\eta) \szp{]} \;+\; \szp{[}  F_{j2}^{\otimes}( \gamma_{1},\gamma_{2},\eta )\; -\;\frac{\gamma_{2}}{\gamma_{1}} \,F_{j1}^{\otimes}( \gamma_{1},\gamma_{2},\eta) \szp{]} 
\end{eqnarray}

\begin{equation}\label{eqn31}
\gamma_{1} U_{j1}^{\otimes}( \gamma_{1},\gamma_{2},\eta )\;+\; \gamma_{2} U_{j2}^{\otimes} ( \gamma_{1},\gamma_{2},\eta )\;=\;0
\end{equation}

\begin{equation}\label{eqn32}
U_{j1}(\gamma_{1},\gamma_{2},0)\;=\; U_{1}^{0}(\gamma_{1},\gamma_{2})
\end{equation}

\begin{equation}\label{eqn33}
U_{j2}(\gamma_{1},\gamma_{2},0)\;=\; U_{2}^{0}(\gamma_{1},\gamma_{2})
\end{equation}

The solution of the system of equations $(\ref{eqn30})\;-\; (\ref{eqn33})\;$  is: 

\begin{equation}\label{eqn34}
U_{j1}^{\otimes}( \gamma_{1},\gamma_{2},\eta )\;=\;\frac{[\gamma_{2}^{2} F_{j1}^{\otimes}( \gamma_{1},\gamma_{2},\eta)- \gamma_{1}\gamma_{2} F_{j2}^{\otimes}( \gamma_{1},\gamma_{2},\eta)+\gamma_{2}^{2} U_{1}^{0}(\gamma_{1},\gamma_{2}) -\gamma_{1}\gamma_{2} U_{2}^{0}(\gamma_{1},\gamma_{2})]}{ (\gamma_{1}^{2}+\gamma_{2}^{2}) [\eta+\nu (\gamma_{1}^{2}+\gamma_{2}^{2})]} 
\end{equation}

\begin{equation}\label{eqn35}
U_{j2}^{\otimes}( \gamma_{1},\gamma_{2},\eta )\;=\;\frac{[\gamma_{1}^{2} F_{j2}^{\otimes}( \gamma_{1},\gamma_{2},\eta)- \gamma_{1}\gamma_{2} F_{j1}^{\otimes}( \gamma_{1},\gamma_{2},\eta)+\gamma_{1}^{2} U_{2}^{0}(\gamma_{1},\gamma_{2}) -\gamma_{1}\gamma_{2} U_{1}^{0}(\gamma_{1},\gamma_{2})]}{ (\gamma_{1}^{2}+\gamma_{2}^{2}) [\eta+\nu (\gamma_{1}^{2}+\gamma_{2}^{2})]} 
\end{equation}

Then we use the convolution formula $(\ref{A6})$ and integral $(\ref{A7})$ for $(\ref{eqn34}) \;,\; (\ref{eqn35})\;$  and obtain:

\begin{eqnarray}\label{eqn36}
U_{j1}(\gamma_{1},\gamma_{2},t)\;=\; \int_{0}^{t} \sz{e} ^{-\nu (\gamma_{1}^{2}+\gamma_{2}^{2}) (t-\tau)} \frac{[\gamma_{2}^{2} F_{j1}( \gamma_{1}, \gamma_{2},\tau)-\gamma_{1}\gamma_{2} F_{j2}( \gamma_{1}, \gamma_{2},\tau )]}{ (\gamma_{1}^{2}+\gamma_{2}^{2}) }\,d\tau \;+ 
\nonumber\\
\nonumber\\
+\;\sz{e}  ^{-\nu (\gamma_{1}^{2}+\gamma_{2}^{2}) t} \;U_{1}^{0}(\gamma_{1},\gamma_{2})
\quad\quad\quad\quad\quad\quad \quad\quad \quad\quad 
\end{eqnarray}

\begin{eqnarray}\label{eqn37}
U_{j2}(\gamma_{1},\gamma_{2},t)\;=\;\int_{0}^{t} \sz{e}  ^{-\nu (\gamma_{1}^{2}+\gamma_{2}^{2}) (t-\tau)} \frac{ [\gamma_{1}^{2} F_{j2}( \gamma_{1}, \gamma_{2},\tau)-\gamma_{1}\gamma_{2} F_{j1}( \gamma_{1}, \gamma_{2},\tau )] }{ (\gamma_{1}^{2}+\gamma_{2}^{2}) }\,d\tau\;+
\nonumber\\
\nonumber\\
+\; \sz{e}  ^{-\nu (\gamma_{1}^{2}+\gamma_{2}^{2}) t} \;U_{2}^{0}(\gamma_{1},\gamma_{2})
\quad\quad\quad\quad\quad\quad \quad\quad \quad\quad 
\end{eqnarray}
\\
$ P_{j}(\gamma_{1},\gamma_{2},t) $ is obtained from $(\ref{eqn21})\;$or$\; (\ref{eqn22})\;$:
\begin{equation}\label{eqn38}
P_{j}(\gamma_{1},\gamma_{2},t) \;=\;i\frac{[\gamma_{1} F_{j1}( \gamma_{1}, \gamma_{2},t)+\gamma_{2} F_{j2}( \gamma_{1}, \gamma_{2},t )]}{ (\gamma_{1}^{2}+\gamma_{2}^{2}) }
\end{equation}
\\
Use of the Fourier inversion formula $(\ref{A2})$ and find:

\begin{eqnarray}\label{eqn39}
u_{j1}(x_{1},x_{2},t)\;=\;\frac{1}{2\pi} \int_{-\infty}^{\infty}\int_{-\infty}^{\infty} \biggl[  \int_{0}^{t} \sz{e}  ^{-\nu (\gamma_{1}^{2}+\gamma_{2}^{2}) (t-\tau)} \frac{[\gamma_{2}^{2} F_{j1}( \gamma_{1}, \gamma_{2},\tau)-\gamma_{1}\gamma_{2} F_{j2}( \gamma_{1}, \gamma_{2},\tau )]}{ (\gamma_{1}^{2}+\gamma_{2}^{2}) } \,d\tau\;+
\nonumber\\
\nonumber\\
\nonumber\\
+\; \sz{e}  ^{-\nu (\gamma_{1}^{2}+\gamma_{2}^{2}) t} \;U_{1}^{0}(\gamma_{1},\gamma_{2})\biggr] \;\sz{e}  ^{-i(x_{1}\gamma_{1}+x_{2}\gamma_{2})}\,d\gamma_{1}d\gamma_{2}\;=
\quad\quad\quad\quad\quad\quad \quad\quad \quad\quad 
\nonumber\\
\nonumber\\
\nonumber\\
= \;\frac{1}{4\pi^2} \int_{-\infty}^{\infty}\int_{-\infty}^{\infty}\frac{\gamma_{2}^{2} }{(\gamma_{1}^{2}+\gamma_{2}^{2}) } \int_{0}^{t} \sz{e}  ^{-\nu (\gamma_{1}^{2}+\gamma_{2}^{2}) (t-\tau)} \int_{-\infty}^{\infty}\int_{-\infty}^{\infty}\sz{e}  ^{i(\tilde x_{1}\gamma_{1}+\tilde x_{2}\gamma_{2})} f_{j1}(\tilde x_{1},\tilde x_{2},\tau)\,d\tilde x_{1}d\tilde x_{2}d\tau\cdot
\nonumber\\
\nonumber\\
\nonumber\\
 \;\cdot\;\sz{e}  ^{-i(x_{1}\gamma_{1}+x_{2}\gamma_{2})}\,d\gamma_{1}d\gamma_{2}\;-
\quad\quad\quad\quad\quad\quad \quad\quad \quad\quad
\nonumber\\
\nonumber\\
\nonumber\\
- \;\frac{1}{4\pi^2} \int_{-\infty}^{\infty}\int_{-\infty}^{\infty}\frac{\gamma_{1} \gamma_{2} }{(\gamma_{1}^{2}+\gamma_{2}^{2}) } \int_{0}^{t} \sz{e}  ^{-\nu (\gamma_{1}^{2}+\gamma_{2}^{2}) (t-\tau)} \int_{-\infty}^{\infty}\int_{-\infty}^{\infty}\sz{e}  ^{i(\tilde x_{1}\gamma_{1}+\tilde x_{2}\gamma_{2})} f_{j2}(\tilde x_{1},\tilde x_{2},\tau)\,d\tilde x_{1}d\tilde x_{2}d\tau\cdot
\nonumber\\
\nonumber\\
\nonumber\\
 \;\cdot\;\sz{e}  ^{-i(x_{1}\gamma_{1}+x_{2}\gamma_{2})}\,d\gamma_{1}d\gamma_{2}\;+
\quad\quad\quad\quad\quad\quad \quad\quad \quad\quad
\nonumber\\
\nonumber\\
\nonumber\\
+ \;\frac{1}{4\pi^2} \int_{-\infty}^{\infty}\int_{-\infty}^{\infty}\sz{e}  ^{-\nu (\gamma_{1}^{2}+\gamma_{2}^{2}) t} \int_{-\infty}^{\infty}\int_{-\infty}^{\infty}\sz{e}  ^{i(\tilde x_{1}\gamma_{1}+\tilde x_{2}\gamma_{2})} u^{0}_{1}(\tilde x_{1},\tilde x_{2})\,d\tilde x_{1}d\tilde x_{2}\sz{e}  ^{-i(x_{1}\gamma_{1}+x_{2}\gamma_{2})}\,d\gamma_{1}d\gamma_{2}\;=
\nonumber\\
\nonumber\\
\nonumber\\
= \;S_{11}(f_{j1})\;+\; S_{12}(f_{j2})\;+\;B(u^{0}_{1})
\quad\quad\quad\quad\quad\quad \quad\quad \quad\quad 
 \end{eqnarray}

\begin{eqnarray}\label{eqn40}
u_{j2} (x_{1},x_{2},t)\;=\;\frac{1}{2\pi} \int_{-\infty}^{\infty}\int_{-\infty}^{\infty} \biggl[ \int_{0}^{t} \sz{e} ^{-\nu (\gamma_{1}^{2}+\gamma_{2}^{2}) (t-\tau)} \frac{[\gamma_{1}^{2} F_{j2}( \gamma_{1}, \gamma_{2},\tau)-\gamma_{1}\gamma_{2} F_{j1}( \gamma_{1}, \gamma_{2},\tau )]}{ (\gamma_{1}^{2}+\gamma_{2}^{2}) }\,d\tau\;+ \nonumber\\
\nonumber\\
\nonumber\\
+\; \sz{e} ^{-\nu (\gamma_{1}^{2}+\gamma_{2}^{2}) t} \;U_{2}^{0}(\gamma_{1},\gamma_{2})\biggr]  \;\sz{e} ^{-i(x_{1}\gamma_{1}+x_{2}\gamma_{2})}\,d\gamma_{1}d\gamma_{2}\;=
\quad\quad\quad\quad\quad\quad \quad\quad \quad\quad 
\nonumber\\
\nonumber\\
\nonumber\\
= \;-\;\frac{1}{4\pi^2} \int_{-\infty}^{\infty}\int_{-\infty}^{\infty}\frac{\gamma_{1} \gamma_{2} }{(\gamma_{1}^{2}+\gamma_{2}^{2}) } \int_{0}^{t} \sz{e}  ^{-\nu (\gamma_{1}^{2}+\gamma_{2}^{2}) (t-\tau)} \int_{-\infty}^{\infty}\int_{-\infty}^{\infty}\sz{e}  ^{i(\tilde x_{1}\gamma_{1}+\tilde x_{2}\gamma_{2})} f_{j1}(\tilde x_{1},\tilde x_{2},\tau)\,d\tilde x_{1}d\tilde x_{2}d\tau\cdot
\nonumber\\
\nonumber\\
\nonumber\\
 \;\cdot\;\sz{e}  ^{-i(x_{1}\gamma_{1}+x_{2}\gamma_{2})}\,d\gamma_{1}d\gamma_{2}\;+
\quad\quad\quad\quad\quad\quad \quad\quad \quad\quad
\nonumber\\
\nonumber\\
\nonumber\\
+ \;\frac{1}{4\pi^2} \int_{-\infty}^{\infty}\int_{-\infty}^{\infty}\frac{\gamma_{1}^{2}  }{(\gamma_{1}^{2}+\gamma_{2}^{2}) } \int_{0}^{t} \sz{e}  ^{-\nu (\gamma_{1}^{2}+\gamma_{2}^{2}) (t-\tau)} \int_{-\infty}^{\infty}\int_{-\infty}^{\infty}\sz{e}  ^{i(\tilde x_{1}\gamma_{1}+\tilde x_{2}\gamma_{2})} f_{j2}(\tilde x_{1},\tilde x_{2},\tau)\,d\tilde x_{1}d\tilde x_{2}d\tau\cdot
\nonumber\\
\nonumber\\
\nonumber\\
 \;\cdot\;\sz{e}  ^{-i(x_{1}\gamma_{1}+x_{2}\gamma_{2})}\,d\gamma_{1}d\gamma_{2}\;+
\quad\quad\quad\quad\quad\quad \quad\quad \quad\quad
\nonumber\\
\nonumber\\
\nonumber\\
+ \;\frac{1}{4\pi^2} \int_{-\infty}^{\infty}\int_{-\infty}^{\infty}\sz{e}  ^{-\nu (\gamma_{1}^{2}+\gamma_{2}^{2}) t}  \int_{-\infty}^{\infty}\int_{-\infty}^{\infty}\sz{e}  ^{i(\tilde x_{1}\gamma_{1}+\tilde x_{2}\gamma_{2})} u^{0}_{2}(\tilde x_{1},\tilde x_{2})\,d\tilde x_{1}d\tilde x_{2}\sz{e}  ^{-i(x_{1}\gamma_{1}+x_{2}\gamma_{2})}\,d\gamma_{1}d\gamma_{2}\;=
\nonumber\\
\nonumber\\
\nonumber\\
= \;S_{21}(f_{j1})\;+\; S_{22}(f_{j2})\;+\;B(u^{0}_{2})
\quad\quad\quad\quad\quad\quad \quad\quad \quad\quad 
\end{eqnarray}

\begin{eqnarray}\label{eqn41}
p_{j}\,(x_{1},x_{2},t)\;=\; \frac{i}{2\pi} \int_{-\infty}^{\infty}\int_{-\infty}^{\infty}\frac{[\gamma_{1} F_{j1}( \gamma_{1}, \gamma_{2},t)+ \gamma_{2} F_{j2}( \gamma_{1}, \gamma_{2},t )]}{ (\gamma_{1}^{2}+\gamma_{2}^{2}) }\;\sz{e} ^{-i(x_{1}\gamma_{1}+x_{2}\gamma_{2})}\,d\gamma_{1}d\gamma_{2}\;=\;
\nonumber\\
\nonumber\\
\nonumber\\
= \;-\;\frac{i}{4\pi^2} \int_{-\infty}^{\infty}\int_{-\infty}^{\infty}\frac{\gamma_{1}}{(\gamma_{1}^{2}+\gamma_{2}^{2}) } \int_{-\infty}^{\infty}\int_{-\infty}^{\infty}\sz{e}  ^{i(\tilde x_{1}\gamma_{1}+\tilde x_{2}\gamma_{2})} f_{j1}(\tilde x_{1},\tilde x_{2},t)\,d\tilde x_{1}d\tilde x_{2}\;
\sz{e}  ^{-i(x_{1}\gamma_{1}+x_{2}\gamma_{2})}\,d\gamma_{1}d\gamma_{2}\;+
\nonumber\\
\nonumber\\
\nonumber\\
+ \;\frac{i}{4\pi^2} \int_{-\infty}^{\infty}\int_{-\infty}^{\infty}\frac{\gamma_{2}}{(\gamma_{1}^{2}+\gamma_{2}^{2}) }  \int_{-\infty}^{\infty}\int_{-\infty}^{\infty}\sz{e}  ^{i(\tilde x_{1}\gamma_{1}+\tilde x_{2}\gamma_{2})} f_{j2}(\tilde x_{1},\tilde x_{2},t)\,d\tilde x_{1}d\tilde x_{2}\;
\sz{e}  ^{-i(x_{1}\gamma_{1}+x_{2}\gamma_{2})}\,d\gamma_{1}d\gamma_{2}\;=
\nonumber\\
\nonumber\\
=\;\tilde S_{1}(f_{j1})\;+\; \tilde S_{2}(f_{j2})
\quad\quad\quad\quad\quad\quad \quad\quad \quad\quad\quad\quad\quad\quad\quad\quad \quad\quad  
\end{eqnarray}

So, the integrals $(\ref{eqn39})\;-\; (\ref{eqn41})\;$exist by the restrictions $(\ref{eqn17})\;, (\ref{eqn18})\;$.

Here $S_{11}(), S_{12}(), S_{21}(), S_{22}(), B(), \tilde S_{1}(), \tilde S_{2}()$ are the integral - operators.

\[S_{12}()\;= \;S_{21}() \]

We have for the vector $\vec{u}_{j}$ from the equations $(\ref{eqn39})\;-\; (\ref{eqn40})\;$:

\begin{equation}\label{eqn42}
\vec{u}_{j}\;=\;\bar{\bar{S}}\;\cdot\;\vec{f}_{j}\;+\;B(\vec{u}^{0})\;,
\end{equation}

where $\;\bar{\bar{S}} \; $ is the matrix - operator:
\[ \left( \begin{array}{ccc}
S_{11} & S_{12} \\
S_{21} & S_{22} \end{array} \right)\]

We put $\vec{f}_{j}$ from equation $(\ref{eqn20})$ into equation $(\ref{eqn42})$ and have:

\begin{eqnarray}\label{eqn44}
\vec{u}_{j} = \bar{\bar{S}}\cdot(\;\vec{f}\;-\;(\;\vec{u}_{j-1}\cdot\nabla)\vec{u}_{j-1})\;+\;B(\vec{u}^{0})\;=
\nonumber\\
\nonumber\\
=\;\bar{\bar{S}}\cdot\vec{f} \;-\;\bar{\bar{S}}\cdot(\vec{u}_{j-1}\;\cdot\;\nabla\;)\;\vec{u}_{j-1}\;+\;B(\vec{u}^{0})\; =
\nonumber\\
\nonumber\\
=\;\vec{u}_{1}\;-\;\bar{\bar{S}}\cdot(\vec{u}_{j-1}\;\cdot\;\nabla)\;\vec{u}_{j-1}
\quad\quad\quad\quad\quad\quad \quad\quad 
\end{eqnarray}

Here $\vec{u}_{1}\;$ is the solution of the system of equations $(\ref{eqn13})\; - \;(\ref{eqn20})$ with condition:

\[\sum_{n=1}^{2} u_{n}\frac{\partial u_{k}}{\partial x_{n}}\;=\;0\;\;\;\;\;\;\;       \rm{k=1,2}\;\;\;  \] 

For j = 1 formula $(\ref{eqn42})$ can be written as follows:

\begin{equation}\label{eqn47}
\vec{u}_{1}\;=\;\bar{\bar{S}}\;\cdot\;\vec{f}_{1}\;+\;B(\vec{u}^{0})\;,\;\;\;\;\;\\\vec{f}_{1}(x,t)\; = \; \vec{f}(x,t)
\end{equation}

If t $\rightarrow$ 0 then $\vec{u}_{1} \rightarrow \vec{u}^{0}$ (look at integral-operators $\bar{\bar{S}}, B()\;\;$ -   integrals $\;(\ref{eqn39})\; , \;(\ref{eqn40})$).

For j = 2  we define from equation $(\ref{eqn20})$:

\begin{equation}
\vec{f}_{2}(x,t)\; = \; \vec{f}_{1}(x,t) \; - \;(\;\vec{u}_{1}\;\cdot\;\nabla\;)\;\vec{u}_{1}\;
\end{equation}

We denote: 

\begin{equation}\label{eqn48}
\vec{f}_{2}^{*}\;=\;(\vec{u}_{1}\;\cdot\;\nabla)\;\vec{u}_{1}
\end{equation}

and then we have: 

\begin{equation}
\vec{f}_{2}(x,t)\; = \; \vec{f}_{1}(x,t) \; - \vec{f}_{2}^{*}
\end{equation}

Then we get $\vec{u}_{2}$ from $(\ref{eqn42}),(\ref{eqn47})$:

\begin{equation}\label{eqn50}
\vec{u}_{2}\;=\;\bar{\bar{S}}\;\cdot\;\vec{f}_{2}  \;+\;B(\vec{u}^{0})\;=\;\bar{\bar{S}}\;\cdot\;(\vec{f}_{1}\;-\;\vec{f}_{2}^{*})  \;+\;B(\vec{u}^{0})\;=\;\vec{u}_{1}\;-\;\vec{u}_{2}^{*}
\end{equation}

Here we have:

\begin{equation}\label{eqn49}
\vec{u}_{2}^{*}\;=\;\bar{\bar{S}}\;\cdot\;\vec{f}_{2}^{*}
\end{equation}

If t $\rightarrow$ 0 then $\vec{u}_{2}^{*} \rightarrow$ 0 (look at integral-operator $\bar{\bar{S}}\;\;$- integrals $\;(\ref{eqn39})\; , \;(\ref{eqn40})$).

Continue for j = 3. We define from equation $(\ref{eqn20})$:

\begin{equation}
\vec{f}_{3}(x,t)\; = \; \vec{f}_{1}(x,t) \; - \;(\;\vec{u}_{2}\;\cdot\;\nabla\;)\;\vec{u}_{2}\;
\end{equation}

Here we have:

\begin{equation}\label{eqn51}
(\vec{u}_{2}\;\cdot\;\nabla)\;\vec{u}_{2}\;=
\;((\vec{u}_{1}\;-\;\vec{u}_{2}^{*})\;\cdot\;\nabla\;)
\;(\vec{u}_{1}\;-\;\vec{u}_{2}^{*})\;=
\;\vec{f}_{2}^{*}\;+\;\vec{f}_{3}^{*}
\end{equation}

We denote in $(\ref{eqn51})$:

\[\vec{f}_{3}^{*}\;=\;-\;  (\vec{u}_{1}\;\cdot\;\nabla)\;\vec{u}_{2}^{*}\; -\; (\vec{u}_{2}^{*}\;\cdot\;\nabla)\;\vec{u}_{1}\;+\;
(\vec{u}_{2}^{*}\;\cdot\;\nabla)\;\vec{u}_{2}^{*}\]

and then we have: 

\begin{equation}
\vec{f}_{3}(x,t)\; = \; \vec{f}_{1}(x,t) \; - \vec{f}_{2}^{*}\; - \vec{f}_{3}^{*}
\end{equation}

Then we get $\vec{u}_{3}$ from $(\ref{eqn42}) , (\ref{eqn47}) , (\ref{eqn49})$:

\begin{equation}\label{eqn53}
\vec{u}_{3}\;=\;\bar{\bar{S}}\;\cdot\;\vec{f}_{3}  \;+\;B(\vec{u}^{0})\;=\;\bar{\bar{S}}\;\cdot\;(\vec{f}_{1} \;-\;\vec{f}_{2}^{*}\;-\;\vec{f}_{3}^{*})\;+\;B(\vec{u}^{0})\;=\;\vec{u}_{1} \;-\;\vec{u}_{2}^{*}\;-\;\vec{u}_{3}^{*}
\end{equation}

Here we denote:

\begin{equation}\label{eqn52}
\vec{u}_{3}^{*}\;=\;\bar{\bar{S}}\;\cdot\;\vec{f}_{3}^{*}
\end{equation}

If t $\rightarrow$ 0 then $\vec{u}_{3}^{*} \rightarrow$ 0 (look at integral-operator $\bar{\bar{S}}\;\;$- integrals $\;(\ref{eqn39})\; , \;(\ref{eqn40})$).

For j = 4. We define from equation $(\ref{eqn20})$:

\begin{equation}
\vec{f}_{4}(x,t)\; = \; \vec{f}_{1}(x,t) \; - \;(\;\vec{u}_{3}\;\cdot\;\nabla\;)\;\vec{u}_{3}\;
\end{equation}

Here we have:

\begin{equation}\label{eqn54}
(\vec{u}_{3}\;\cdot\;\nabla)\;\vec{u}_{3}\;=
\;((\vec{u}_{2}\;-\;\vec{u}_{3}^{*})\;\cdot\;\nabla\;)
\;(\vec{u}_{2}\;-\;\vec{u}_{3}^{*})\;=
\;\vec{f}_{2}^{*}\;+\;\vec{f}_{3}^{*}\;+\;\vec{f}_{4}^{*}
\end{equation}
\\
We denote in $(\ref{eqn54})$:

\[\vec{f}_{4}^{*}\;=\;-\;  (\vec{u}_{2}\;\cdot\;\nabla)\;\vec{u}_{3}^{*}\; -\; (\vec{u}_{3}^{*}\;\cdot\;\nabla)\;\vec{u}_{2}\;+\;
(\vec{u}_{3}^{*}\;\cdot\;\nabla)\;\vec{u}_{3}^{*}\]

and then we have: 

\begin{equation}
\vec{f}_{4}(x,t)\; = \; \vec{f}_{1}(x,t) \; - \vec{f}_{2}^{*}\; - \vec{f}_{3}^{*}\; - \vec{f}_{4}^{*}
\end{equation}

Then we get $\vec{u}_{4}$ from $(\ref{eqn42}) , (\ref{eqn47}) , (\ref{eqn49}) , (\ref{eqn52})$:

\begin{equation}\label{eqn56}
\vec{u}_{4}\;=\;\bar{\bar{S}}\;\cdot\;\vec{f}_{4}\;+\;B(\vec{u}^{0})\;=\;\bar{\bar{S}}\;\cdot\;(\vec{f}_{1} \; -\;\vec{f}_{2}^{*}\; -\;\vec{f}_{3}^{*}\;-\;\vec{f}_{4}^{*})\;+\;B(\vec{u}^{0})\;=\;\vec{u}_{1} \;-\;\vec{u}_{2}^{*}\;-\;\vec{u}_{3}^{*}\;-\;\vec{u}_{4}^{*}
\end{equation}

Here we denote:

\begin{equation}\label{eqn55}
\vec{u}_{4}^{*}\;=\;\bar{\bar{S}}\;\cdot\;\vec{f}_{4}^{*}
\end{equation}

If t $\rightarrow$ 0 then $\vec{u}_{4}^{*} \rightarrow 0$ (look at integral-operator $\bar{\bar{S}}\;\;$- integrals $\;(\ref{eqn39})\; , \;(\ref{eqn40})$).

For arbitrary number j $(j \geq 2)$. We define from equation $(\ref{eqn20})$:

\begin{equation}
\vec{f}_{j}(x,t)\; = \; \vec{f}_{1}(x,t) \; - \;(\;\vec{u}_{j-1}\;\cdot\;\nabla\;)\;\vec{u}_{j-1}\;
\end{equation}

Here we have:

\begin{equation}\label{eqn60}
(\vec{u}_{j-1}\;\cdot\;\nabla)\;\vec{u}_{j-1}\;=\;\sum_{l=2}^{j} \vec{f}_{l}^{*}
\end{equation}

and it follows:

\begin{equation}\label{eqn60a}
\vec{f}_{j}\;=\;\vec{f}_{1}\;-\; \sum_{l=2}^{j} \vec{f}_{l}^{*}
\end{equation}

Then we get $\vec{u}_{j}$ from $(\ref{eqn42}) , (\ref{eqn47}) $:

\begin{equation}\label{eqn57}
\vec{u}_{j}\;=\;\bar{\bar{S}}\;\cdot\; \vec{f}_{j}\;+ \;B(\vec{u}^{0})\;=\;\bar{\bar{S}}\;\cdot\;( \vec{f}_{1}\;-\; \sum_{l=2}^{j} \vec{f}_{l}^{*})\;+ \;B(\vec{u}^{0})\;=\;\vec{u}_{1}\;-\;\sum_{l=2}^{j} \vec{u}_{l}^{*}\
\end{equation}

Here we denote:

\begin{equation}\label{eqn59}
\vec{u}_{l}^{*}\;=\;\bar{\bar{S}}\;\cdot\;\vec{f}_{l}^{*}\;\;\;\;\;\;\;\;\;\;\;\;\;\;(2\;\leq\;l\;\leq\;j)
\end{equation}

If t $\rightarrow$ 0 then $\vec{u}_{l}^{*} \rightarrow$ 0 (look at integral-operator $\bar{\bar{S}}\;\;$- integrals $\;(\ref{eqn39})\; , \;(\ref{eqn40})$).

We consider the equations $(\ref{eqn47})$ - $(\ref{eqn59})$ and see that the series $(\ref{eqn57})$ converge for $j \rightarrow \infty$ 

with the conditions for the first step (j = 1) of the iterative process:

\[\;\;\;\sum_{n=1}^{2} u_{0n}\frac{\partial u_{0k}}{\partial x_{n}} = 0\;\;\;\;\;\;\;\;       \rm{k=1,2}\]

and conditions

\begin{equation}\label{eqn61}
C_{\alpha K}\leq\;\frac{1}{2}\;\;,\;\;C_{\alpha \beta K}\leq\;\frac{1}{2}
\end{equation}

 Here $\;\;C_{\alpha K}\;\;$and$\;\;C_{\alpha \beta K}\;$ are received from $(\ref{eqn17})$ , $(\ref{eqn18}).$

Hence, we receive from equation $(\ref{eqn44})\;$ when $j \rightarrow \infty$:

\begin{equation}\label{eqn45}
\vec{u}_{\infty}=\;\vec{u}_{1}\;-\;\bar{\bar{S}}\cdot(\vec{u}_{\infty}\;\cdot\;\nabla)\;\vec{u}_{\infty}
\end{equation}

Equation $(\ref{eqn45})$ describes the converging iterative process. 

Then we have from formula $(\ref{eqn41})\;$:

\begin{equation}\label{eqn46a}
p_{\infty}\,\;=\; \;\tilde S_{1}(f_{\infty 1})\;+\; \tilde S_{2}(f_{\infty 2})
\end{equation}

Here $\vec{f}_{\infty}$ = ($f_{\infty 1} , f_{\infty 2}$) is received from formula $(\ref{eqn60a})\;$.

On the other hand we can transform the original system of differential equations $(\ref{eqn7})\; - \;(\ref{eqn9})$ to the equivalent system of integral equations by the scheme of iterative process $(\ref{eqn42})\;, \;(\ref{eqn44})$ for vector $\vec{u}$:

\begin{equation}\label{eqn46}
\vec{u}\;=\;\vec{u}_{1}\;-\;\bar{\bar{S}}\cdot(\vec{u}\;\cdot\;\nabla)\;\vec{u},
\end{equation}

where $\vec{u}_{1}$ is from formula $(\ref{eqn47})$.
We compare the equations $(\ref{eqn45})$ and $(\ref{eqn46})$ and see that the iterative process $(\ref{eqn45})$ converge to the solution of the system $(\ref{eqn46})$ and hence to the solution of the differential equations $(\ref{eqn7})\; - \;(\ref{eqn9})$ with conditions $(\ref{eqn61})$. 

\textbf{In other words there exist smooth functions}  $\mathbf{p_{\infty}(x, t)}$, $\mathbf{u_{\infty i}(x, t)}$  \textbf{(i = 1, 2) on} $\mathbf{R^{2} \times [0,\infty)}$ \textbf{that satisfy} $\mathbf{(\ref{eqn1}), (\ref{eqn2}), (\ref{eqn3})}$ \textbf{and}

\begin{equation}\label{eqn46b}
\mathbf{p_{\infty}, \;u_{\infty i} \in  C^{\infty}(R^{2} \times [0,\infty)),}
\nonumber\\
\nonumber\\
\end{equation}

\begin{equation}\label{eqn4aa}
\mathbf{\int_{R^{2}}|\vec{u}_{\infty}(x, t)|^{2}dx < C } 
\end{equation}

\textbf{for all t} $\mathbf{\geq 0}$.

\section{Solution. Case N = 3}\ 

We use Fourier transform $(\ref{A3})$ for equations $(\ref{eqn13})\; - \;(\ref{eqn20})$ and get:

\[\\U_{jk}( \gamma_{1}, \gamma_{2}, \gamma_{3}, t)\;=\;F[u_{jk} ( x_{1}, x_{2}, x_{3}, t)]\]

\[F\sz{[} \frac{\partial^{2}u_{jk}( x_{1}, x_{2}, x_{3}, t)}{\partial x^{2}_{s}}\sz{]} \;=\;-\gamma^{2}_{s}U_{jk}( \gamma_{1}, \gamma_{2}, \gamma_{3}, t) \;\;\;\;\; \rm{[use (\ref{eqn16})]}\]

\[\\U_{k}^{0}( \gamma_{1} ,\gamma_{2} ,\gamma_{3})\;=\;F[u_{k}^{0} ( x_{1}, x_{2}, x_{3})]\]

\[\\P_{j}( \gamma_{1}, \gamma_{2}, \gamma_{3}, t)\;=\;F[p_{j}\, ( x_{1}, x_{2}, x_{3}, t)]\]

\[\\F_{jk}( \gamma_{1}, \gamma_{2}, \gamma_{3}, t)\;=\;F[f_{jk} ( x_{1}, x_{2}, x_{3}, t)]\]

\[ \\  k,s\;=\;1,2,3 \]

and then:

\begin{equation}\label{eqn134}
\frac{d U_{j1}( \gamma_{1}, \gamma_{2}, \gamma_{3}, t )}{d t}  \;=\;-\nu
( \gamma_{1}^{2} +\gamma_{2}^{2} +\gamma_{3}^{2}) U_{j1}( \gamma_{1}, \gamma_{2}, \gamma_{3}, t )\;+\;i\gamma_{1} P_{j}( \gamma_{1}, \gamma_{2}, \gamma_{3}, t )\;+\; F_{j1}( \gamma_{1}, \gamma_{2}, \gamma_{3},t )
\end{equation}
\begin{equation}\label{eqn135}
\frac{d U_{j2}( \gamma_{1}, \gamma_{2}, \gamma_{3}, t )}{d t}  \;=\;-\nu
( \gamma_{1}^{2} +\gamma_{2}^{2} +\gamma_{3}^{2}) U_{j2}( \gamma_{1}, \gamma_{2}, \gamma_{3}, t )\;+\;i\gamma_{2} P_{j}( \gamma_{1}, \gamma_{2}, \gamma_{3}, t )\;+\; F_{j2}( \gamma_{1}, \gamma_{2}, \gamma_{3},t )
\end{equation}
\begin{equation}\label{eqn136}
\frac{d U_{j3}( \gamma_{1}, \gamma_{2}, \gamma_{3}, t )}{d t}  \;=\;-\nu
( \gamma_{1}^{2} +\gamma_{2}^{2} +\gamma_{3}^{2}) U_{j3}( \gamma_{1}, \gamma_{2}, \gamma_{3}, t )\;+\;i\gamma_{3} P_{j}( \gamma_{1}, \gamma_{2}, \gamma_{3}, t )\;+\; F_{j3}( \gamma_{1}, \gamma_{2}, \gamma_{3},t )
\end{equation}

\begin{equation}\label{eqn137}
\gamma_{1} U_{j1}( \gamma_{1}, \gamma_{2}, \gamma_{3}, t ) \;+\; \gamma_{2}\, U_{j2}( \gamma_{1}, \gamma_{2}, \gamma_{3}, t ) \;+\; \gamma_{3}\, U_{j3}( \gamma_{1}, \gamma_{2}, \gamma_{3}, t ) \;=\;0
\end{equation}

\begin{equation}\label{eqn138}
U_{j1}(\gamma_{1}, \gamma_{2},  \gamma_{3},  0)\;=\; U_{1}^{0}(\gamma_{1} ,\gamma_{2} ,\gamma_{3})
\end{equation}

\begin{equation}\label{eqn139}
U_{j2}(\gamma_{1}, \gamma_{2},  \gamma_{3},  0)\;=\; U_{2}^{0}(\gamma_{1} ,\gamma_{2} ,\gamma_{3})
\end{equation}

\begin{equation}\label{eqn140}
U_{j3}(\gamma_{1}, \gamma_{2},  \gamma_{3},  0)\;=\; U_{3}^{0}(\gamma_{1} ,\gamma_{2} ,\gamma_{3})
\end{equation}

Hence eliminate $P_{j}(\gamma_{1}, \gamma_{2}, \gamma_{3}, t)$ from equations $(\ref{eqn134})\;-\; (\ref{eqn136})\;$ and find:

\begin{eqnarray}\label{eqn141}
\frac{d}{dt} \szp{[} U_{j2}( \gamma_{1}, \gamma_{2}, \gamma_{3}, t )\; -\;\frac{\gamma_{2}}{\gamma_{1}} \,U_{j1}( \gamma_{1}, \gamma_{2}, \gamma_{3}, t) \szp{]} \;= -\nu( \gamma_{1}^{2} +\gamma_{2}^{2} +\gamma_{3}^{2}) \szp{[} U_{j2}( \gamma_{1}, \gamma_{2}, \gamma_{3}, t )\; -\; \nonumber \\ 
\nonumber\\
-\;\frac{\gamma_{2}}{\gamma_{1}}\, U_{j1}( \gamma_{1}, \gamma_{2}, \gamma_{3}, t) \szp{]}   + \; \szp{[} F_{j2}( \gamma_{1}, \gamma_{2}, \gamma_{3}, t )\; -\;\frac{\gamma_{2}}{\gamma_{1}} \, F_{j1}( \gamma_{1}, \gamma_{2}, \gamma_{3}, t) \szp{]} 
\quad\quad\quad\quad\quad\quad 
\end{eqnarray}

\begin{eqnarray}\label{eqn142}
\frac{d}{dt} \szp{[}  U_{j3}( \gamma_{1}, \gamma_{2}, \gamma_{3}, t )\; -\;\frac{\gamma_{3}}{\gamma_{1}}\, U_{j1}( \gamma_{1}, \gamma_{2}, \gamma_{3}, t) \szp{]} \;= -\nu( \gamma_{1}^{2} +\gamma_{2}^{2} +\gamma_{3}^{2}) \szp{[} U_{j3}( \gamma_{1}, \gamma_{2}, \gamma_{3}, t )\; -\; \nonumber \\ 
\nonumber\\
-\;\frac{\gamma_{3}}{\gamma_{1}}\, U_{j1}( \gamma_{1}, \gamma_{2}, \gamma_{3}, t) \szp{]}   + \;\szp{[} F_{j3}( \gamma_{1}, \gamma_{2}, \gamma_{3}, t )\; -\;\frac{\gamma_{3}}{\gamma_{1}}\, F_{j1}( \gamma_{1}, \gamma_{2}, \gamma_{3}, t) \szp{]} 
\quad\quad\quad\quad\quad\quad 
\end{eqnarray}

\begin{equation}\label{eqn143}
\gamma_{1} U_{j1}( \gamma_{1}, \gamma_{2}, \gamma_{3}, t ) \;+\; \gamma_{2}\, U_{j2}( \gamma_{1}, \gamma_{2}, \gamma_{3}, t ) \;+\; \gamma_{3}\, U_{j3}( \gamma_{1}, \gamma_{2}, \gamma_{3}, t ) \;=\;0
\end{equation}

\begin{equation}\label{eqn144}
U_{j1}(\gamma_{1}, \gamma_{2},  \gamma_{3},  0)\;=\; U_{1}^{0}(\gamma_{1} ,\gamma_{2} ,\gamma_{3})
\end{equation}

\begin{equation}\label{eqn145}
U_{j2}(\gamma_{1}, \gamma_{2},  \gamma_{3},  0)\;=\; U_{2}^{0}(\gamma_{1} ,\gamma_{2} ,\gamma_{3})
\end{equation}

\begin{equation}\label{eqn146}
U_{j3}(\gamma_{1}, \gamma_{2},  \gamma_{3},  0)\;=\; U_{3}^{0}(\gamma_{1} ,\gamma_{2} ,\gamma_{3})
\end{equation}

We use Laplace transform $(\ref{A4}), (\ref{A5})$  for equations $(\ref{eqn141})\;-\; (\ref{eqn143})\;$  and have:

\[U_{jk}^{\otimes} (\gamma_{1}, \gamma_{2}, \gamma_{3}, \eta) \;=\;L[U_{jk}(\gamma_{1}, \gamma_{2}, \gamma_{3}, t)] \;\;\;\;\;\;\;      \rm{k=1,2,3}\] 

\begin{eqnarray}\label{eqn147}
\eta \szp{[}  U_{j2}^{\otimes}( \gamma_{1}, \gamma_{2}, \gamma_{3}, \eta )\; -\;\frac{\gamma_{2}}{\gamma_{1}} U_{j1}^{\otimes}( \gamma_{1}, \gamma_{2}, \gamma_{3}, \eta) \szp{]}  \;-\; \szp{[}   U_{j2}( \gamma_{1}, \gamma_{2}, \gamma_{3}, 0 )\; -\;\frac{\gamma_{2}}{\gamma_{1}} U_{j1}( \gamma_{1}, \gamma_{2}, \gamma_{3}, 0) \szp{]}  \;=  \nonumber
\\
\nonumber\\
-\nu( \gamma_{1}^{2} +\gamma_{2}^{2} +\gamma_{3}^{2})\szp{[}   U_{j2}^{\otimes}                          ( \gamma_{1}, \gamma_{2}, \gamma_{3}, \eta )\; -\;\frac{\gamma_{2}}{\gamma_{1}} U_{j1}^{\otimes}( \gamma_{1}, \gamma_{2}, \gamma_{3}, \eta) \szp{]}  \;+ 
\quad\quad\quad\quad\quad\quad 
\nonumber
\\
\nonumber\\
+\; \szp{[}   F_{j2}^{\otimes}( \gamma_{1}, \gamma_{2}, \gamma_{3}, \eta )\; -\;\frac{\gamma_{2}}{\gamma_{1}} F_{j1}^{\otimes}( \gamma_{1}, \gamma_{2}, \gamma_{3}, \eta) \szp{]}
\quad\quad\quad\quad\quad\quad \quad\quad\quad \quad\quad\quad 
\end{eqnarray}

\begin{eqnarray}\label{eqn148}
\eta  \szp{[}   U_{j3}^{\otimes}( \gamma_{1}, \gamma_{2}, \gamma_{3}, \eta )\; -\;\frac{\gamma_{3}}{\gamma_{1}} U_{j1}^{\otimes}( \gamma_{1}, \gamma_{2}, \gamma_{3}, \eta) \szp{]}   \;-\; \szp{[}    U_{j3}( \gamma_{1}, \gamma_{2}, \gamma_{3}, 0 )\; -\;\frac{\gamma_{3}}{\gamma_{1}} U_{j1}( \gamma_{1}, \gamma_{2}, \gamma_{3}, 0) \szp{]}   \;=  \nonumber
\\
\nonumber\\
-\nu( \gamma_{1}^{2} +\gamma_{2}^{2} +\gamma_{3}^{2})\szp{[}    U_{j3}^{\otimes}                          ( \gamma_{1}, \gamma_{2}, \gamma_{3}, \eta )\; -\;\frac{\gamma_{3}}{\gamma_{1}} U_{j1}^{\otimes}( \gamma_{1}, \gamma_{2}, \gamma_{3}, \eta) \szp{]}   \;+ 
\quad\quad\quad\quad\quad\quad 
\nonumber
\\
\nonumber\\
+\; \szp{[}    F_{j3}^{\otimes}( \gamma_{1}, \gamma_{2}, \gamma_{3}, \eta )\; -\;\frac{\gamma_{3}}{\gamma_{1}} F_{j1}^{\otimes}( \gamma_{1}, \gamma_{2}, \gamma_{3}, \eta) \szp{]} 
\quad\quad\quad\quad\quad\quad \quad\quad\quad \quad\quad\quad 
\end{eqnarray}

\begin{equation}\label{eqn149}
\gamma_{1} U_{j1}^{\otimes}( \gamma_{1}, \gamma_{2}, \gamma_{3}, \eta ) \;+\; \gamma_{2}\, U_{j2}^{\otimes}( \gamma_{1}, \gamma_{2}, \gamma_{3}, \eta ) \;+\; \gamma_{3}\, U_{j3}^{\otimes}( \gamma_{1}, \gamma_{2}, \gamma_{3}, \eta ) \;=\;0
\end{equation}

\begin{equation}\label{eqn150}
U_{j1}(\gamma_{1}, \gamma_{2},  \gamma_{3},  0)\;=\; U_{1}^{0}(\gamma_{1} ,\gamma_{2} ,\gamma_{3})
\end{equation}

\begin{equation}\label{eqn151}
U_{j2}(\gamma_{1}, \gamma_{2},  \gamma_{3},  0)\;=\; U_{2}^{0}(\gamma_{1} ,\gamma_{2} ,\gamma_{3})
\end{equation}

\begin{equation}\label{eqn152}
U_{j3}(\gamma_{1}, \gamma_{2},  \gamma_{3},  0)\;=\; U_{3}^{0}(\gamma_{1} ,\gamma_{2} ,\gamma_{3})
\end{equation}

In the usual way the solution of the system of equations $(\ref{eqn147})\;-\; (\ref{eqn149})\;$with formulas $(\ref{eqn150})\;-\; (\ref{eqn152})\;$can be rewritten in the following form:

\begin{eqnarray}\label{eqn153}
U_{j1}^{\otimes}( \gamma_{1}, \gamma_{2}, \gamma_{3}, \eta )\;=\;\frac{[( \gamma_{2}^{2} +\gamma_{3}^{2})  F_{j1}^{\otimes}( \gamma_{1}, \gamma_{2}, \gamma_{3}, \eta) - \gamma_{1}\gamma_{2} F_{j2}^{\otimes}( \gamma_{1}, \gamma_{2}, \gamma_{3}, \eta) - \gamma_{1}\gamma_{3} F_{j3}^{\otimes}( \gamma_{1}, \gamma_{2}, \gamma_{3}, \eta)]}{ (\gamma_{1}^{2} +\gamma_{2}^{2} +\gamma_{3}^{2}) [\eta+\nu (\gamma_{1}^{2} +\gamma_{2}^{2} +\gamma_{3}^{2})] }\;+\nonumber
\\
\nonumber\\
+\; \frac{ U_{1}^{0}(\gamma_{1} , \gamma_{2} , \gamma_{3})}{[\eta+\nu (\gamma_{1}^{2} +\gamma_{2}^{2} +\gamma_{3}^{2})] } 
\quad\quad\quad\quad\quad\quad \quad\quad\quad \quad\quad\quad 
\quad\quad\quad\quad\quad\quad 
\end{eqnarray}

\begin{eqnarray}\label{eqn154}
U_{j2}^{\otimes}( \gamma_{1}, \gamma_{2}, \gamma_{3}, \eta )\;=\;\frac{[( \gamma_{3}^{2} +\gamma_{1}^{2})  F_{j2}^{\otimes}( \gamma_{1}, \gamma_{2}, \gamma_{3}, \eta) - \gamma_{2}\gamma_{3} F_{j3}^{\otimes}( \gamma_{1}, \gamma_{2}, \gamma_{3}, \eta) - \gamma_{2}\gamma_{1} F_{j1}^{\otimes}( \gamma_{1}, \gamma_{2}, \gamma_{3}, \eta)]}{ (\gamma_{1}^{2} +\gamma_{2}^{2} +\gamma_{3}^{2}) [\eta+\nu (\gamma_{1}^{2} +\gamma_{2}^{2} +\gamma_{3}^{2})] }\;+\nonumber
\\
\nonumber\\
+\; \frac{ U_{2}^{0}(\gamma_{1} , \gamma_{2} , \gamma_{3})}{[\eta+\nu (\gamma_{1}^{2} +\gamma_{2}^{2} +\gamma_{3}^{2})] } 
\quad\quad\quad\quad\quad\quad \quad\quad\quad \quad\quad\quad 
\quad\quad\quad\quad\quad\quad 
\end{eqnarray}

\begin{eqnarray}\label{eqn155}
U_{j3}^{\otimes}( \gamma_{1}, \gamma_{2}, \gamma_{3}, \eta )\;=\;\frac{[( \gamma_{1}^{2} +\gamma_{2}^{2})  F_{j3}^{\otimes}( \gamma_{1}, \gamma_{2}, \gamma_{3}, \eta) - \gamma_{3}\gamma_{1} F_{j1}^{\otimes}( \gamma_{1}, \gamma_{2}, \gamma_{3}, \eta) - \gamma_{3}\gamma_{2} F_{j2}^{\otimes}( \gamma_{1}, \gamma_{2}, \gamma_{3}, \eta)]}{ (\gamma_{1}^{2} +\gamma_{2}^{2} +\gamma_{3}^{2}) [\eta+\nu (\gamma_{1}^{2} +\gamma_{2}^{2} +\gamma_{3}^{2})] }\;+ \nonumber
\\
\nonumber\\
+\; \frac{ U_{3}^{0}(\gamma_{1} , \gamma_{2} , \gamma_{3})}{[\eta+\nu (\gamma_{1}^{2} +\gamma_{2}^{2} +\gamma_{3}^{2})] } 
\quad\quad\quad\quad\quad\quad \quad\quad\quad \quad\quad\quad 
\quad\quad\quad\quad\quad\quad 
\end{eqnarray}

Then we use the convolution formula (\ref{A6}) and integral (\ref{A7}) for $(\ref{eqn153})\;-\; (\ref{eqn155})\;$and obtain:

\begin{eqnarray}\label{eqn156}
U_{j1}(\gamma_{1}, \gamma_{2}, \gamma_{3}, t)\;=\;
\quad\quad\quad\quad\quad\quad \quad\quad\quad \quad\quad\quad 
\quad\quad\quad\quad\quad\quad
 \nonumber
\\
\nonumber\\
\int_{0}^{t} \sz{e} ^{-\nu (\gamma_{1}^{2} +\gamma_{2}^{2} +\gamma_{3}^{2}) (t-\tau)} \frac{[( \gamma_{2}^{2} +\gamma_{3}^{2})  F_{j1}( \gamma_{1}, \gamma_{2}, \gamma_{3}, \tau) -\gamma_{1}\gamma_{2} F_{j2}( \gamma_{1}, \gamma_{2}, \gamma_{3}, \tau ) -\gamma_{1}\gamma_{3} F_{j3}( \gamma_{1}, \gamma_{2}, \gamma_{3}, \tau ) ]}{ (\gamma_{1}^{2} +\gamma_{2}^{2}+\gamma_{3}^{2} ) }\,d\tau\;+ \nonumber
\\
\nonumber\\
+\;\sz{e} ^{-\nu (\gamma_{1}^{2} +\gamma_{2}^{2} +\gamma_{3}^{2}) t} \;U_{1}^{0}(\gamma_{1} ,\gamma_{2} ,\gamma_{3})
\quad\quad\quad\quad\quad\quad \quad\quad\quad \quad\quad\quad 
\quad\quad\quad\quad 
\end{eqnarray}

\begin{eqnarray}\label{eqn157}
U_{j2}(\gamma_{1}, \gamma_{2}, \gamma_{3}, t)\;=\; 
\quad\quad\quad\quad\quad\quad \quad\quad\quad \quad\quad\quad 
\quad\quad\quad\quad\quad\quad
\nonumber
\\
\nonumber\\
\int_{0}^{t} \sz{e} ^{-\nu (\gamma_{1}^{2} +\gamma_{2}^{2} +\gamma_{3}^{2}) (t-\tau)} \frac{[( \gamma_{3}^{2} +\gamma_{1}^{2})  F_{j2}( \gamma_{1}, \gamma_{2}, \gamma_{3}, \tau) -\gamma_{2}\gamma_{3} F_{j3}( \gamma_{1}, \gamma_{2}, \gamma_{3}, \tau ) -\gamma_{2}\gamma_{1} F_{j1}( \gamma_{1}, \gamma_{2}, \gamma_{3}, \tau ) ]}{ (\gamma_{1}^{2} +\gamma_{2}^{2}+\gamma_{3}^{2} ) }\,d\tau\;+ \nonumber
\\
\nonumber\\
+\;\sz{e} ^{-\nu (\gamma_{1}^{2} +\gamma_{2}^{2} +\gamma_{3}^{2}) t} \;U_{2}^{0}(\gamma_{1} ,\gamma_{2} ,\gamma_{3})
\quad\quad\quad\quad\quad\quad \quad\quad\quad \quad\quad\quad 
\quad\quad\quad\quad 
\end{eqnarray}

\begin{eqnarray}\label{eqn158}
U_{j3}(\gamma_{1}, \gamma_{2}, \gamma_{3}, t)\;=\; 
\quad\quad\quad\quad\quad\quad \quad\quad\quad \quad\quad\quad 
\quad\quad\quad\quad\quad\quad
\nonumber
\\
\nonumber\\
\int_{0}^{t} \sz{e} ^{-\nu (\gamma_{1}^{2} +\gamma_{2}^{2} +\gamma_{3}^{2}) (t-\tau)} \frac{[( \gamma_{1}^{2} +\gamma_{2}^{2})  F_{j3}( \gamma_{1}, \gamma_{2}, \gamma_{3}, \tau) -\gamma_{3}\gamma_{1} F_{j1}( \gamma_{1}, \gamma_{2}, \gamma_{3}, \tau ) -\gamma_{3}\gamma_{2} F_{j2}( \gamma_{1}, \gamma_{2}, \gamma_{3}, \tau ) ]}{ (\gamma_{1}^{2} +\gamma_{2}^{2}+\gamma_{3}^{2} ) }\,d\tau\;+ 
\nonumber
\\
\nonumber\\
+\;\sz{e} ^{-\nu (\gamma_{1}^{2} +\gamma_{2}^{2} +\gamma_{3}^{2}) t} \;U_{3}^{0}(\gamma_{1} ,\gamma_{2} ,\gamma_{3})
\quad\quad\quad\quad\quad\quad \quad\quad\quad \quad\quad\quad 
\quad\quad\quad\quad 
\end{eqnarray}

$ P_{j}(\gamma_{1}, \gamma_{2}, \gamma_{3}, t) $ is obtained from $(\ref{eqn134})\;[(\ref{eqn135})\;$or$\; (\ref{eqn136})]\;$:

\begin{equation}\label{eqn159}
P_{j}(\gamma_{1}, \gamma_{2}, \gamma_{3}, t) \;=\;i\frac{[\gamma_{1} F_{j1}( \gamma_{1}, \gamma_{2}, \gamma_{3}, t) +\gamma_{2} F_{j2}( \gamma_{1}, \gamma_{2}, \gamma_{3}, t ) +\gamma_{3} F_{j3}( \gamma_{1}, \gamma_{2}, \gamma_{3}, t )]}{ (\gamma_{1}^{2} +\gamma_{2}^{2} +\gamma_{3}^{2}) }
\end{equation}
\\
Use of the Fourier inversion formula $(\ref{A3})$ and find:

\begin{eqnarray}\label{eqn160}
u_{j1}(x_{1}, x_{2}, x_{3}, t)\;=\; 
\frac{1}{(2\pi)^{3/2}} \int_{-\infty}^{\infty} \int_{-\infty}^{\infty} \int_{-\infty}^{\infty} \biggl[ \int_{0}^{t} \sz{e} ^{-\nu (\gamma_{1}^{2} +\gamma_{2}^{2} +\gamma_{3}^{2}) (t-\tau)} \frac{[ ( \gamma_{2}^{2} +\gamma_{3}^{2})  F_{j1}( \gamma_{1}, \gamma_{2}, \gamma_{3}, \tau)]} { (\gamma_{1}^{2} +\gamma_{2}^{2}+\gamma_{3}^{2} ) } \,d\tau \;- \nonumber
\\
\nonumber\\
-\; \int_{0}^{t} \sz{e} ^{-\nu (\gamma_{1}^{2} +\gamma_{2}^{2} +\gamma_{3}^{2}) (t-\tau)}\frac{ [\gamma_{1}\gamma_{2} F_{j2}( \gamma_{1}, \gamma_{2}, \gamma_{3}, \tau ) + \gamma_{1}\gamma_{3} F_{j3}( \gamma_{1}, \gamma_{2}, \gamma_{3}, \tau )]} { (\gamma_{1}^{2} +\gamma_{2}^{2}+\gamma_{3}^{2} ) }\,d\tau\;   +
\quad
\nonumber\\
\nonumber\\
+\;\sz{e} ^{-\nu (\gamma_{1}^{2} +\gamma_{2}^{2} +\gamma_{3}^{2}) t} \;U_{1}^{0}(\gamma_{1} ,\gamma_{2} ,\gamma_{3})\biggr] \;\sz{e} ^{-i(x_{1}\gamma_{1}+x_{2}\gamma_{2}+x_{3}\gamma_{3})}\,d\gamma_{1}d\gamma_{2}d\gamma_{3}\;=
\nonumber\\
\nonumber\\
=\;\frac{1}{8\pi^{3}} \int_{-\infty}^{\infty} \int_{-\infty}^{\infty} \int_{-\infty}^{\infty}\frac{( \gamma_{2}^{2} +\gamma_{3}^{2})} { (\gamma_{1}^{2} +\gamma_{2}^{2}+\gamma_{3}^{2} ) }  \biggl[ \int_{0}^{t} \sz{e} ^{-\nu (\gamma_{1}^{2} +\gamma_{2}^{2} +\gamma_{3}^{2}) (t-\tau)} \int_{-\infty}^{\infty}\int_{-\infty}^{\infty}\int_{-\infty}^{\infty}\sz{e}  ^{i(\tilde x_{1}\gamma_{1}+\tilde x_{2}\gamma_{2}+\tilde x_{3}\gamma_{3})} \cdot 
\nonumber\\
\nonumber\\
\cdot f_{j1}(\tilde x_{1},\tilde x_{2}, \tilde x_{3},\tau)\,d\tilde x_{1}d\tilde x_{2} d\tilde x_{3}d\tau\biggr]\sz{e} ^{-i(x_{1}\gamma_{1}+x_{2}\gamma_{2}+x_{3}\gamma_{3})}\,d\gamma_{1}d\gamma_{2}d\gamma_{3}\;-
\nonumber\\
\nonumber\\
-\;\frac{1}{8\pi^{3}} \int_{-\infty}^{\infty} \int_{-\infty}^{\infty} \int_{-\infty}^{\infty}\frac{ \gamma_{1}\gamma_{2}} { (\gamma_{1}^{2} +\gamma_{2}^{2}+\gamma_{3}^{2} ) }  \biggl[ \int_{0}^{t} \sz{e} ^{-\nu (\gamma_{1}^{2} +\gamma_{2}^{2} +\gamma_{3}^{2}) (t-\tau)} \int_{-\infty}^{\infty}\int_{-\infty}^{\infty}\int_{-\infty}^{\infty}\sz{e}  ^{i(\tilde x_{1}\gamma_{1}+\tilde x_{2}\gamma_{2}+\tilde x_{3}\gamma_{3})} \cdot 
\nonumber\\
\nonumber\\
\cdot f_{j2}(\tilde x_{1},\tilde x_{2}, \tilde x_{3},\tau)\,d\tilde x_{1}d\tilde x_{2} d\tilde x_{3}d\tau\biggr]\sz{e} ^{-i(x_{1}\gamma_{1}+x_{2}\gamma_{2}+x_{3}\gamma_{3})}\,d\gamma_{1}d\gamma_{2}d\gamma_{3}\;-
\nonumber\\
\nonumber\\
-\;\frac{1}{8\pi^{3}} \int_{-\infty}^{\infty} \int_{-\infty}^{\infty} \int_{-\infty}^{\infty}\frac{ \gamma_{1}\gamma_{3}} { (\gamma_{1}^{2} +\gamma_{2}^{2}+\gamma_{3}^{2} ) }  \biggl[ \int_{0}^{t} \sz{e} ^{-\nu (\gamma_{1}^{2} +\gamma_{2}^{2} +\gamma_{3}^{2}) (t-\tau)} \int_{-\infty}^{\infty}\int_{-\infty}^{\infty}\int_{-\infty}^{\infty}\sz{e}  ^{i(\tilde x_{1}\gamma_{1}+\tilde x_{2}\gamma_{2}+\tilde x_{3}\gamma_{3})} \cdot 
\nonumber\\
\nonumber\\
\cdot f_{j3}(\tilde x_{1},\tilde x_{2}, \tilde x_{3},\tau)\,d\tilde x_{1}d\tilde x_{2} d\tilde x_{3}d\tau\biggr]\sz{e} ^{-i(x_{1}\gamma_{1}+x_{2}\gamma_{2}+x_{3}\gamma_{3})}\,d\gamma_{1}d\gamma_{2}d\gamma_{3}\;+
\nonumber\\
\nonumber\\
+\;\frac{1}{8\pi^{3}} \int_{-\infty}^{\infty} \int_{-\infty}^{\infty} \int_{-\infty}^{\infty} \sz{e} ^{-\nu (\gamma_{1}^{2} +\gamma_{2}^{2} +\gamma_{3}^{2}) t}\biggl[  \int_{-\infty}^{\infty}\int_{-\infty}^{\infty}\int_{-\infty}^{\infty}\sz{e}  ^{i(\tilde x_{1}\gamma_{1}+\tilde x_{2}\gamma_{2}+\tilde x_{3}\gamma_{3})}
 \cdot \quad\quad\quad\quad\quad\quad\quad\quad\quad
\nonumber\\
\nonumber\\
\cdot \; u_{1}^{0}(\tilde x_{1},\tilde x_{2}, \tilde x_{3})\,d\tilde x_{1}d\tilde x_{2} d\tilde x_{3}\biggr]\sz{e} ^{-i(x_{1}\gamma_{1}+x_{2}\gamma_{2}+x_{3}\gamma_{3})}\,d\gamma_{1}d\gamma_{2}d\gamma_{3} \;= 
\nonumber\\
\nonumber\\
=\; S_{11}(f_{j1})\;+\; S_{12}(f_{j2})\;+\; S_{13}(f_{j3})\;+\;B(u_{1}^0) \quad\quad\quad\quad\quad\quad \quad\quad\quad \quad\quad\quad 
\end{eqnarray}

\begin{eqnarray}\label{eqn161}
u_{j2}(x_{1}, x_{2}, x_{3}, t)\;=\;
\frac{1}{(2\pi)^{3/2}} \int_{-\infty}^{\infty} \int_{-\infty}^{\infty} \int_{-\infty}^{\infty} \biggl[ \int_{0}^{t} \sz{e} ^{-\nu (\gamma_{1}^{2} +\gamma_{2}^{2} +\gamma_{3}^{2}) (t-\tau)} \frac{[( \gamma_{3}^{2} +\gamma_{1}^{2})  F_{j2}( \gamma_{1}, \gamma_{2}, \gamma_{3}, \tau)]} { (\gamma_{1}^{2} +\gamma_{2}^{2}+\gamma_{3}^{2} ) } \,d\tau \;-
\nonumber\\
\nonumber\\
-\; \int_{0}^{t} \sz{e} ^{-\nu (\gamma_{1}^{2} +\gamma_{2}^{2} +\gamma_{3}^{2}) (t-\tau)}\frac{[\gamma_{2}\gamma_{3} F_{j3}( \gamma_{1}, \gamma_{2}, \gamma_{3}, \tau ) +\gamma_{2}\gamma_{1} F_{j1}( \gamma_{1}, \gamma_{2}, \gamma_{3}, \tau )]} { (\gamma_{1}^{2} +\gamma_{2}^{2}+\gamma_{3}^{2} ) }\,d\tau\;   +
\quad
\nonumber\\
\nonumber\\
+\;\sz{e} ^{-\nu (\gamma_{1}^{2} +\gamma_{2}^{2} +\gamma_{3}^{2}) t} \;U_{2}^{0}(\gamma_{1} ,\gamma_{2} ,\gamma_{3})\biggr] \;\sz{e} ^{-i(x_{1}\gamma_{1}+x_{2}\gamma_{2}+x_{3}\gamma_{3})}\,d\gamma_{1}d\gamma_{2}d\gamma_{3}\;=
\nonumber\\
\nonumber\\
=\;-\;\frac{1}{8\pi^{3}} \int_{-\infty}^{\infty} \int_{-\infty}^{\infty} \int_{-\infty}^{\infty}\frac{ \gamma_{2}\gamma_{1}} { (\gamma_{1}^{2} +\gamma_{2}^{2}+\gamma_{3}^{2} ) }  \biggl[ \int_{0}^{t} \sz{e} ^{-\nu (\gamma_{1}^{2} +\gamma_{2}^{2} +\gamma_{3}^{2}) (t-\tau)} \int_{-\infty}^{\infty}\int_{-\infty}^{\infty}\int_{-\infty}^{\infty}\sz{e}  ^{i(\tilde x_{1}\gamma_{1}+\tilde x_{2}\gamma_{2}+\tilde x_{3}\gamma_{3})} \cdot 
\nonumber\\
\nonumber\\
\cdot f_{j1}(\tilde x_{1},\tilde x_{2}, \tilde x_{3},\tau)\,d\tilde x_{1}d\tilde x_{2} d\tilde x_{3}d\tau\biggr]\sz{e} ^{-i(x_{1}\gamma_{1}+x_{2}\gamma_{2}+x_{3}\gamma_{3})}\,d\gamma_{1}d\gamma_{2}d\gamma_{3}\;+
\nonumber\\
\nonumber\\
+\;\frac{1}{8\pi^{3}} \int_{-\infty}^{\infty} \int_{-\infty}^{\infty} \int_{-\infty}^{\infty}\frac{ (\gamma_{3}^2 + \gamma_{1}^2)} { (\gamma_{1}^{2} +\gamma_{2}^{2}+\gamma_{3}^{2} ) }  \biggl[ \int_{0}^{t} \sz{e} ^{-\nu (\gamma_{1}^{2} +\gamma_{2}^{2} +\gamma_{3}^{2}) (t-\tau)} \int_{-\infty}^{\infty}\int_{-\infty}^{\infty}\int_{-\infty}^{\infty}\sz{e}  ^{i(\tilde x_{1}\gamma_{1}+\tilde x_{2}\gamma_{2}+\tilde x_{3}\gamma_{3})} \cdot 
\nonumber\\
\nonumber\\
\cdot f_{j2}(\tilde x_{1},\tilde x_{2}, \tilde x_{3},\tau)\,d\tilde x_{1}d\tilde x_{2} d\tilde x_{3}d\tau\biggr]\sz{e} ^{-i(x_{1}\gamma_{1}+x_{2}\gamma_{2}+x_{3}\gamma_{3})}\,d\gamma_{1}d\gamma_{2}d\gamma_{3}\;-
\nonumber\\
\nonumber\\
-\;\frac{1}{8\pi^{3}} \int_{-\infty}^{\infty} \int_{-\infty}^{\infty} \int_{-\infty}^{\infty}\frac{ \gamma_{2}\gamma_{3}} { (\gamma_{1}^{2} +\gamma_{2}^{2}+\gamma_{3}^{2} ) }  \biggl[ \int_{0}^{t} \sz{e} ^{-\nu (\gamma_{1}^{2} +\gamma_{2}^{2} +\gamma_{3}^{2}) (t-\tau)} \int_{-\infty}^{\infty}\int_{-\infty}^{\infty}\int_{-\infty}^{\infty}\sz{e}  ^{i(\tilde x_{1}\gamma_{1}+\tilde x_{2}\gamma_{2}+\tilde x_{3}\gamma_{3})} \cdot 
\nonumber\\
\nonumber\\
\cdot f_{j3}(\tilde x_{1},\tilde x_{2}, \tilde x_{3},\tau)\,d\tilde x_{1}d\tilde x_{2} d\tilde x_{3}d\tau\biggr]\sz{e} ^{-i(x_{1}\gamma_{1}+x_{2}\gamma_{2}+x_{3}\gamma_{3})}\,d\gamma_{1}d\gamma_{2}d\gamma_{3}\;+
\nonumber\\
\nonumber\\
+\;\frac{1}{8\pi^{3}} \int_{-\infty}^{\infty} \int_{-\infty}^{\infty} \int_{-\infty}^{\infty} \sz{e} ^{-\nu (\gamma_{1}^{2} +\gamma_{2}^{2} +\gamma_{3}^{2}) t}\biggl[  \int_{-\infty}^{\infty}\int_{-\infty}^{\infty}\int_{-\infty}^{\infty}\sz{e}  ^{i(\tilde x_{1}\gamma_{1}+\tilde x_{2}\gamma_{2}+\tilde x_{3}\gamma_{3})}
 \cdot \quad\quad\quad\quad\quad\quad\quad\quad\quad
\nonumber\\
\nonumber\\
\cdot \; u_{2}^{0}(\tilde x_{1},\tilde x_{2}, \tilde x_{3})\,d\tilde x_{1}d\tilde x_{2} d\tilde x_{3}\biggr]\sz{e} ^{-i(x_{1}\gamma_{1}+x_{2}\gamma_{2}+x_{3}\gamma_{3})}\,d\gamma_{1}d\gamma_{2}d\gamma_{3} \;= 
\nonumber\\
\nonumber\\
=\; S_{21}(f_{j1})\;+\; S_{22}(f_{j2})\;+\; S_{23}(f_{j3})\;+\;B(u_{2}^0) 
\quad\quad\quad\quad\quad\quad \quad\quad\quad \quad\quad\quad 
\end{eqnarray}

\begin{eqnarray}\label{eqn162}
u_{j3}(x_{1}, x_{2}, x_{3}, t)\;=\;
\frac{1}{(2\pi)^{3/2}} \int_{-\infty}^{\infty} \int_{-\infty}^{\infty} \int_{-\infty}^{\infty} \biggl[ \int_{0}^{t} \sz{e} ^{-\nu (\gamma_{1}^{2} +\gamma_{2}^{2} +\gamma_{3}^{2}) (t-\tau)} \frac{[( \gamma_{1}^{2} +\gamma_{2}^{2})  F_{j3}( \gamma_{1}, \gamma_{2}, \gamma_{3}, \tau)]} { (\gamma_{1}^{2} +\gamma_{2}^{2}+\gamma_{3}^{2} ) } \,d\tau \;-
\nonumber\\
\nonumber\\
-\;\int_{0}^{t} \sz{e} ^{-\nu (\gamma_{1}^{2} +\gamma_{2}^{2} +\gamma_{3}^{2}) (t-\tau)}\frac{[\gamma_{3}\gamma_{1} F_{j1}( \gamma_{1}, \gamma_{2}, \gamma_{3}, \tau ) +\gamma_{3}\gamma_{2} F_{j2}( \gamma_{1}, \gamma_{2}, \gamma_{3}, \tau )]} { (\gamma_{1}^{2} +\gamma_{2}^{2}+\gamma_{3}^{2} ) }\,d\tau\;   +
\quad
\nonumber\\
\nonumber\\
+\; \sz{e} ^{-\nu (\gamma_{1}^{2} +\gamma_{2}^{2} +\gamma_{3}^{2}) t} \;U_{3}^{0}(\gamma_{1} ,\gamma_{2} ,\gamma_{3})\biggr] \;\sz{e} ^{-i(x_{1}\gamma_{1}+x_{2}\gamma_{2}+x_{3}\gamma_{3})}\,d\gamma_{1}d\gamma_{2}d\gamma_{3}\;=
\nonumber\\
\nonumber\\
=\;-\;\frac{1}{8\pi^{3}} \int_{-\infty}^{\infty} \int_{-\infty}^{\infty} \int_{-\infty}^{\infty}\frac{ \gamma_{3}\gamma_{1}} { (\gamma_{1}^{2} +\gamma_{2}^{2}+\gamma_{3}^{2} ) }  \biggl[ \int_{0}^{t} \sz{e} ^{-\nu (\gamma_{1}^{2} +\gamma_{2}^{2} +\gamma_{3}^{2}) (t-\tau)} \int_{-\infty}^{\infty}\int_{-\infty}^{\infty}\int_{-\infty}^{\infty}\sz{e}  ^{i(\tilde x_{1}\gamma_{1}+\tilde x_{2}\gamma_{2}+\tilde x_{3}\gamma_{3})} \cdot 
\nonumber\\
\nonumber\\
\cdot f_{j1}(\tilde x_{1},\tilde x_{2}, \tilde x_{3},\tau)\,d\tilde x_{1}d\tilde x_{2} d\tilde x_{3}d\tau\biggr]\sz{e} ^{-i(x_{1}\gamma_{1}+x_{2}\gamma_{2}+x_{3}\gamma_{3})}\,d\gamma_{1}d\gamma_{2}d\gamma_{3}\;-
\nonumber\\
\nonumber\\
-\;\frac{1}{8\pi^{3}} \int_{-\infty}^{\infty} \int_{-\infty}^{\infty} \int_{-\infty}^{\infty}\frac{ \gamma_{3}\gamma_{2}} { (\gamma_{1}^{2} +\gamma_{2}^{2}+\gamma_{3}^{2} ) }  \biggl[ \int_{0}^{t} \sz{e} ^{-\nu (\gamma_{1}^{2} +\gamma_{2}^{2} +\gamma_{3}^{2}) (t-\tau)} \int_{-\infty}^{\infty}\int_{-\infty}^{\infty}\int_{-\infty}^{\infty}\sz{e}  ^{i(\tilde x_{1}\gamma_{1}+\tilde x_{2}\gamma_{2}+\tilde x_{3}\gamma_{3})} \cdot 
\nonumber\\
\nonumber\\
\cdot f_{j2}(\tilde x_{1},\tilde x_{2}, \tilde x_{3},\tau)\,d\tilde x_{1}d\tilde x_{2} d\tilde x_{3}d\tau\biggr]\sz{e} ^{-i(x_{1}\gamma_{1}+x_{2}\gamma_{2}+x_{3}\gamma_{3})}\,d\gamma_{1}d\gamma_{2}d\gamma_{3}\;+
\nonumber\\
\nonumber\\
+\;\frac{1}{8\pi^{3}} \int_{-\infty}^{\infty} \int_{-\infty}^{\infty} \int_{-\infty}^{\infty}\frac{ (\gamma_{1}^2 + \gamma_{2}^2)} { (\gamma_{1}^{2} +\gamma_{2}^{2}+\gamma_{3}^{2} ) }  \biggl[ \int_{0}^{t} \sz{e} ^{-\nu (\gamma_{1}^{2} +\gamma_{2}^{2} +\gamma_{3}^{2}) (t-\tau)} \int_{-\infty}^{\infty}\int_{-\infty}^{\infty}\int_{-\infty}^{\infty}\sz{e}  ^{i(\tilde x_{1}\gamma_{1}+\tilde x_{2}\gamma_{2}+\tilde x_{3}\gamma_{3})} \cdot 
\nonumber\\
\nonumber\\
\cdot f_{j3}(\tilde x_{1},\tilde x_{2}, \tilde x_{3},\tau)\,d\tilde x_{1}d\tilde x_{2} d\tilde x_{3}d\tau\biggr]\sz{e} ^{-i(x_{1}\gamma_{1}+x_{2}\gamma_{2}+x_{3}\gamma_{3})}\,d\gamma_{1}d\gamma_{2}d\gamma_{3}\;+
\nonumber\\
\nonumber\\
+\;\frac{1}{8\pi^{3}} \int_{-\infty}^{\infty} \int_{-\infty}^{\infty} \int_{-\infty}^{\infty} \sz{e} ^{-\nu (\gamma_{1}^{2} +\gamma_{2}^{2} +\gamma_{3}^{2}) t}\biggl[  \int_{-\infty}^{\infty}\int_{-\infty}^{\infty}\int_{-\infty}^{\infty}\sz{e}  ^{i(\tilde x_{1}\gamma_{1}+\tilde x_{2}\gamma_{2}+\tilde x_{3}\gamma_{3})}
 \cdot \quad\quad\quad\quad\quad\quad\quad\quad\quad
\nonumber\\
\nonumber\\
\cdot \; u_{3}^{0}(\tilde x_{1},\tilde x_{2}, \tilde x_{3})\,d\tilde x_{1}d\tilde x_{2} d\tilde x_{3}\biggr]\sz{e} ^{-i(x_{1}\gamma_{1}+x_{2}\gamma_{2}+x_{3}\gamma_{3})}\,d\gamma_{1}d\gamma_{2}d\gamma_{3} \;= 
\nonumber\\
\nonumber\\
=\; S_{31}(f_{j1})\;+\; S_{32}(f_{j2})\;+\; S_{33}(f_{j3})\;+\;B(u_{3}^0) 
\quad\quad\quad\quad\quad\quad \quad\quad\quad \quad\quad\quad 
\end{eqnarray}

\begin{eqnarray}\label{eqn163}
p_{j}\,(x_{1}, x_{2}, x_{3}, t)\;=\; 
\frac{i}{(2\pi)^{3/2}} \int_{-\infty}^{\infty} \int_{-\infty}^{\infty} \int_{-\infty}^{\infty}\biggl[ \;\frac{[\gamma_{1} F_{j1}( \gamma_{1}, \gamma_{2}, \gamma_{3}, t) + \gamma_{2} F_{j2}( \gamma_{1}, \gamma_{2}, \gamma_{3}, t )]} { (\gamma_{1}^{2} +\gamma_{2}^{2} +\gamma_{3}^{2}) }\;+ 
\nonumber\\
\nonumber\\
+\;\frac{ \gamma_{3} F_{j3}( \gamma_{1}, \gamma_{2}, \gamma_{3}, t )}{  (\gamma_{1}^{2} +\gamma_{2}^{2} +\gamma_{3}^{2}) }\; \biggr] \;\sz{e} ^{-i(x_{1}\gamma_{1} +x_{2}\gamma_{2} +x_{3}\gamma_{3})}\,d\gamma_{1} d\gamma_{2} d\gamma_{3}\; = \; 
\nonumber\\
\nonumber\\
=\frac{i}{8\pi^{3}} \int_{-\infty}^{\infty} \int_{-\infty}^{\infty} \int_{-\infty}^{\infty}\frac{ \gamma_{1}} { (\gamma_{1}^{2} +\gamma_{2}^{2}+\gamma_{3}^{2} ) }  \biggl[ \int_{-\infty}^{\infty}\int_{-\infty}^{\infty}\int_{-\infty}^{\infty}\sz{e}  ^{i(\tilde x_{1}\gamma_{1}+\tilde x_{2}\gamma_{2}+\tilde x_{3}\gamma_{3})} \cdot 
\nonumber\\
\nonumber\\
\cdot f_{j1}(\tilde x_{1},\tilde x_{2}, \tilde x_{3},t)\,d\tilde x_{1}d\tilde x_{2} d\tilde x_{3}\biggr]\sz{e} ^{-i(x_{1}\gamma_{1}+x_{2}\gamma_{2}+x_{3}\gamma_{3})}\,d\gamma_{1}d\gamma_{2}d\gamma_{3}\;+
\nonumber\\
\nonumber\\
+\;\frac{i}{8\pi^{3}} \int_{-\infty}^{\infty} \int_{-\infty}^{\infty} \int_{-\infty}^{\infty}\frac{ \gamma_{2}} { (\gamma_{1}^{2} +\gamma_{2}^{2}+\gamma_{3}^{2} ) } \biggl [ \int_{-\infty}^{\infty}\int_{-\infty}^{\infty}\int_{-\infty}^{\infty}\sz{e}  ^{i(\tilde x_{1}\gamma_{1}+\tilde x_{2}\gamma_{2}+\tilde x_{3}\gamma_{3})} \cdot 
\nonumber\\
\nonumber\\
\cdot f_{j2}(\tilde x_{1},\tilde x_{2}, \tilde x_{3},t)\,d\tilde x_{1}d\tilde x_{2} d\tilde x_{3}\biggr]\sz{e} ^{-i(x_{1}\gamma_{1}+x_{2}\gamma_{2}+x_{3}\gamma_{3})}\,d\gamma_{1}d\gamma_{2}d\gamma_{3}\;+
\nonumber\\
\nonumber\\
+\;\frac{i}{8\pi^{3}} \int_{-\infty}^{\infty} \int_{-\infty}^{\infty} \int_{-\infty}^{\infty}\frac{ \gamma_{3}} { (\gamma_{1}^{2} +\gamma_{2}^{2}+\gamma_{3}^{2} ) }  \biggl[ \int_{-\infty}^{\infty}\int_{-\infty}^{\infty}\int_{-\infty}^{\infty}\sz{e}  ^{i(\tilde x_{1}\gamma_{1}+\tilde x_{2}\gamma_{2}+\tilde x_{3}\gamma_{3})} \cdot 
\nonumber\\
\nonumber\\
\cdot f_{j3}(\tilde x_{1},\tilde x_{2}, \tilde x_{3},t)\,d\tilde x_{1}d\tilde x_{2} d\tilde x_{3}\biggr]\sz{e} ^{-i(x_{1}\gamma_{1}+x_{2}\gamma_{2}+x_{3}\gamma_{3})}\,d\gamma_{1}d\gamma_{2}d\gamma_{3}\;=
\nonumber\\
\nonumber\\
=\;\tilde S_{1}(f_{j1})\;+\; \tilde S_{2}(f_{j2})\;+\; \tilde S_{3}(f_{j3})
\quad\quad\quad\quad\quad\quad \quad\quad \quad\quad\quad\quad\quad\quad\quad\quad \quad\quad  
\end{eqnarray}

So, the integrals $(\ref{eqn160})\;-\; (\ref{eqn163})\;$exist by the restrictions $(\ref{eqn17})\;, (\ref{eqn18})\;$.

Here $S_{11}(), S_{12}(), S_{13}(), S_{21}(), S_{22}(), S_{23}(), S_{31}(), S_{32}(), S_{33}(), B(), \tilde S_{1}(), \tilde S_{2}(), \tilde S_{3}()$ are the integral - operators.

\[S_{12}()\;= \;S_{21}() \]\[S_{13}()\;= \;S_{31}() \]\[S_{23}()\;= \;S_{32}() \]

We have for the vector $\vec{u}_{j}$ from the equations $(\ref{eqn160})\;-\; (\ref{eqn162})\;$:

\begin{equation}\label{eqn164}
\vec{u}_{j}\;=\;\bar{\bar{S}}\;\cdot\;\vec{f}_{j}\;+\;B(\vec{u}^{0})\;,
\end{equation}

where $\;\bar{\bar{S}} \; $ is the matrix - operator:
\[ \left( \begin{array}{ccc}
S_{11} & S_{12} & S_{13} \\
S_{21} & S_{22} & S_{23} \\
S_{31} & S_{32} & S_{33}  \end{array} \right)\]

We put $\vec{f}_{j}$ from equation $(\ref{eqn20})$ into equation $(\ref{eqn164})$ and have:

\begin{eqnarray}\label{eqn165}
\vec{u}_{j} = \bar{\bar{S}}\cdot(\;\vec{f}\;-\;(\;\vec{u}_{j-1}\cdot\nabla)\vec{u}_{j-1})\;+\;B(\vec{u}^{0})\;=
\nonumber\\
\nonumber\\
=\;\bar{\bar{S}}\cdot\vec{f} \;-\;\bar{\bar{S}}\cdot(\vec{u}_{j-1}\;\cdot\;\nabla\;)\;\vec{u}_{j-1}\;+\;B(\vec{u}^{0})\; =
\nonumber\\
\nonumber\\
=\;\vec{u}_{1}\;-\;\bar{\bar{S}}\cdot(\vec{u}_{j-1}\;\cdot\;\nabla)\;\vec{u}_{j-1}
\quad\quad\quad\quad\quad\quad \quad\quad 
\end{eqnarray}

Here $\vec{u}_{1}\;$ is the solution of the system of equations $(\ref{eqn13})\; - \;(\ref{eqn20})$ with condition:

\[\sum_{n=1}^{3} u_{n}\frac{\partial u_{k}}{\partial x_{n}}\;=\;0\;\;\;\;\;\;\;       \rm{k=1,2,3}\;\;\] 

For j = 1 formula $(\ref{eqn164})$ can be written as follows:

\begin{equation}\label{eqn168}
\vec{u}_{1}\;=\;\bar{\bar{S}}\;\cdot\;\vec{f}_{1}\;+\;B(\vec{u}^{0})\;,\;\;\;\;\;\\\vec{f}_{1}(x,t)\; = \; \vec{f}(x,t)
\end{equation}

If t $\rightarrow$ 0 then $\vec{u}_{1} \rightarrow \vec{u}^{0}$ (look at integral-operators $\bar{\bar{S}}, B()\;\;$- integrals $\;(\ref{eqn160})\; - \;(\ref{eqn162})$).

For j = 2  we define from equation $(\ref{eqn20})$:

\begin{equation}
\vec{f}_{2}(x,t)\; = \; \vec{f}_{1}(x,t) \; - \;(\;\vec{u}_{1}\;\cdot\;\nabla\;)\;\vec{u}_{1}\;
\end{equation}

We denote: 

\begin{equation}\label{eqn169}
\vec{f}_{2}^{*}\;=\;(\vec{u}_{1}\;\cdot\;\nabla)\;\vec{u}_{1}
\end{equation}

and then we have: 

\begin{equation}
\vec{f}_{2}(x,t)\; = \; \vec{f}_{1}(x,t) \; - \vec{f}_{2}^{*}
\end{equation}

Then we get $\vec{u}_{2}$ from $(\ref{eqn164}),(\ref{eqn168})$:

\begin{equation}\label{eqn171}
\vec{u}_{2}\;=\;\bar{\bar{S}}\;\cdot\;\vec{f}_{2}  \;+\;B(\vec{u}^{0})\;=\;\bar{\bar{S}}\;\cdot\;(\vec{f}_{1}\;-\;\vec{f}_{2}^{*})  \;+\;B(\vec{u}^{0})\;=\;\vec{u}_{1}\;-\;\vec{u}_{2}^{*}
\end{equation}

Here we have:

\begin{equation}\label{eqn170}
\vec{u}_{2}^{*}\;=\;\bar{\bar{S}}\;\cdot\;\vec{f}_{2}^{*}
\end{equation}

If t $\rightarrow$ 0 then $\vec{u}_{2}^{*} \rightarrow$ 0 (look at integral-operator $\bar{\bar{S}}\;\;$- integrals $\;(\ref{eqn160})\; - \;(\ref{eqn162})$).

Continue for j = 3. We define from equation $(\ref{eqn20})$:

\begin{equation}
\vec{f}_{3}(x,t)\; = \; \vec{f}_{1}(x,t) \; - \;(\;\vec{u}_{2}\;\cdot\;\nabla\;)\;\vec{u}_{2}\;
\end{equation}

Here we have:

\begin{equation}\label{eqn172}
(\vec{u}_{2}\;\cdot\;\nabla)\;\vec{u}_{2}\;=
\;((\vec{u}_{1}\;-\;\vec{u}_{2}^{*})\;\cdot\;\nabla\;)
\;(\vec{u}_{1}\;-\;\vec{u}_{2}^{*})\;=
\;\vec{f}_{2}^{*}\;+\;\vec{f}_{3}^{*}
\end{equation}

We denote in $(\ref{eqn172})$:

\[\vec{f}_{3}^{*}\;=\;-\;  (\vec{u}_{1}\;\cdot\;\nabla)\;\vec{u}_{2}^{*}\; -\; (\vec{u}_{2}^{*}\;\cdot\;\nabla)\;\vec{u}_{1}\;+\;
(\vec{u}_{2}^{*}\;\cdot\;\nabla)\;\vec{u}_{2}^{*}\]

and then we have: 

\begin{equation}
\vec{f}_{3}(x,t)\; = \; \vec{f}_{1}(x,t) \; - \vec{f}_{2}^{*}\; - \vec{f}_{3}^{*}
\end{equation}

Then we get $\vec{u}_{3}$ from $(\ref{eqn164}) , (\ref{eqn168}) ,(\ref{eqn170})$:

\begin{equation}\label{eqn174}
\vec{u}_{3}\;=\;\bar{\bar{S}}\;\cdot\;\vec{f}_{3} \;+\;B(\vec{u}^{0})\;=\;\bar{\bar{S}}\;\cdot\;(\vec{f}_{1} \;-\;\vec{f}_{2}^{*}\;-\;\vec{f}_{3}^{*})\;+\;B(\vec{u}^{0})\;=\;\vec{u}_{1} \;-\;\vec{u}_{2}^{*}\;-\;\vec{u}_{3}^{*}
\end{equation}

Here we denote:

\begin{equation}\label{eqn173}
\vec{u}_{3}^{*}\;=\;\bar{\bar{S}}\;\cdot\;\vec{f}_{3}^{*}
\end{equation}

If t $\rightarrow$ 0 then $\vec{u}_{3}^{*} \rightarrow$ 0 (look at integral-operator $\bar{\bar{S}}\;\;$- integrals $\;(\ref{eqn160})\; - \;(\ref{eqn162})$).

For j = 4. We define from equation $(\ref{eqn20})$:

\begin{equation}
\vec{f}_{4}(x,t)\; = \; \vec{f}_{1}(x,t) \; - \;(\;\vec{u}_{3}\;\cdot\;\nabla\;)\;\vec{u}_{3}\;
\end{equation}

Here we have:

\begin{equation}\label{eqn175}
(\vec{u}_{3}\;\cdot\;\nabla)\;\vec{u}_{3}\;=
\;((\vec{u}_{2}\;-\;\vec{u}_{3}^{*})\;\cdot\;\nabla\;)
\;(\vec{u}_{2}\;-\;\vec{u}_{3}^{*})\;=
\;\vec{f}_{2}^{*}\;+\;\vec{f}_{3}^{*}\;+\;\vec{f}_{4}^{*}
\end{equation}

We denote in $(\ref{eqn175})$:

\[\vec{f}_{4}^{*}\;=\;-\;  (\vec{u}_{2}\;\cdot\;\nabla)\;\vec{u}_{3}^{*}\; -\; (\vec{u}_{3}^{*}\;\cdot\;\nabla)\;\vec{u}_{2}\;+\;
(\vec{u}_{3}^{*}\;\cdot\;\nabla)\;\vec{u}_{3}^{*}\]

and then we have: 

\begin{equation}
\vec{f}_{4}(x,t)\; = \; \vec{f}_{1}(x,t) \; - \vec{f}_{2}^{*}\; - \vec{f}_{3}^{*}\; - \vec{f}_{4}^{*}
\end{equation}

Then we get $\vec{u}_{4}$ from $(\ref{eqn164}) , (\ref{eqn168}) ,(\ref{eqn170}) ,(\ref{eqn173})$:

\begin{equation}\label{eqn177}
\vec{u}_{4}\;=\;\bar{\bar{S}}\;\cdot\;(\vec{f}_{1} \; -\;\vec{f}_{2}^{*}\; -\;\vec{f}_{3}^{*}\;-\;\vec{f}_{4}^{*})\;+\;B(\vec{u}^{0})\;=\;\vec{u}_{1} \;-\;\vec{u}_{2}^{*}\;-\;\vec{u}_{3}^{*}\;-\;\vec{u}_{4}^{*}
\end{equation}

Here we denote:

\begin{equation}\label{eqn176}
\vec{u}_{4}^{*}\;=\;\bar{\bar{S}}\;\cdot\;\vec{f}_{4}^{*}
\end{equation}

If t $\rightarrow$ 0 then $\vec{u}_{4}^{*} \rightarrow$ 0 (look at integral-operator $\bar{\bar{S}}\;\;$- integrals $\;(\ref{eqn160})\; - \;(\ref{eqn162})$).

For arbitrary number j $(j \geq 2)$. We define from equation $(\ref{eqn20})$:

\begin{equation}
\vec{f}_{j}(x,t)\; = \; \vec{f}_{1}(x,t) \; - \;(\;\vec{u}_{j-1}\;\cdot\;\nabla\;)\;\vec{u}_{j-1}\;
\end{equation}

Here we have:

\begin{equation}\label{eqn181}
(\vec{u}_{j-1}\;\cdot\;\nabla)\;\vec{u}_{j-1}\;=\;\sum_{l=2}^{j} \vec{f}_{l}^{*}
\end{equation}

and it follows:

\begin{equation}\label{eqn181a}
\vec{f}_{j}\;=\;\vec{f}_{1}\;-\; \sum_{l=2}^{j} \vec{f}_{l}^{*}
\end{equation}

Then we get $\vec{u}_{j}$ from $(\ref{eqn164}) , (\ref{eqn168}) $

\begin{equation}\label{eqn182}
\vec{u}_{j}\;=\;\bar{\bar{S}}\;\cdot\; \vec{f}_{j}\;+ \;B(\vec{u}^{0})\;=\;\bar{\bar{S}}\;\cdot\;( \vec{f}_{1}\;-\; \sum_{l=2}^{j} \vec{f}_{l}^{*})\;+ \;B(\vec{u}^{0})\;=\;\vec{u}_{1}\;-\;\sum_{l=2}^{j} \vec{u}_{l}^{*}\
\end{equation}

Here we denote:

\begin{equation}\label{eqn183}
\vec{u}_{l}^{*}\;=\;\bar{\bar{S}}\;\cdot\;\vec{f}_{l}^{*}\;\;\;\;\;\;\;\;\;\;\;\;\;\;(2\;\leq\;l\;\leq\;j)
\end{equation}

If t $\rightarrow$ 0 then $\vec{u}_{l}^{*} \rightarrow$ 0 (look at integral-operator $\bar{\bar{S}}\;\;$- integrals $\;(\ref{eqn160})\; , \;(\ref{eqn162})$).

We consider the equations $(\ref{eqn168})$ - $(\ref{eqn183})$ and see that the series $(\ref{eqn182})$ converge for $j \rightarrow \infty$ 

with the conditions for the first step (j = 1) of the iterative process:

\[\;\;\;\sum_{n=1}^{3} u_{0n}\frac{\partial u_{0k}}{\partial x_{n}} = 0\;\;\;\;\;\;\;\;       \rm{k=1,2,3}\]

and conditions

\begin{equation}\label{eqn184}
C_{\alpha K}\leq\;\frac{1}{2}\;\;,\;\;C_{\alpha \beta K}\leq\;\frac{1}{2}.
\end{equation}
 Here $\;\;C_{\alpha K}\;\;$and$\;\;C_{\alpha \beta K}\;$ are received from $(\ref{eqn17})$ , $(\ref{eqn18}).$

Hence, we receive from equation $(\ref{eqn165})\;$ when $j \rightarrow \infty$:

\begin{equation}\label{eqn185}
\vec{u}_{\infty}=\;\vec{u}_{1}\;-\;\bar{\bar{S}}\cdot(\vec{u}_{\infty}\;\cdot\;\nabla)\;\vec{u}_{\infty}
\end{equation}

Equation $(\ref{eqn185})$ describes the converging iterative process. 

Then we have from formula $(\ref{eqn163})\;$:

\begin{equation}\label{eqn185a}
p_{\infty}\,\;=\; \;\tilde S_{1}(f_{\infty 1})\;+\; \tilde S_{2}(f_{\infty 2})\;+\; \tilde S_{3}(f_{\infty 3})
\end{equation}

Here $\vec{f}_{\infty}$ = ($f_{\infty 1} , f_{\infty 2}, f_{\infty 3}$) is received from formula $(\ref{eqn181a})\;$.

On the other hand we can transform the original system of differential equations $(\ref{eqn7})\; - \;(\ref{eqn9})$ to the equivalent system of integral equations by the scheme of iterative process $(\ref{eqn164})\;, \;(\ref{eqn165})$ for vector $\vec{u}$:

\begin{equation}\label{eqn186}
\vec{u}\;=\;\vec{u}_{1}\;-\;\bar{\bar{S}}\cdot(\vec{u}\;\cdot\;\nabla)\;\vec{u},
\end{equation}

where $\vec{u}_{1}$ is from formula $(\ref{eqn168})$.
We compare the equations $(\ref{eqn185})$ and $(\ref{eqn186})$ and see that the iterative process $(\ref{eqn185})$ converge to the solution of the system $(\ref{eqn186})$ and hence to the solution of the differential equations $(\ref{eqn7})\; - \;(\ref{eqn9})$ with conditions $(\ref{eqn184})$. 

\textbf{In other words there exist smooth functions}  $\mathbf{p_{\infty}(x, t)}$, $\mathbf{u_{\infty i}(x, t)}$  \textbf{(i = 1, 2, 3) on} $\mathbf{R^{3} \times [0,\infty)}$ \textbf{that satisfy} $\mathbf{(\ref{eqn1}), (\ref{eqn2}), (\ref{eqn3})}$ \textbf{and}

\begin{equation}\label{eqn186b}
\mathbf{p_{\infty}, \;u_{\infty i} \in  C^{\infty}(R^{3} \times [0,\infty)),}
\nonumber\\
\nonumber\\
\end{equation}

\begin{equation}\label{eqn186c}
\mathbf{\int_{R^{3}}|\vec{u}_{\infty}(x, t)|^{2}dx < C } 
\end{equation}

\textbf{for all t} $\mathbf{\geq 0}$.
\nonumber\\

In the following chapters 4 and 5 we describe in further details examples of the solutions for the Navier-Stocks and Euler problems with various applied forces and different values of the viscosity coefficient ${\nu}$.
\nonumber\\

\section{Example of the solution of the Cauchy problem  for the Euler equations by the described iterative method with a particular applied force (N = 2)}\

We will consider an example of the solution of the Cauchy problem  for the Euler equations ${(\nu = 0)}$ for  N = 2 and with initial conditions:

\begin{equation}\label{eqn400}
\vec{u}(x,0)\; = \; \vec{u}^{0}(x)\; = \;0\;\;\;\;\;\;\;\; (x\in R^{2})
\end{equation}

Hence, and from formulas $(\ref{eqn39}), (\ref{eqn40})$ for arbitrary step j of the iterative process, it follows:

\begin{eqnarray}\label{eqn401}
u_{j1}(x_{1},x_{2},t)\;=\;\frac{1}{4\pi^2} \int_{-\infty}^{\infty}\int_{-\infty}^{\infty}\frac{\gamma_{2}^{2} }{(\gamma_{1}^{2}+\gamma_{2}^{2}) }   \int_{0}^{t} \int_{-\infty}^{\infty}\int_{-\infty}^{\infty}\sz{e}  ^{i(\tilde x_{1}\gamma_{1}+\tilde x_{2}\gamma_{2})} f_{j1}(\tilde x_{1},\tilde x_{2},\tau)\,d \tilde x_{1}d \tilde  x_{2}d\tau\cdot
\nonumber\\
\nonumber\\
\nonumber\\
 \;\cdot\;\sz{e}  ^{-i(x_{1}\gamma_{1}+x_{2}\gamma_{2})}\,d\gamma_{1}d\gamma_{2}\;-
\quad\quad\quad\quad\quad\quad \quad\quad \quad\quad
\nonumber\\
\nonumber\\
\nonumber\\
- \;\frac{1}{4\pi^2} \int_{-\infty}^{\infty}\int_{-\infty}^{\infty}\frac{\gamma_{1} \gamma_{2} }{(\gamma_{1}^{2}+\gamma_{2}^{2}) }  \int_{0}^{t} \int_{-\infty}^{\infty}\int_{-\infty}^{\infty}\sz{e}  ^{i(\tilde x_{1}\gamma_{1}+\tilde x_{2}\gamma_{2})} f_{j2}(\tilde x_{1},\tilde x_{2},\tau)\,d \tilde x_{1}d \tilde x_{2}d\tau\cdot
\nonumber\\
\nonumber\\
\nonumber\\
 \;\cdot\;\sz{e}  ^{-i(x_{1}\gamma_{1}+x_{2}\gamma_{2})}\,d\gamma_{1}d\gamma_{2}
\quad\quad\quad\quad\quad\quad \quad\quad \quad\quad
\nonumber\\
\quad\quad\quad\quad\quad\quad \quad\quad \quad\quad 
 \end{eqnarray}

\begin{eqnarray}\label{eqn402}
u_{j2}(x_{1},x_{2},t)\;=\;-\frac{1}{4\pi^2} \int_{-\infty}^{\infty}\int_{-\infty}^{\infty}\frac{\gamma_{1} \gamma_{2} }{(\gamma_{1}^{2}+\gamma_{2}^{2}) } \int_{0}^{t} \int_{-\infty}^{\infty}\int_{-\infty}^{\infty}\sz{e}  ^{i(\tilde x_{1}\gamma_{1}+\tilde x_{2}\gamma_{2})} f_{j1}(\tilde x_{1},\tilde x_{2},\tau)\,d \tilde x_{1}d \tilde x_{2}d\tau\cdot
\nonumber\\
\nonumber\\
\nonumber\\
 \;\cdot\;\sz{e}  ^{-i(x_{1}\gamma_{1}+x_{2}\gamma_{2})}\,d\gamma_{1}d\gamma_{2}\;+
\quad\quad\quad\quad\quad\quad \quad\quad \quad\quad
\nonumber\\
\nonumber\\
\nonumber\\
+ \;\frac{1}{4\pi^2} \int_{-\infty}^{\infty}\int_{-\infty}^{\infty}\frac{\gamma_{1}^{2} }{(\gamma_{1}^{2}+\gamma_{2}^{2}) } \int_{0}^{t} \int_{-\infty}^{\infty}\int_{-\infty}^{\infty}\sz{e}  ^{i(\tilde x_{1}\gamma_{1}+\tilde x_{2}\gamma_{2})} f_{j2}(\tilde x_{1},\tilde x_{2},\tau)\,d \tilde x_{1}d \tilde x_{2}d\tau \cdot
\nonumber\\
\nonumber\\
\nonumber\\
 \;\cdot\;\sz{e}  ^{-i(x_{1}\gamma_{1}+x_{2}\gamma_{2})}\,d\gamma_{1}d\gamma_{2}
\quad\quad\quad\quad\quad\quad \quad\quad \quad\quad
\nonumber\\
\quad\quad\quad\quad\quad\quad \quad\quad \quad\quad 
\end{eqnarray} 

We convert the Cartesian coordinates to the polar coordinates by formulas:

$x_{1}$ = $r\;\cdot\;$cos$\varphi\;$;$\;\;x_{2}$ = $r\;\cdot\;$sin$\varphi\;$;$\;\;\gamma_{1}$ = $\rho\;\cdot\;$cos$\psi\;$;$\;\;\gamma_{2}$ = $\rho\;\cdot\;$sin$\psi\;$;$\;\;\tilde x_{1}$ = $\tilde r\;\cdot\;$cos$\tilde \varphi\;$;$\;\;\tilde x_{2}$ = $\tilde r\;\cdot\;$sin$\tilde \varphi$;

and obtain from formulas $(\ref{eqn401}), (\ref{eqn402})$:

\begin{eqnarray}\label{eqn403}
u_{j1}(r,\varphi,t)\;=\;\frac{1}{4\pi^2} \int_{0}^{\infty}\int_{0}^{2\pi}{sin^{2}\psi } \int_{0}^{t} \int_{0}^{\infty}\int_{0}^{2\pi}\sz{e}  ^{i\tilde r\rho cos(\tilde \varphi-\psi)} f_{j1}(\tilde r,\tilde \varphi,\tau)\,\tilde r d \tilde r d \tilde  \varphi d\tau\cdot
\nonumber\\
\nonumber\\
\nonumber\\
 \;\cdot\;\sz{e}  ^{-ir\rho cos(\psi-\varphi)}\,\rho d\rho d\psi\;-
\quad\quad\quad\quad\quad\quad\quad
\nonumber\\
\nonumber\\
\nonumber\\
- \;\frac{1}{4\pi^2} \int_{0}^{\infty}\int_{0}^{2\pi}{sin\psi cos\psi }  \int_{0}^{t} \int_{0}^{\infty}\int_{0}^{2\pi}\sz{e}  ^{i\tilde r\rho cos(\tilde \varphi-\psi)} f_{j2}(\tilde r,\tilde \varphi,\tau)\,\tilde r d \tilde r d \tilde  \varphi d\tau\cdot
\nonumber\\
\nonumber\\
\nonumber\\
 \;\cdot\;
\sz{e}  ^{-ir\rho cos(\psi-\varphi)}\,\rho d\rho d\psi
\quad\quad\quad\quad\quad\quad\quad\quad\quad
\nonumber\\
\nonumber\\
\nonumber\\
\quad\quad\quad\quad\quad\quad\quad\quad\quad
\end{eqnarray}

\begin{eqnarray}\label{eqn404}
u_{j2}(r,\varphi,t)\;=\;-\frac{1}{4\pi^2} \int_{0}^{\infty}\int_{0}^{2\pi}{sin\psi cos\psi } \int_{0}^{t} \int_{0}^{\infty}\int_{0}^{2\pi}\sz{e}  ^{i\tilde r\rho cos(\tilde \varphi-\psi)} f_{j1}(\tilde r,\tilde \varphi,\tau)\,\tilde r d \tilde r d \tilde  \varphi d\tau\cdot
\nonumber\\
\nonumber\\
\nonumber\\
 \;\;\;\;\cdot\;\sz{e}  ^{-ir\rho cos(\psi-\varphi)}\,\rho d\rho d\psi\;+
\quad\quad\quad\quad\quad\quad \quad\quad \quad\quad
\nonumber\\
\nonumber\\
\nonumber\\
+ \;\frac{1}{4\pi^2} \int_{0}^{\infty}\int_{0}^{2\pi}cos^{2}\psi \int_{0}^{t} \int_{0}^{\infty}\int_{0}^{2\pi}\sz{e}  ^{i\tilde r\rho cos(\tilde \varphi-\psi)} f_{j2}(\tilde r,\tilde \varphi,\tau)\,\tilde r d \tilde r d \tilde  \varphi d\tau\cdot
\nonumber\\
\nonumber\\
\nonumber\\
 \;\;\;\;\cdot\;\sz{e}  ^{-ir\rho cos(\psi-\varphi)}\,\rho d\rho d\psi
\quad\quad\quad\quad\quad\quad \quad\quad \quad\quad
\nonumber\\
\quad\quad\quad\quad\quad\quad \quad\quad \quad\quad 
\end{eqnarray} 

We have the applied force $\vec{f}_{j}$ for arbitrary step j of the iterative process:

\begin{equation}\label{eqn405}
f_{j\tilde r}(\tilde r,\tilde \varphi,\tau) = f_{j\tilde r}(\tilde r)\sz{e}^{in_{j}\tilde \varphi}f_{j\tau}(\tau)\;\;\;,\;\;\;f_{j\tilde \varphi}(\tilde r,\tilde \varphi,\tau) \equiv 0
\end{equation}

or

\begin{equation}\label{eqn406}
f_{j\tilde r}(\tilde r,\tilde \varphi,\tau)\equiv 0\;\;\;,\;\;\;f_{j\tilde \varphi}(\tilde r,\tilde \varphi,\tau) = f_{j\tilde \varphi}(\tilde r)\sz{e}^{in_{j}\tilde \varphi}f_{j\tau}(\tau)
\end{equation}

where $\;\;f_{j\tilde r}(\tilde r,\tilde \varphi,\tau)\;\;\;,\;\;\;f_{j\tilde \varphi}(\tilde r,\tilde \varphi,\tau)\;\;\; - \;\;\;$  radial and tangential components of  the applied force. 

$n_{j}$ - separate circumferential mode, $n_{j}$ = 0,1,2,3,...

We take the radial and tangential components of  the applied force $(\ref{eqn405}), (\ref{eqn406})$ with condition $(\ref{eqn18})\;$.  For the radial component of  the applied force we use De Moivre's formulas (\ref{A8}) and have:

\begin{eqnarray}\label{eqn407}
f_{j1}(\tilde r,\tilde \varphi,\tau) = f_{j\tilde r}(\tilde r)\sz{e}^{in_{j}\tilde \varphi}cos \tilde \varphi f_{j\tau}(\tau) = \frac{1}{2}f_{j\tilde r}(\tilde r)\bigl(\sz{e}^{i(n_{j}-1)\tilde \varphi} + \sz{e}^{i(n_{j}+1)\tilde \varphi}\bigr)f_{j\tau}(\tau)
\nonumber\\
\nonumber\\
f_{j2}(\tilde r,\tilde \varphi,\tau) = f_{j\tilde r}(\tilde r)\sz{e}^{in_{j}\tilde \varphi}sin \tilde \varphi f_{j\tau}(\tau) = \frac{i}{2}f_{j\tilde r}(\tilde r)\bigl(\sz{e}^{i(n_{j}-1)\tilde \varphi} - \sz{e}^{i(n_{j}+1)\tilde \varphi}\bigr)f_{j\tau}(\tau)
\end{eqnarray}

We put the applied force components $(\ref{eqn407})$ in formulas $(\ref{eqn403}), (\ref{eqn404})$ , change the order of integration and find:

\begin{eqnarray}\label{eqn408}
u_{jr1}(r,\varphi,t)\;=\;\frac{1}{8\pi^2} \biggl[ \int_{0}^{\infty}\int_{0}^{2\pi}{sin^{2}\psi } \int_{0}^{\infty} f_{j\tilde r}(\tilde r) \int_{0}^{2\pi}\sz{e}  ^{i\tilde r\rho cos(\tilde \varphi-\psi)} (\sz{e}^{i(n_{j}-1)\tilde \varphi} + \sz{e}^{i(n_{j}+1)\tilde \varphi})d \tilde  \varphi \tilde r d \tilde r  \cdot
\nonumber\\
\nonumber\\
\nonumber\\
\cdot\;\sz{e}  ^{-ir\rho cos(\psi-\varphi)}\,\rho d\rho d\psi\;-
\quad\quad\quad\quad\quad\quad\quad
\nonumber\\
\nonumber\\
\nonumber\\
-\;\; i \int_{0}^{\infty}\int_{0}^{2\pi}{sin\psi cos\psi }  \int_{0}^{\infty} f_{j\tilde r}(\tilde r) \int_{0}^{2\pi}\sz{e}  ^{i\tilde r\rho cos(\tilde \varphi-\psi)}(\sz{e}^{i(n_{j}-1)\tilde \varphi} - \sz{e}^{i(n_{j}+1)\tilde \varphi})d \tilde  \varphi \tilde r d \tilde r  \cdot
\nonumber\\
\nonumber\\
\nonumber\\
 \;\cdot\;\sz{e}  ^{-ir\rho cos(\psi-\varphi)}\,\rho d\rho d\psi \biggr] \int_{0}^{t}f_{j\tau}(\tau) d\tau
\quad\quad\quad\quad\quad\quad\quad\quad\quad
\nonumber\\
\nonumber\\
\nonumber\\
\quad\quad\quad\quad\quad\quad\quad\quad\quad
\end{eqnarray}

\begin{eqnarray}\label{eqn409}
u_{jr2}(r,\varphi,t)\;=\;\frac{1}{8\pi^2} \biggl[ - \int_{0}^{\infty}\int_{0}^{2\pi}{sin\psi cos\psi} \int_{0}^{\infty} f_{j\tilde r}(\tilde r) \int_{0}^{2\pi}\sz{e}  ^{i\tilde r\rho cos(\tilde \varphi-\psi)} (\sz{e}^{i(n_{j}-1)\tilde \varphi} + \sz{e}^{i(n_{j}+1)\tilde \varphi})d \tilde  \varphi \tilde r d \tilde r  \cdot
\nonumber\\
\nonumber\\
\nonumber\\
\cdot\;\sz{e}  ^{-ir\rho cos(\psi-\varphi)}\,\rho d\rho d\psi\;+
\quad\quad\quad\quad\quad\quad\quad
\nonumber\\
\nonumber\\
\nonumber\\
+\;\; i \int_{0}^{\infty}\int_{0}^{2\pi}{ cos^{2}\psi }  \int_{0}^{\infty} f_{j\tilde r}(\tilde r) \int_{0}^{2\pi}\sz{e}  ^{i\tilde r\rho cos(\tilde \varphi-\psi)}(\sz{e}^{i(n_{j}-1)\tilde \varphi} - \sz{e}^{i(n_{j}+1)\tilde \varphi})d \tilde  \varphi \tilde r d \tilde r  \cdot
\nonumber\\
\nonumber\\
\nonumber\\
 \;\cdot\;
\sz{e}  ^{-ir\rho cos(\psi-\varphi)}\,\rho d\rho d\psi \biggr] \int_{0}^{t}f_{j\tau}(\tau) d\tau
\quad\quad\quad\quad\quad\quad\quad\quad\quad
\nonumber\\
\nonumber\\
\nonumber\\
\quad\quad\quad\quad\quad\quad\quad\quad\quad
\end{eqnarray}

We denote:

\begin{equation}\label{eqn410}
u_{jt}(t) = \int_{0}^{t}f_{j\tau}(\tau) d\tau
\end{equation}

and from formulas $(\ref{eqn408}), (\ref{eqn409})$ it follows:

\begin{eqnarray}\label{eqn411}
u_{jr1}(r,\varphi,t)\;=\;u_{jr1}(r,\varphi)u_{jt}(t)
\nonumber\\
\nonumber\\
u_{jr2}(r,\varphi,t)\;=\;u_{jr2}(r,\varphi)u_{jt}(t)
\end{eqnarray}

where

\begin{eqnarray}\label{eqn412}
u_{jr1}(r,\varphi)\;=\;\frac{1}{8\pi^2} \biggl[ \int_{0}^{\infty}\int_{0}^{2\pi}{sin^{2}\psi } \int_{0}^{\infty} f_{j\tilde r}(\tilde r) \int_{0}^{2\pi}\sz{e}  ^{i\tilde r\rho cos(\tilde \varphi-\psi)} (\sz{e}^{i(n_{j}-1)\tilde \varphi} + \sz{e}^{i(n_{j}+1)\tilde \varphi})d \tilde  \varphi \tilde r d \tilde r  \cdot
\nonumber\\
\nonumber\\
\nonumber\\
\cdot\;\sz{e}  ^{-ir\rho cos(\psi-\varphi)}\,\rho d\rho d\psi\;-
\quad\quad\quad\quad\quad\quad\quad
\nonumber\\
\nonumber\\
\nonumber\\
-\;\; i \int_{0}^{\infty}\int_{0}^{2\pi}{sin\psi cos\psi }  \int_{0}^{\infty} f_{j\tilde r}(\tilde r) \int_{0}^{2\pi}\sz{e}  ^{i\tilde r\rho cos(\tilde \varphi-\psi)}(\sz{e}^{i(n_{j}-1)\tilde \varphi} - \sz{e}^{i(n_{j}+1)\tilde \varphi})d \tilde  \varphi \tilde r d \tilde r  \cdot
\nonumber\\
\nonumber\\
\nonumber\\
 \;\cdot\;
\sz{e}  ^{-ir\rho cos(\psi-\varphi)}\,\rho d\rho d\psi \biggr] 
\quad\quad\quad\quad\quad\quad\quad\quad\quad
\nonumber\\
\nonumber\\
\nonumber\\
\quad\quad\quad\quad\quad\quad\quad\quad\quad
\end{eqnarray}

\begin{eqnarray}\label{eqn413}
u_{jr2}(r,\varphi)\;=\;\frac{1}{8\pi^2} \biggl[ - \int_{0}^{\infty}\int_{0}^{2\pi}{sin\psi cos\psi} \int_{0}^{\infty} f_{j\tilde r}(\tilde r) \int_{0}^{2\pi}\sz{e}  ^{i\tilde r\rho cos(\tilde \varphi-\psi)} (\sz{e}^{i(n_{j}-1)\tilde \varphi} + \sz{e}^{i(n_{j}+1)\tilde \varphi})d \tilde  \varphi \tilde r d \tilde r  \cdot
\nonumber\\
\nonumber\\
\nonumber\\
\cdot\;\sz{e}  ^{-ir\rho cos(\psi-\varphi)}\,\rho d\rho d\psi\;+
\quad\quad\quad\quad\quad\quad\quad
\nonumber\\
\nonumber\\
\nonumber\\
+\;\; i \int_{0}^{\infty}\int_{0}^{2\pi}{ cos^{2}\psi }  \int_{0}^{\infty} f_{j\tilde r}(\tilde r) \int_{0}^{2\pi}\sz{e}  ^{i\tilde r\rho cos(\tilde \varphi-\psi)}(\sz{e}^{i(n_{j}-1)\tilde \varphi} - \sz{e}^{i(n_{j}+1)\tilde \varphi})d \tilde  \varphi \tilde r d \tilde r  \cdot
\nonumber\\
\nonumber\\
\nonumber\\
 \;\cdot\;
\sz{e}  ^{-ir\rho cos(\psi-\varphi)}\,\rho d\rho d\psi \biggr] 
\quad\quad\quad\quad\quad\quad\quad\quad\quad
\nonumber\\
\nonumber\\
\nonumber\\
\quad\quad\quad\quad\quad\quad\quad\quad\quad
\end{eqnarray}

Let us denote  internal integrals in $(\ref{eqn412}), (\ref{eqn413})$ as $\;\;I_{\underline{+}}(\tilde r,\rho,\psi)\;\;$: 

\begin{equation}\label{eqn414}
I_{\underline{+}}(\tilde r,\rho,\psi) = \int_{0}^{2\pi}\sz{e}  ^{i\tilde r\rho cos(\tilde \varphi-\psi)} (\sz{e}^{i(n_{j}-1)\tilde \varphi\;\;} \underline{+}\;\; \sz{e}^{i(n_{j}+1)\tilde \varphi})d \tilde  \varphi
\end{equation}

 We have two integrals here. Plus (+) is for the first part of each integral $(\ref{eqn412}), (\ref{eqn413})$ and minus (-) is for the second part.

We substitute $\tilde \theta$ for $\tilde \varphi$:$\;\;\;\tilde \theta\;$ = $\;\tilde \varphi\;$ - $\psi\;\;$ , d$\tilde \theta\;$ = d$\tilde \varphi\;$ and receive:

\begin{equation}\label{eqn415}
I_{\underline{+}}(\tilde r,\rho,\psi) = \sz{e}^{i(n_{j}-1)\psi}\int_{-\psi}^{2\pi-\psi}\sz{e}  ^{i\tilde r\rho cos\tilde \theta + i(n_{j}-1)\tilde \theta} d \tilde \theta \;\;\underline{+}\;\;\sz{e}^{i(n_{j}+1)\psi}\int_{-\psi}^{2\pi-\psi}\sz{e}  ^{i\tilde r\rho cos\tilde \theta + i(n_{j}+1)\tilde \theta} d \tilde \theta
\end{equation}

Then we use the Bessel function's integral representation $(\ref{A9})$ and have:

\begin{equation}\label{eqn416}
I_{\underline{+}}(\tilde r,\rho,\psi) = 2\pi i^{(n_{j}-1)}\sz{e}^{i(n_{j}-1)\psi}J_{n_{j}-1}(\tilde{r}\rho)\;\;\underline{+}\;\;2\pi i^{(n_{j}+1)}\sz{e}^{i(n_{j}+1)\psi}J_{n_{j}+1}(\tilde{r}\rho)
\end{equation}

Put $\;\;I_{\underline{+}}(\tilde r,\rho,\psi)\;\;$ from $(\ref{eqn416})$ in formulas $(\ref{eqn412}), (\ref{eqn413})\;$ , change order of integration and obtain:

\begin{equation}\label{eqn417}
u_{jr1}(r,\varphi)\;=\;\frac{1}{8\pi^2}\int_{0}^{\infty}\int_{0}^{\infty}f_{j\tilde r}(\tilde r)\int_{0}^{2\pi}\biggl[{sin^{2} \psi} I_{+}(\tilde r,\rho,\psi)-i\;{sin\psi cos\psi}I_{-}(\tilde r,\rho,\psi)\biggr]\;\sz{e}  ^{-ir\rho cos(\psi-\varphi)}d\psi \tilde r d \tilde r \rho d\rho 
\nonumber\\
\nonumber\\
\quad\quad
\end{equation}

\begin{equation}\label{eqn418}
u_{jr2}(r,\varphi)\;=\;\frac{1}{8\pi^2}\int_{0}^{\infty}\int_{0}^{\infty}f_{j\tilde r}(\tilde r)\int_{0}^{2\pi}\biggl[-{sin\psi}{cos\psi} I_{+}(\tilde r,\rho,\psi)+i\;{cos^{2}\psi}I_{-}(\tilde r,\rho,\psi)\biggr]\;\sz{e}  ^{-ir\rho cos(\psi-\varphi)}d\psi \tilde r d \tilde r \rho d\rho 
\nonumber\\
\nonumber\\
\quad\quad
\end{equation}

Then we group parts in brackets of formulas $(\ref{eqn417}), (\ref{eqn418})\;$ , use De Moivre's formulas $(\ref{A8})$ and the Bessel function's properties. And we get:

\begin{equation}\label{eqn419}
u_{jr1}(r,\varphi)\;=\;-\;\frac{n_{j}\;i^{n_{j}}}{2\pi}\int_{0}^{\infty}\int_{0}^{\infty}f_{j\tilde r}(\tilde r)\int_{0}^{2\pi}{sin\psi}\;\sz{e}  ^{-ir\rho cos(\psi-\varphi)+i n_{j} \psi}\;d\psi \;J_{n_{j}}(\tilde r \rho)\;d \tilde r d\rho 
\quad\quad
\end{equation}

\begin{equation}\label{eqn420}
u_{jr2}(r,\varphi)\;=\;\frac{n_{j}\;i^{n_{j}}}{2\pi}\int_{0}^{\infty}\int_{0}^{\infty}f_{j\tilde r}(\tilde r)\int_{0}^{2\pi}{cos\psi}\;\sz{e}  ^{-ir\rho cos(\psi-\varphi)+i n_{j} \psi}\;d\psi \;J_{n_{j}}(\tilde r \rho)\;d \tilde r d\rho 
\quad\quad
\end{equation}

We substitute $\theta$ for $\psi$:$\;\;\;\theta\;$ = $\;\psi\;$ - $\varphi\;\;$ , d$\theta\;$ = d$\psi\;$ in the internal integrals of formulas $(\ref{eqn419}), (\ref{eqn420})\;$, use De Moivre's formulas $(\ref{A8})$ and the Bessel function's integral representation $(\ref{A9})$ and have from formulas $(\ref{eqn419}), (\ref{eqn420})$:

\begin{equation}\label{eqn421}
u_{jr1}(r,\varphi)\;=\;\frac{n_{j}}{2}\;\sz{e}  ^{i n_{j} \varphi}\int_{0}^{\infty}\biggl[\sz{e}  ^{i \varphi}J_{n_{j+1}}(r \rho)\;+\;\sz{e}  ^{-i \varphi}J_{n_{j-1}}(r \rho)\biggr]\int_{0}^{\infty}f_{j\tilde r}(\tilde r)\;J_{n_{j}}(\tilde r \rho)\;d \tilde r d\rho 
\quad\quad
\end{equation}

\begin{equation}\label{eqn422}
u_{jr2}(r,\varphi)\;=\;\frac{i\; n_{j}}{2}\;\sz{e}  ^{i n_{j} \varphi}\int_{0}^{\infty}\biggl[\sz{e}  ^{i \varphi}J_{n_{j+1}}(r \rho)\;-\;\sz{e}  ^{-i \varphi}J_{n_{j-1}}(r \rho)\biggr]\int_{0}^{\infty}f_{j\tilde r}(\tilde r)\;J_{n_{j}}(\tilde r \rho)\;d \tilde r d\rho 
\quad\quad
\end{equation}

Let us denote:

\begin{equation}\label{eqn423}
R_{j,n_{j}-1,r}(r)\;=\;\int_{0}^{\infty}\int_{0}^{\infty}f_{j\tilde r}(\tilde r)\;J_{n_{j}}(\tilde r \rho)J_{n_{j}-1}(r \rho)\;d \tilde r d\rho 
\quad\quad
\end{equation}

\begin{equation}\label{eqn424}
R_{j,n_{j}+1,r}(r)\;=\;\int_{0}^{\infty}\int_{0}^{\infty}f_{j\tilde r}(\tilde r)\;J_{n_{j}}(\tilde r \rho)J_{n_{j}+1}(r \rho)\;d \tilde r d\rho 
\quad\quad
\end{equation}

Then we have from formulas $(\ref{eqn421}), (\ref{eqn422})$:

\begin{equation}\label{eqn425}
u_{jr1}(r,\varphi)\;=\;\frac{n_{j}}{2}\bigl[R_{j,n_{j}-1,r}(r)\;\sz{e}  ^{i (n_{j}-1) \varphi}\;+\;R_{j,n_{j}+1,r}(r)\;\sz{e}  ^{i (n_{j}+1) \varphi}\bigr]
\end{equation}

\begin{equation}\label{eqn426}
u_{jr2}(r,\varphi)\;=\;\frac{i\;n_{j}}{2}\bigl[R_{j,n_{j}-1,r}(r)\;\sz{e}  ^{i (n_{j}-1) \varphi}\;-\;R_{j,n_{j}+1,r}(r)\;\sz{e}  ^{i (n_{j}+1) \varphi}\bigr]
\end{equation}

Then if $n_{j}\;=\;0$ it follows from $(\ref{eqn423}), (\ref{eqn424}) , (\ref{eqn425}), (\ref{eqn426})$ that
$u_{jr1}(r,\varphi)\;=\;u_{jr2}(r,\varphi)\;=\;0\;$ and hence 
$u_{1}\;\;=\;\;u_{2}\;\;=\;\;0$.

In the equations bellow we will consider $n_{j}\;\geq\;1$.

We change the order of integration in formulas $(\ref{eqn423}), (\ref{eqn424}) $ and obtain:

\begin{equation}\label{eqn427}
R_{j,n_{j}-1,r}(r)\;=\;\int_{0}^{\infty}f_{j\tilde r}(\tilde r)\int_{0}^{\infty}J_{n_{j}}(\tilde r \rho)J_{n_{j}-1}(r \rho)\;d\rho\; d \tilde r  
\quad\quad
\end{equation}

\begin{equation}\label{eqn428}
R_{j,n_{j}+1,r}(r)\;=\;\int_{0}^{\infty}f_{j\tilde r}(\tilde r)\int_{0}^{\infty}J_{n_{j}}(\tilde r \rho)J_{n_{j}+1}(r \rho)\;d\rho\; d \tilde r  
\quad\quad
\end{equation}

Internal integrals in formulas $(\ref{eqn427}), (\ref{eqn428})$ are established by the discontinuous integral of Weber and Schafheitlin $(\ref{A17})\; \cite{gW44}$. Then we have from $(\ref{eqn427}), (\ref{eqn428})$:

\begin{equation}\label{eqn429}
R_{j,n_{j}-1,r}(r)\;=\;r^{n_{j}-1}  \int_{r}^{\infty}\frac{f_{j\tilde r}(\tilde r)}{\tilde r^{n_{j}}} d \tilde r  
\quad\quad
\end{equation}

\begin{equation}\label{eqn430}
R_{j,n_{j}+1,r}(r)\;=\;\frac{1}{r^{n_{j}+1}} \int_{0}^{r}\tilde r^{n_{j}} f_{j\tilde r}(\tilde r) d \tilde r  
\quad\quad
\end{equation}

Now we integrate the solution $(\ref{eqn403}), (\ref{eqn404})$ by the tangential component of  the applied force $(\ref{eqn406})$ for ($n_{j}\;\geq\;1$). Then we use De Moivre's formulas $(\ref{A8})$ and have:

\begin{eqnarray}\label{eqn431}
f_{j1}(\tilde r,\tilde \varphi,\tau) = - f_{j\tilde \varphi}(\tilde r)\sz{e}^{in_{j}\tilde \varphi}sin \tilde \varphi f_{j\tau}(\tau) = - \frac{i}{2}f_{j\tilde \varphi}(\tilde r) \bigl( \sz{e}^{i(n_{j}-1)\tilde \varphi} - \sz{e}^{i(n_{j}+1)\tilde \varphi} \bigr) f_{j\tau}(\tau)
\nonumber\\
\nonumber\\
f_{j2}(\tilde r,\tilde \varphi,\tau) \;\;=\;\; f_{j\tilde \varphi}(\tilde r)\sz{e}^{in_{j}\tilde \varphi}cos \tilde \varphi f_{j\tau}(\tau) \;\;= \; \frac{1}{2} f_{j\tilde \varphi}(\tilde r)\bigl(\sz{e}^{i(n_{j}-1)\tilde \varphi} + \sz{e}^{i(n_{j}+1)\tilde \varphi}\bigr)f_{j\tau}(\tau)
\nonumber\\
\nonumber\\
\end{eqnarray}

Hence formulas $(\ref{eqn431})$ are the components $f_{j1}$ and $f_{j2}$ from the tangential component of  the applied force $(\ref{eqn406})$, while formulas $(\ref{eqn407})$ are the components $f_{j1}$ and $f_{j2}$ from the radial component of  the applied force $(\ref{eqn405})$.

Let us put $(\ref{eqn431})$ in formulas $(\ref{eqn403}), (\ref{eqn404})$ and do the operations as we did in $(\ref{eqn408}) - (\ref{eqn426})\;  (\;n_{j}\;\geq\;1)$. We consider that $f_{j\tilde \varphi}(\tilde r)\cdot f_{j\tau}(\tau)$ is restricted by condition $(\ref{eqn18})$ and get:

\begin{equation}\label{eqn432}
R_{j,n_{j}-1,\varphi}(r)\;=\;-\;\int_{0}^{\infty}\int_{0}^{\infty}\bigl(f_{j\tilde \varphi}(\tilde r)\;\tilde r\bigr)^{'}_{\tilde r}\;J_{n_{j}}(\tilde r \rho)J_{n_{j}-1}(r \rho)\;d \tilde r d\rho 
\quad\quad
\end{equation}

\begin{equation}\label{eqn433}
R_{j,n_{j}+1,\varphi}(r)\;=\;-\;\int_{0}^{\infty}\int_{0}^{\infty}\bigl(f_{j\tilde \varphi}(\tilde r)\;\tilde r\bigr)^{'}_{\tilde r}\;J_{n_{j}}(\tilde r \rho)J_{n_{j}+1}(r \rho)\;d \tilde r d\rho \;,
\quad\quad
\end{equation}

Here $\bigl(\bigr)^{'}_{\tilde r} \equiv \frac{\partial}{\partial \tilde r}$. Hence we have:

\begin{equation}\label{eqn434}
u_{j\varphi1}(r,\varphi)\;=\;-\;\frac{i}{2}\bigl[R_{j,n_{j}-1,\varphi}(r)\;\sz{e}  ^{i (n_{j}-1) \varphi}\;+\;R_{j,n_{j}+1,\varphi}(r)\;\sz{e}  ^{i (n_{j}+1) \varphi}\bigr]
\end{equation}

\begin{equation}\label{eqn435}
u_{j\varphi2}(r,\varphi)\;=\;\frac{1}{2}\bigl[R_{j,n_{j}-1,\varphi}(r)\;\sz{e}  ^{i (n_{j}-1) \varphi}\;-\;R_{j,n_{j}+1,\varphi}(r)\;\sz{e}  ^{i (n_{j}+1) \varphi}\bigr]
\end{equation}

We change the order of integration in formulas $(\ref{eqn432}), (\ref{eqn433}) $ and obtain:

\begin{equation}\label{eqn436}
R_{j,n_{j}-1,\varphi}(r)\;=\;-\;\int_{0}^{\infty}\bigl(f_{j\tilde \varphi}(\tilde r)\;\tilde r\bigr)^{'}_{\tilde r}\int_{0}^{\infty}\;J_{n_{j}}(\tilde r \rho)J_{n_{j}-1}(r \rho)\;d\rho d \tilde r  
\quad\quad
\end{equation}

\begin{equation}\label{eqn437}
R_{j,n_{j}+1,\varphi}(r)\;=\;-\;\int_{0}^{\infty}\bigl(f_{j\tilde \varphi}(\tilde r)\;\tilde r\bigr)^{'}_{\tilde r}\int_{0}^{\infty}\;J_{n_{j}}(\tilde r \rho)J_{n_{j}+1}(r \rho)\; d\rho d \tilde r
\quad\quad
\end{equation}

Internal integrals in formulas $(\ref{eqn436}), (\ref{eqn437})$ are established by the discontinuous integral of Weber and Schafheitlin $(\ref{A17})\cite{gW44}$. Then we have from $(\ref{eqn436}), (\ref{eqn437})$:

\begin{equation}\label{eqn438}
R_{j,n_{j}-1,\varphi}(r)\;=\; -\; r^{n_{j}-1}  \int_{r}^{\infty}\bigl(f_{j\tilde \varphi}(\tilde r)\;\tilde r\bigr)^{'}_{\tilde r}\;\frac{d \tilde r}{\tilde r^{n_{j}}}   
\quad\quad
\end{equation}

\begin{equation}\label{eqn439}
R_{j,n_{j}+1,\varphi}(r)\;=\; -\; \frac{1}{r^{n_{j}+1}} \int_{0}^{r}\bigl(f_{j\tilde \varphi}(\tilde r)\;\tilde r\bigr)^{'}_{\tilde r}\;\tilde r^{n_{j}} d \tilde r  
\quad\quad
\end{equation}

We have obtain formulas $(\ref{eqn407})\;$ - $\;(\ref{eqn439})$ for an arbitrary step j of the iterative process and the applied forces $(\ref{eqn405})$ or $(\ref{eqn406})$. 

Now we investigate the first step (j = 1 , $n_{1} = n = 1,2,3,...$) of the iterative process with the particular radial component of  the applied force ${f}_{1}(x,t)\;\;$[ look at $\;(\ref{eqn47})$]:

\begin{eqnarray}\label{eqn440}
f_{1\tilde r}(\tilde r,\tilde \varphi,\tau) = f_{1\tilde r}(\tilde r)\sz{e}^{in\tilde \varphi}f_{1\tau}(\tau)\;\;\;,\;\;\;f_{1\tilde \varphi}(\tilde r,\tilde \varphi,\tau) \equiv 0
\nonumber\\
\nonumber\\
f_{1\tilde r}(\tilde r) = F_{n}\tilde r^{n+1}\sz{e}^{-\mu_{n}\tilde r}\;\;,\;\;f_{1\tau}(\tau) = \sz{e}^{-\sigma_{n}\tau}
\quad\quad
\end{eqnarray}

$F_{n}, \mu_{n}, \sigma_{n}$ - constants.  $0 < F_{n}< \infty$  , $1 < \mu_{n}< \infty$  ,  $1 < \sigma_{n}< \infty$.

Let us put the applied force $(\ref{eqn440})$ in formulas $(\ref{eqn429}), (\ref{eqn430})$, integrate and then we have:

\begin{eqnarray}\label{eqn441}
R_{1,n-1,r}(r)\;=\;r^{n-1}  \int_{r}^{\infty}\frac{F_{n}\tilde r^{n+1}\sz{e}^{-\mu_{n}\tilde r}}{\tilde r^{n}} d \tilde r\;=\;F_{n} r^{n-1}  \int_{r}^{\infty}\sz{e}^{-\mu_{n}\tilde r} \tilde r d \tilde r\;=\;\frac{F_{n} r^{n-1}}{\mu_{n}^2}\Gamma (2, \mu_{n}r)  
\nonumber\\
\nonumber\\
\quad\quad
\end{eqnarray}

\begin{eqnarray}\label{eqn442}
R_{1,n+1,r}(r)\;=\;\frac{1}{r^{n+1}} \int_{0}^{r} F_{n}\tilde r^{2n+1}\sz{e}^{-\mu_{n}\tilde r} d \tilde r\;=\;\;\frac{F_{n}}{r^{n+1}\mu_{n}^{2n+2}}\gamma(2n+2, \mu_{n}r)  
\nonumber\\
\nonumber\\
\quad\quad
\end{eqnarray}

$\Gamma (\alpha, x) , \gamma(\alpha, x)$ are the incomplete gamma functions $\cite{BE253}$.

Hence, and from formulas $(\ref{eqn425}), (\ref{eqn426})$ it follows by j = 1:

\begin{equation}\label{eqn443}
u_{1r1}(r,\varphi)\;=\;\frac{n}{2}\bigl[R_{1,n-1,r}(r)\;\sz{e}  ^{i (n-1) \varphi}\;+\;R_{1,n+1,r}(r)\;\sz{e}  ^{i (n+1) \varphi}\bigr]
\end{equation}

\begin{equation}\label{eqn444}
u_{1r2}(r,\varphi)\;=\;\frac{i\;n}{2}\bigl[R_{1,n-1,r}(r)\;\sz{e}  ^{i (n-1) \varphi}\;-\;R_{1,n+1,r}(r)\;\sz{e}  ^{i (n+1) \varphi}\bigr]
\end{equation}

and from formula $(\ref{eqn410})$:

\begin{equation}\label{eqn445}
u_{1t}(t) = \int_{0}^{t}f_{1\tau}(\tau) d\tau = \int_{0}^{t}\sz{e}^{-\sigma_{n}\tau}d\tau = \frac{1}{\sigma_{n}}\gamma(1,\sigma_{n}t)
\end{equation}

And we have the velocity $\vec u_{1}$ [look at $\;(\ref{eqn47})$] from formulas $(\ref{eqn411})$:

\begin{eqnarray}\label{eqn446}
u_{1r1}(r,\varphi,t)\;=\;u_{1r1}(r,\varphi)\;u_{1t}(t)
\nonumber\\
\nonumber\\
u_{1r2}(r,\varphi,t)\;=\;u_{1r2}(r,\varphi)\;u_{1t}(t)
\end{eqnarray}

We obtain the following equations by performing appropriate transformations:

\begin{equation}\label{eqn447}
u_{1 r}(r,\varphi, t)\;=\;\frac{n}{2}\bigl[R_{1,n-1,r}(r)\;+\;R_{1,n+1,r}(r)\bigr]\;\sz{e}  ^{i n \varphi}\;u_{1t}(t)
\end{equation}

\begin{equation}\label{eqn448}
u_{1\varphi}(r,\varphi, t)\;=\;\frac{i\;n}{2}\bigl[R_{1,n-1,r}(r)\;-\;R_{1,n+1,r}(r)\bigr]\;\sz{e}^{i n \varphi}\;u_{1t}(t)
\end{equation}

$u_{1 r}(r,\varphi, t) ,\;\; u_{1\varphi}(r,\varphi, t)$ are the radial and tangential components of the velocity $\vec u_{1}$.

We have from formula $(\ref{eqn445})$:

\begin{equation}
  \begin{array}{ll} 
				lim\;\;u_{1t}(t)= 0 \\
				t \rightarrow 0
					\end{array}     
\end{equation}
					
Hence, and from formulas $(\ref{eqn446})-(\ref{eqn448})$:

\begin{equation}
\begin{array}{ll} 
lim \;\;u_{1r1}(r,\varphi,t)= 0;\;\;\;   lim \;\;u_{1r2}(r,\varphi,t)= 0;\\
t \rightarrow 0; \;\;\;\;\;\;\;\;\;\;\;\;\;\;\;\;\;\;\;\;\;\;\;\;\;\; t \rightarrow 0;\\
lim \;\;u_{1r}(r,\varphi,t)= 0;\;\;\;   lim \;\;u_{1\varphi}(r,\varphi,t)= 0;\\
t \rightarrow 0; \;\;\;\;\;\;\;\;\;\;\;\;\;\;\;\;\;\;\;\;\;\;\;\;\;\; t \rightarrow 0;\\
\end{array}
\end{equation}

In other words the velocity $\vec u_{1}$ satisfies the initial conditions $(\ref{eqn400})$.

We use the asymptotic properties of the incomplete gamma functions $\Gamma (\alpha, x) ,\;\; \gamma(\alpha, x)$ and from formulas $(\ref{eqn441}) - (\ref{eqn444})\;,$ $(\ref{eqn447}),\;\; (\ref{eqn448})$ we have the velocity $\vec u_{1}$ satisfies conditions  $(\ref{eqn16})$ (for $r \;\rightarrow\; \infty )$.

Let us continue investigation for the second step (j = 2) of the iterative process. 

Find $\vec{f}_{2}^{*}(r,\varphi, t) = \{f_{21}^{*}, f_{22}^{*}\}$ - the first correction of the particular radial applied force ${f}_{1}(x,t)$ $(\ref{eqn440})$.

We have for $\vec{f}_{2}^{*}$ from formula $(\ref{eqn48})$:

\begin{equation}\label{eqn449}
f_{21}^{*} =  u_{1r1}\;\frac{\partial u_{1r1}}{\partial x_{1}}\;+\;u_{1r2}\;\frac{\partial u_{1r1}}{\partial x_{2}}
\end{equation}

\begin{equation}\label{eqn450}
f_{22}^{*} =  u_{1r1}\;\frac{\partial u_{1r2}}{\partial x_{1}}\;+\;u_{1r2}\;\frac{\partial u_{1r2}}{\partial x_{2}}
\end{equation}

where $u_{1r1},\;\;u_{1r2}$ are the components of $\vec{u_{1}}$ and were taken from formulas $(\ref{eqn446})$.

We have here:

\begin{eqnarray}\label{eqn451}
\frac{\partial u_{1r1}(r,\varphi, t)}{\partial x_{1}} \; = \; \frac{\partial u_{1r1}(r,\varphi, t)}{\partial r}\; \frac{\partial r}{\partial x_{1}} \; + \; \frac{\partial u_{1r1}(r,\varphi, t)}{\partial \varphi}\; \frac{\partial \varphi}{\partial x_{1}}
\nonumber\\
\nonumber\\
\frac{\partial u_{1r1}(r,\varphi, t)}{\partial x_{2}}\;  = \; \frac{\partial u_{1r1}(r,\varphi, t)}{\partial r}\; \frac{\partial r}{\partial x_{2}} \; + \; \frac{\partial u_{1r1}(r,\varphi, t)}{\partial \varphi}\; \frac{\partial \varphi}{\partial x_{2}}
\nonumber\\
\nonumber\\
\frac{\partial u_{1r2}(r,\varphi, t)}{\partial x_{1}}\;  = \; \frac{\partial u_{1r2}(r,\varphi, t)}{\partial r}\; \frac{\partial r}{\partial x_{1}} \; + \; \frac{\partial u_{1r2}(r,\varphi, t)}{\partial \varphi}\; \frac{\partial \varphi}{\partial x_{1}}
\nonumber\\
\nonumber\\
\frac{\partial u_{1r2}(r,\varphi, t)}{\partial x_{2}}\;  = \; \frac{\partial u_{1r2}(r,\varphi, t)}{\partial r}\; \frac{\partial r}{\partial x_{2}} \; + \; \frac{\partial u_{1r2}(r,\varphi, t)}{\partial \varphi}\; \frac{\partial \varphi}{\partial x_{2}}
\nonumber\\
\nonumber\\
\frac{\partial r}{\partial x_{1}} = cos \varphi \; , \; \frac{\partial \varphi}{\partial x_{1}} = \;-\; \frac{sin \varphi}{r}\;,\; \frac{\partial r}{\partial x_{2}} = sin \varphi \; , \;\frac{\partial \varphi}{\partial x_{2}} =  \frac{cos \varphi}{r}
\nonumber\\
\nonumber\\
\end{eqnarray}

Hence, we use formulas $(\ref{eqn443}) - (\ref{eqn446})$ for $u_{1r1}(r,\varphi, t),\;\; u_{1r2}(r,\varphi, t)$ and have from $(\ref{eqn451})$:

\begin{eqnarray}\label{eqn452}
\frac{\partial u_{1r1}(r,\varphi, t)}{\partial x_{1}} \; = \; \frac{n}{2}\biggl\{\bigl[R_{1,n-1,r}^{'}(r)\;\sz{e}  ^{i (n-1) \varphi}\;+\;R_{1,n+1,r}^{'}(r)\;\sz{e}  ^{i (n+1) \varphi}\bigr]\;cos \varphi\;+
\nonumber\\
+\;i\;\bigl[(n-1)R_{1,n-1,r}(r)\;\sz{e}  ^{i (n-1) \varphi}\;+\;(n+1)R_{1,n+1,r}(r)\;\sz{e}  ^{i (n+1) \varphi}\bigr]\bigl(-\frac{sin \varphi}{r} \bigr ) \biggr \}\;u_{1t}(t)
\nonumber\\
\nonumber\\
\frac{\partial u_{1r1}(r,\varphi, t)}{\partial x_{2}}\;  = \; \frac{n}{2}\biggl\{\bigl[R_{1,n-1,r}^{'}(r)\;\sz{e}  ^{i (n-1) \varphi}\;+\;R_{1,n+1,r}^{'}(r)\;\sz{e}  ^{i (n+1) \varphi}\bigr]\;sin \varphi\;+
\nonumber\\
+\;i\;\bigl[(n-1)R_{1,n-1,r}(r)\;\sz{e}  ^{i (n-1) \varphi}\;+\;(n+1)R_{1,n+1,r}(r)\;\sz{e}  ^{i (n+1) \varphi}\bigr]\frac{cos \varphi}{r}  \biggr \}\;u_{1t}(t)
\nonumber\\
\nonumber\\
\frac{\partial u_{1r2}(r,\varphi, t)}{\partial x_{1}}\;  = \; \frac{i n}{2}\biggl\{\bigl[R_{1,n-1,r}^{'}(r)\;\sz{e}  ^{i (n-1) \varphi}\;-\;R_{1,n+1,r}^{'}(r)\;\sz{e}  ^{i (n+1) \varphi}\bigr]\;cos \varphi\;+
\nonumber\\
+\;i\;\bigl[(n-1)R_{1,n-1,r}(r)\;\sz{e}  ^{i (n-1) \varphi}\;-\;(n+1)R_{1,n+1,r}(r)\;\sz{e}  ^{i (n+1) \varphi}\bigr]\bigl(-\frac{sin \varphi}{r} \bigr ) \biggr \}\;u_{1t}(t)
\nonumber\\
\nonumber\\
\frac{\partial u_{1r2}(r,\varphi, t)}{\partial x_{2}}\;  = \; \frac{i n}{2}\biggl\{\bigl[R_{1,n-1,r}^{'}(r)\;\sz{e}  ^{i (n-1) \varphi}\;-\;R_{1,n+1,r}^{'}(r)\;\sz{e}  ^{i (n+1) \varphi}\bigr]\;sin \varphi\;+
\nonumber\\
+\;i\;\bigl[(n-1)R_{1,n-1,r}(r)\;\sz{e}  ^{i (n-1) \varphi}\;-\;(n+1)R_{1,n+1,r}(r)\;\sz{e}  ^{i (n+1) \varphi}\bigr]\frac{cos \varphi}{r}  \biggr \}\;u_{1t}(t)
\nonumber\\
\nonumber\\
\end{eqnarray}

where 

\begin{eqnarray}\label{eqn453}
R_{1,n-1,r}^{'}(r)\;=\;\frac{d R_{1,n-1,r}(r)}{d r}\;=\;F_{n}\bigl [ \frac{(n-1) r^{n-2}}{\mu_{n}^2}\; \Gamma (2, \mu_{n}r) - r^{n} \sz{e}  ^{- \mu_{n}r}\bigr ] 
\quad\quad
\nonumber\\
\nonumber\\
R_{1,n+1,r}^{'}(r)\;=\;\frac{d R_{1,n+1,r}(r)}{d r}\;=\;F_{n}\bigl [- \frac{ (n+1) }{r^{n+2} \mu_{n}^{2n+2}}\; \gamma (2n+2, \mu_{n}r) + r^{n} \sz{e}  ^{- \mu_{n}r}\bigr ] 
\nonumber\\
\nonumber\\
\quad\quad
\end{eqnarray}

Let us put $u_{1r1},\; u_{1r2},\; \frac{\partial u_{1r1}}{\partial x_{1}},\; \frac{\partial u_{1r1}}{\partial x_{2}},\; \frac{\partial u_{1r2}}{\partial x_{1}},\; \frac{\partial u_{1r2}}{\partial x_{2}}\;$ from formulas $(\ref{eqn446}), (\ref{eqn452})$ in formulas $(\ref{eqn449}), (\ref{eqn450})$ for $f_{21}^{*} \;,\; f_{22}^{*}\;$. 
\\
\\
After completing appropriate operations we have:

\begin{equation}\label{eqn454}
f_{21}^{*} (r,\varphi, t) =  \frac{n^{2}}{2^{2}}\bigl [ T_{2,2n-1,r}(r)\;\sz{e}  ^{i (2n-1) \varphi} + T_{2,2n+1,r}(r)\;\sz{e}  ^{i (2n+1) \varphi} \bigr ] T_{2n}(t)
\end{equation}

\begin{equation}\label{eqn455}
f_{22}^{*} (r,\varphi, t) =  \frac{i n^{2}}{2^{2}}\bigl [ T_{2,2n-1,r}(r)\;\sz{e}  ^{i (2n-1) \varphi} - T_{2,2n+1,r}(r)\;\sz{e}  ^{i (2n+1) \varphi} \bigr ] T_{2n}(t)
\end{equation}

where

\begin{eqnarray}\label{eqn456}
T_{2,2n-1,r}(r) = \bigl [ R_{1,n-1,r}(r) + R_{1,n+1,r}(r) \bigr ] R_{1,n-1,r}^{'}(r) - \frac{(n-1)R_{1,n-1,r}(r)}{r} \bigl [ R_{1,n-1,r}(r) - R_{1,n+1,r}(r) \bigr ]
\quad\quad
\nonumber\\
\nonumber\\
T_{2,2n+1,r}(r) = \bigl [ R_{1,n-1,r}(r) + R_{1,n+1,r}(r) \bigr ] R_{1,n+1,r}^{'}(r) - \frac{(n+1)R_{1,n+1,r}(r)}{r} \bigl [ R_{1,n-1,r}(r) - R_{1,n+1,r}(r) \bigr ]
\quad\quad
\nonumber\\
\nonumber\\
\quad\quad
\end{eqnarray}

\begin{equation}\label{eqn457}
T_{2n} (t) = u_{1t}^{2} (t)
\end{equation}

We use formulas $(\ref{eqn441}) , (\ref{eqn442})$ for $R_{1,n-1,r}(r)\;,\;R_{1,n+1,r}(r)$ and $(\ref{eqn453})$ for  $R_{1,n-1,r}^{'}(r)\;,\;R_{1,n+1,r}^{'}(r)$ then do appropriate operations for $T_{2,2n-1,r}(r)\;,\;T_{2,2n+1,r}(r)$ and get:

\begin{eqnarray}\label{eqn458}
T_{2,2n-1,r}(r) = F_{n}^{2}\;\bigl [-\; \frac{ r^{2n-1}}{\mu_{n}^2}\;\sz{e}  ^{- \mu_{n}r}\; \Gamma (2, \mu_{n}r)\;+\;\frac{2(n-1)}{r^{3}\mu_{n}^{2n+4}}\;\gamma (2n+2, \mu_{n}r)\;\Gamma (2, \mu_{n}r)\;-
\quad\quad\quad\quad\quad\quad\quad\quad
\nonumber\\
-\;\frac{1}{r\mu_{n}^{2n+2}}\;\sz{e}  ^{- \mu_{n}r}\;\gamma (2n+2, \mu_{n}r) \bigr ]
\quad\quad\quad\quad\quad\quad\quad\quad\quad\quad\quad
\nonumber\\
\nonumber\\
T_{2,2n+1,r}(r) = F_{n}^{2}\;\bigl [ \frac{ r^{2n-1}}{\mu_{n}^2}\;\sz{e}  ^{- \mu_{n}r}\; \Gamma (2, \mu_{n}r)\;-\;\frac{2(n+1)}{r^{3}\mu_{n}^{2n+4}}\;\gamma (2n+2, \mu_{n}r)\;\Gamma (2, \mu_{n}r)\;+
\quad\quad\quad\quad\quad\quad\quad\quad
\nonumber\\
+\;\frac{1}{r\mu_{n}^{2n+2}}\;\sz{e}  ^{- \mu_{n}r}\;\gamma (2n+2, \mu_{n}r) \bigr ]
\quad\quad\quad\quad\quad\quad\quad\quad\quad\quad\quad
\nonumber\\
\nonumber\\
\quad\quad\quad\quad\quad\quad\quad
\end{eqnarray}

For radial $f_{2r}^{*}$ and tangential $f_{2\varphi}^{*}$ components of the first correction $\vec{f_{2}^{*}}(r, \varphi, t)$ of the particular radial applied force  we have:

\begin{equation}\label{eqn459}
f_{2r}^{*} (r,\varphi, t) =  \frac{n^{2}}{2^{2}}\bigl [ T_{2,2n-1,r}(r) + T_{2,2n+1,r}(r) \bigr ]\; \sz{e}  ^{i\; 2n \varphi}\;T_{2n}(t)\;=\;\frac{n^{2}}{2^{2}}T_{2,2n,r}(r)\; \sz{e}  ^{i\; 2n \varphi}\;T_{2n}(t)
\end{equation}

\begin{equation}\label{eqn460}
f_{2\varphi}^{*} (r,\varphi, t) =  \frac{i n^{2}}{2^{2}}\bigl [ T_{2,2n-1,r}(r) - T_{2,2n+1,r}(r) \bigr ]\; \sz{e}  ^{i\; 2n \varphi}\; T_{2n}(t)\;=\;\frac{i n^{2}}{2^{2}}T_{2,2n,\varphi}(r)\; \sz{e}  ^{i\; 2n \varphi}\;T_{2n}(t)
\end{equation}

Here

\begin{equation}\label{eqn461}
T_{2,2n,r}(r) = -\;F_{n}^{2}\;\frac{4}{r^{3}\mu_{n}^{2n+4}}\;\gamma (2n+2, \mu_{n}r)\;\Gamma (2, \mu_{n}r)\quad\quad\quad\quad\quad\quad\quad\quad\quad\quad\quad\quad\quad\quad
\nonumber\\
\nonumber\\
\end{equation}

\begin{eqnarray}\label{eqn462}
T_{2,2n,\varphi}(r) = -\;F_{n}^{2}\;\bigl [ \frac{ 2 r^{2n-1}}{\mu_{n}^2}\;\sz{e}  ^{- \mu_{n}r}\; \Gamma (2, \mu_{n}r)\;-\;\frac{4 n}{r^{3}\mu_{n}^{2n+4}}\;\gamma (2n+2, \mu_{n}r)\;\Gamma (2, \mu_{n}r)\;+
\nonumber\\
+\;\frac{2}{r\mu_{n}^{2n+2}}\;\sz{e}  ^{- \mu_{n}r}\;\gamma (2n+2, \mu_{n}r) \bigr ]
\quad\quad
\nonumber\\
\nonumber\\
\end{eqnarray}

We compare the particular radial applied force $\vec{f_{1}}\;$ from $\;(\ref{eqn440})$ with the first correction $\vec{f_{2}^{*}}\;$ from $\;((\ref{eqn459})- (\ref{eqn462}))$ of this particular radial applied force, and  we have:

\begin{equation}\label{eqn463}
\mid\vec{f_{2}^{*}}\mid\; <<\; \mid\vec{f_{1}}\mid
\end{equation}

with condition 

\begin{equation}\label{eqn464}
F_{n}\;\leq\;\frac{1}{n}
\end{equation}

After the first step of the iterative process (j = 1) we obtained the velocity $\vec{u_{1}}\;$  [ see $\; (\ref{eqn446})$].
\\  
Now we will calculate $\vec{u_{2}^{*}}\;$ - the first correction of the velocity $\vec{u_{1}} $. Solution of this problem has two stages. On the first stage we find the part of the first correction $\vec{u_{2r}^{*}}$, corresponding to the radial component  of the first correction of applied force  $f_{2r}^{*}\;$ from $\;(\ref{eqn459})$: 

\begin{equation}\label{eqn465}
f_{2r}^{*}(r,\varphi,t) = \;\frac{n^{2}}{2^{2}}T_{2,2n,r}(r)\; \sz{e}  ^{i\; 2n \varphi}\;T_{2n}(t)\;\;\;,\;\;\;f_{2\varphi}^{*}(r,\varphi,t) \equiv 0
\end{equation}

On the second stage we calculate the other part of the first correction $\vec{u_{2\varphi}^{*}}\;$, corresponding to the tangential component  of the first correction of applied force  $f_{2\varphi}^{*}\;$ from $\;(\ref{eqn460})$:

\begin{equation}\label{eqn466}
f_{2r}^{*}(r,\varphi,t) \equiv 0\;\;\;,\;\;\; f_{2\varphi}^{*}(r,\varphi,t) = \;\frac{i n^{2}}{2^{2}}T_{2,2n,\varphi}(r)\; \sz{e}  ^{i\; 2n\varphi}\;T_{2n}(t)
\end{equation}

In other words 

\begin{eqnarray}\label{eqn467}
\vec{u_{2}^{*}} = \vec{u_{2r}^{*}} + \vec{u_{2\varphi}^{*}}, \;\;\;\vec{u_{2r}^{*}} = \{u_{2r1}^{*}, u_{2r2}^{*}\}, \;\;\;\vec{u_{2\varphi}^{*}} = \{u_{2\varphi1}^{*}, u_{2\varphi2}^{*}\}.
\end{eqnarray}

First stage: we use formulas $(\ref{eqn425})\;,\;(\ref{eqn426})\;$ for components $\;u_{jr1}(r, \varphi),\; u_{jr2}(r, \varphi)$ and formulas $(\ref{eqn429})\;,\;(\ref{eqn430})$ for $R_{j,n_{j}-1,r}(r),\; R_{j,n_{j}+1,r}(r)$ for j = 2 and then formulas $(\ref{eqn459})\;,\;(\ref{eqn461})$ for $f_{2r}^{*} (r,\varphi, t)\;,\;T_{2,2n,r}(r)$. We do appropriate operations and have:

\begin{equation}\label{eqn468}
u_{2r1}^{*}(r,\varphi)\;=\;n\bigl[R_{2,2n-1,r}(r)\;\sz{e}  ^{i (2n-1) \varphi}\;+\;R_{2,2n+1,r}(r)\;\sz{e}  ^{i (2n+1) \varphi}\bigr]
\end{equation}

\begin{equation}\label{eqn469}
u_{2r2}^{*}(r,\varphi)\;=\;i\;n\bigl[R_{2,2n-1,r}(r)\;\sz{e}  ^{i (2n-1) \varphi}\;-\;R_{2,2n+1,r}(r)\;\sz{e}  ^{i (2n+1) \varphi}\bigr]
\end{equation}

Here 

\begin{equation}\label{eqn470}
R_{2,2n-1,r}(r)\;=\;\frac{n^{2}r^{2n-1}}{2^{2}}  \int_{r}^{\infty}\frac{T_{2,2n,r}(\tilde r)}{\tilde r^{2n}} d \tilde r  
\quad\quad
\end{equation}

\begin{equation}\label{eqn471}
R_{2,2n+1,r}(r)\;=\;\frac{n^{2}}{2^{2} r^{2n+1}} \int_{0}^{r}\tilde r^{2n} T_{2,2n,r}(\tilde r) d \tilde r  
\quad\quad
\end{equation}

Second stage: we use formulas $(\ref{eqn434})\;,\;(\ref{eqn435})$ for components $u_{j\varphi1}(r, \varphi),\; u_{j\varphi2}(r, \varphi)$ and formulas $(\ref{eqn438})\;,\;(\ref{eqn439})$ for $R_{j,n_{j}-1,\varphi}(r),\; R_{j,n_{j}+1,\varphi}(r)$ for j = 2 and then formulas $(\ref{eqn460})\;,\;(\ref{eqn462})$ for $f_{2\varphi}^{*} (r,\varphi, t)\;,\;T_{2,2n,\varphi}(r)$. We do appropriate operations and have:

\begin{equation}\label{eqn472}
u_{2\varphi1}^{*}(r,\varphi)\;=\;- \;\frac{i}{2}\bigl[R_{2,2n-1,\varphi}(r)\;\sz{e}  ^{i (2n-1) \varphi}\;+\;R_{2,2n+1,\varphi}(r)\;\sz{e}  ^{i (2n+1) \varphi}\bigr]
\end{equation}

\begin{equation}\label{eqn473}
u_{2\varphi2}^{*}(r,\varphi)\;=\;\frac{1}{2}\bigl[R_{2,2n-1,\varphi}(r)\;\sz{e}  ^{i (2n-1) \varphi}\;-\;R_{2,2n+1,\varphi}(r)\;\sz{e}  ^{i (2n+1) \varphi}\bigr]
\end{equation}

Here:

\begin{equation}\label{eqn474}
R_{2,2n-1,\varphi}(r)\;=\;-\;\frac{i n^{2}r^{2n-1}}{2^{2}}  \int_{r}^{\infty}\frac{(T_{2,2n,\varphi}(\tilde r)\cdot \tilde r)_{\tilde r}^{'}}{\tilde r^{2n}} d \tilde r  
\quad\quad
\end{equation}

\begin{equation}\label{eqn475}
R_{2,2n+1,\varphi}(r)\;=\;-\;\frac{i n^{2}}{2^{2} r^{2n+1}} \int_{0}^{r}\tilde r^{2n} (T_{2,2n,\varphi}(\tilde r)\cdot  \tilde r)_{\tilde r}^{'} d \tilde r  
\quad\quad
\end{equation}

Then we have:

\begin{eqnarray}\label{eqn476}
u_{21}^{*}(r,\varphi) = u_{2r1}^{*}(r,\varphi) + u_{2\varphi1}^{*}(r,\varphi)\;=
\quad\quad\quad\quad\quad\quad\quad\quad\quad\quad\quad\quad\quad\quad\quad\
\nonumber\\
\nonumber\\
=\;\bigl [ n R_{2,2n-1,r}(r) - \frac{i}{2} R_{2,2n-1,\varphi}(r)\bigr ]\;\sz{e}  ^{i (2n-1) \varphi} + \bigl [ n R_{2,2n+1,r}(r) - \frac{i}{2} R_{2,2n+1,\varphi}(r)\bigr ]\;\sz{e}  ^{i (2n+1) \varphi}
\nonumber\\
\nonumber\\
\end{eqnarray}

\begin{eqnarray}\label{eqn477}
u_{22}^{*}(r,\varphi) = u_{2r2}^{*}(r,\varphi) + u_{2\varphi2}^{*}(r,\varphi)\;=
\quad\quad\quad\quad\quad\quad\quad\quad\quad\quad\quad\quad\quad\quad\quad
\nonumber\\
\nonumber\\
=\;i\;\bigl [ n R_{2,2n-1,r}(r) - \frac{i}{2} R_{2,2n-1,\varphi}(r)\bigr ]\;\sz{e}  ^{i (2n-1) \varphi} - i\;\bigl [ n R_{2,2n+1,r}(r) - \frac{i}{2} R_{2,2n+1,\varphi}(r)\bigr ]\;\sz{e}  ^{i (2n+1) \varphi}
\nonumber\\
\nonumber\\
\end{eqnarray}

From formula $(\ref{eqn410})$ for j = 2 and formula $(\ref{eqn457})$ we have:

\begin{equation}\label{eqn478}
u_{2t}(t) = \int_{0}^{t}T_{2n}(\tau) d\tau = \frac{1}{\sigma_{n}^{2}} \bigl [ t - \frac{2}{\sigma_{n}} \gamma (1, \sigma_{n}t) + \frac{1}{2 \sigma_{n}} \gamma (1,2 \sigma_{n}t) \bigr ]
\end{equation}

Hence, and from equation $(\ref{eqn411})$ it follows:

\begin{eqnarray}\label{eqn479}
u_{21}^{*}(r,\varphi,t)\;=\;u_{21}^{*}(r,\varphi)\;u_{2t}(t)
\nonumber\\
\nonumber\\
u_{22}^{*}(r,\varphi,t)\;=\;u_{22}^{*}(r,\varphi)\;u_{2t}(t)
\end{eqnarray}

After completing appropriate operations we have:

\begin{eqnarray}\label{eqn480}
u_{2 r}^{*}(r,\varphi, t)\;=\;\bigl \{\bigl [ n R_{2,2n-1,r}(r) - \frac{i}{2} R_{2,2n-1,\varphi}(r)\bigr ] + \bigl [ n R_{2,2n+1,r}(r) - \frac{i}{2} R_{2,2n+1,\varphi}(r)\bigr ]\bigr \}\;\sz{e}  ^{i 2n \varphi}\;u_{2t}(t)
\nonumber\\
\nonumber\\
\quad\quad\quad\quad
\end{eqnarray}

\begin{eqnarray}\label{eqn481}
u_{2\varphi}^{*}(r,\varphi, t)\;=\;i\;\bigl \{\bigl [ n R_{2,2n-1,r}(r) - \frac{i}{2} R_{2,2n-1,\varphi}(r)\bigr ] - \bigl [ n R_{2,2n+1,r}(r) - \frac{i}{2} R_{2,2n+1,\varphi}(r)\bigr ]\bigr \}\;\sz{e}  ^{i 2n \varphi}\;u_{2t}(t)
\nonumber\\
\nonumber\\
\quad\quad\quad\quad
\end{eqnarray}

Here $u_{2 r}^{*}(r,\varphi, t) ,\; u_{2\varphi}^{*}(r,\varphi, t)$ are the radial and tangential components of the first correction $\vec u_{2}^{*}$ of the velocity $\vec u_{1}$ and

\begin{eqnarray}\label{eqn482}
\bigl [ n R_{2,2n-1,r}(r) - \frac{i}{2} R_{2,2n-1,\varphi}(r)\bigr ]\;=\;\frac{F_{n}^{2}n^{2}r^{2n-1}}{2^{2}}\biggl [-\; \frac{\Gamma(1, 2\mu_{n}r)}{2\mu_{n}^{2}} -\;\frac{\Gamma(2, 2\mu_{n}r)}{2^{2}\mu_{n}^{2}} +
\nonumber\\
\nonumber\\
+\; n\sum_{l=0}^{\infty}\frac{\Gamma(l+2, 2\mu_{n}r)}{(2n+2)_{l+1}2^{l+1}\mu_{n}^{2}} - \sum_{l=0}^{\infty}\frac{\Gamma(l+3, 2\mu_{n}r)}{(2n+2)_{l+1}2^{l+3}\mu_{n}^{2}} \biggr ]
\nonumber\\
\quad\quad\quad\quad
\end{eqnarray}

\begin{eqnarray}\label{eqn483}
\bigl [ n R_{2,2n+1,r}(r) - \frac{i}{2} R_{2,2n+1,\varphi}(r)\bigr ]\;=\;\frac{F_{n}^{2}n^{2}}{2^{2}r^{2n+1}}\biggl [-\; \frac{n (2n+1)}{2^{4n-2}\mu_{n}^{4n+2}}\;\gamma(4n, 2\mu_{n}r) +
\nonumber\\
\nonumber\\
+\;\frac{2n (2n+1)}{\mu_{n}^{2n+2}}\;\gamma(2n, \mu_{n}r)\;r^{2n}\; \sz{e}  ^{- \mu_{n}r}\biggr ]
\nonumber\\
\quad\quad\quad\quad
\end{eqnarray}

From formulas $(\ref{eqn479})$ or $(\ref{eqn480}), (\ref{eqn481})$ with properties of $u_{2t}(t)\;$ - $\;(\ref{eqn478})$ it follows:

\begin{eqnarray}
  \begin{array}{ll} 								
lim \;\;u_{21}^{*}(r,\varphi,t)= 0;\;\;\;   lim \;\;u_{22}^{*}(r,\varphi,t)= 0;\\
t \rightarrow 0; \;\;\;\;\;\;\;\;\;\;\;\;\;\;\;\;\;\;\;\;\;\;\;\; t \rightarrow 0;
\end{array}
\end{eqnarray}

\begin{eqnarray}
  \begin{array}{ll} 								
lim \;\;u_{2r}^{*}(r,\varphi,t)= 0;\;\;\;   lim \;\;u_{2\varphi}^{*}(r,\varphi,t)= 0;\\
t \rightarrow 0; \;\;\;\;\;\;\;\;\;\;\;\;\;\;\;\;\;\;\;\;\;\;\;\; t \rightarrow 0;
\end{array}
\end{eqnarray}

and we have the velocity $\vec u_{2} = \vec u_{1} - \vec u_{2}^{*}$ [ look at $\;(\ref{eqn50})$] satisfying the initial conditions $(\ref{eqn400})$.
We use the asymptotic properties of the incomplete gamma functions $\Gamma (\alpha, x) ,\; \gamma(\alpha, x)$ and from formulas $(\ref{eqn482}), (\ref{eqn483})\;$ we have: the first correction $\vec u_{2}^{*}$ and therefore the velocity $\vec u_{2}$ satisfies conditions  $(\ref{eqn16})$ (for $r \;\rightarrow\; \infty $).

Let us compare the solution  $(\ref{eqn446})$ or $(\ref{eqn447}), (\ref{eqn448})$ for $\vec u_{1}$ of the first step of iterative process  with the first correction $(\ref{eqn479})$ or $(\ref{eqn480}), (\ref{eqn481})$ for $\vec u_{2}^{*}$  , which is received on the second step of iterative process. We see that

\begin{equation}\label{eqn484}
\mid\vec{u_{2}^{*}}\mid\; <<\; \mid\vec{u_{1}}\mid
\end{equation}

with conditions 

\begin{eqnarray}\label{eqn485}
F_{n}\;\leq\;\frac{1}{n} 
\nonumber\\
\nonumber\\
t\;\leq\;\sigma_{n}
\end{eqnarray}

By continuing this iterative process we can obtain next parts $\;\vec{u_{3}^{*}}\;, \vec{u_{4}^{*}}\;,...$ of the converging series for $\vec{u}$. For  arbitrary step j of the iterative process we have by using formula $(\ref{eqn57})$:
\begin{equation}\label{eqn486}
\vec{u}_{j}\;=\;\vec{u}_{1}\;-\;\sum_{l=2}^{j} \vec{u}_{l}^{*}\
\end{equation}
and then: 
\begin{equation}\label{eqn487}
  \begin{array}{ll} 								
lim \;\;\vec{u_{j}}= \vec{u}\\   
j \rightarrow \infty 
\end{array}
\end{equation}

where $\vec{u}$ is the solution of the problem $(\ref{eqn1}) - (\ref{eqn6})$ for $\nu$ = 0.
\nonumber\\
\nonumber\\
Below we provide numerical analysis of these results for the following values of problem's parameters:
 
Circumferential modes n = 1, 2, 3, 4, 5.

$\sigma_{n}$ = 10. 

${0\leq t \leq 10}$.

Results were obtained for functions $\vec u_{1} - (\ref{eqn446})$ or $(\ref{eqn447}), (\ref{eqn448})$ ; $\vec u_{2}^{*} - (\ref{eqn479})$ or $(\ref{eqn480}), (\ref{eqn481})$ with calculations of the incomplete gamma functions $\cite{BE253}$. 

$\vec u_{2} = \vec u_{1} - \vec u_{2}^{*}$ and is shown in FIG. 4.1 - 4.5.
The vector field $\vec u_{2}$ at distances r = 1, 2, 3, 5, 7 is represented by the dotted curves in left diagrams.
The comparison of $\mid\vec u_{1}\mid$ (dashed plots) and $\mid\vec u_{2}^{*}\mid$ (solid plots) in plane $\varphi$ = [0, $\pi$], at distances 0 $\leq$ r $\leq$ 50 is represented in right diagrams. This comparison shows $\mid\vec{u_{2}^{*}}\mid\; <<\; \mid\vec{u_{1}}\mid$ 
and is corresponding to the conclusion $(\ref{eqn484})$.
$\nonumber\\$
   \begin{center}    
     \includegraphics[height=60mm]{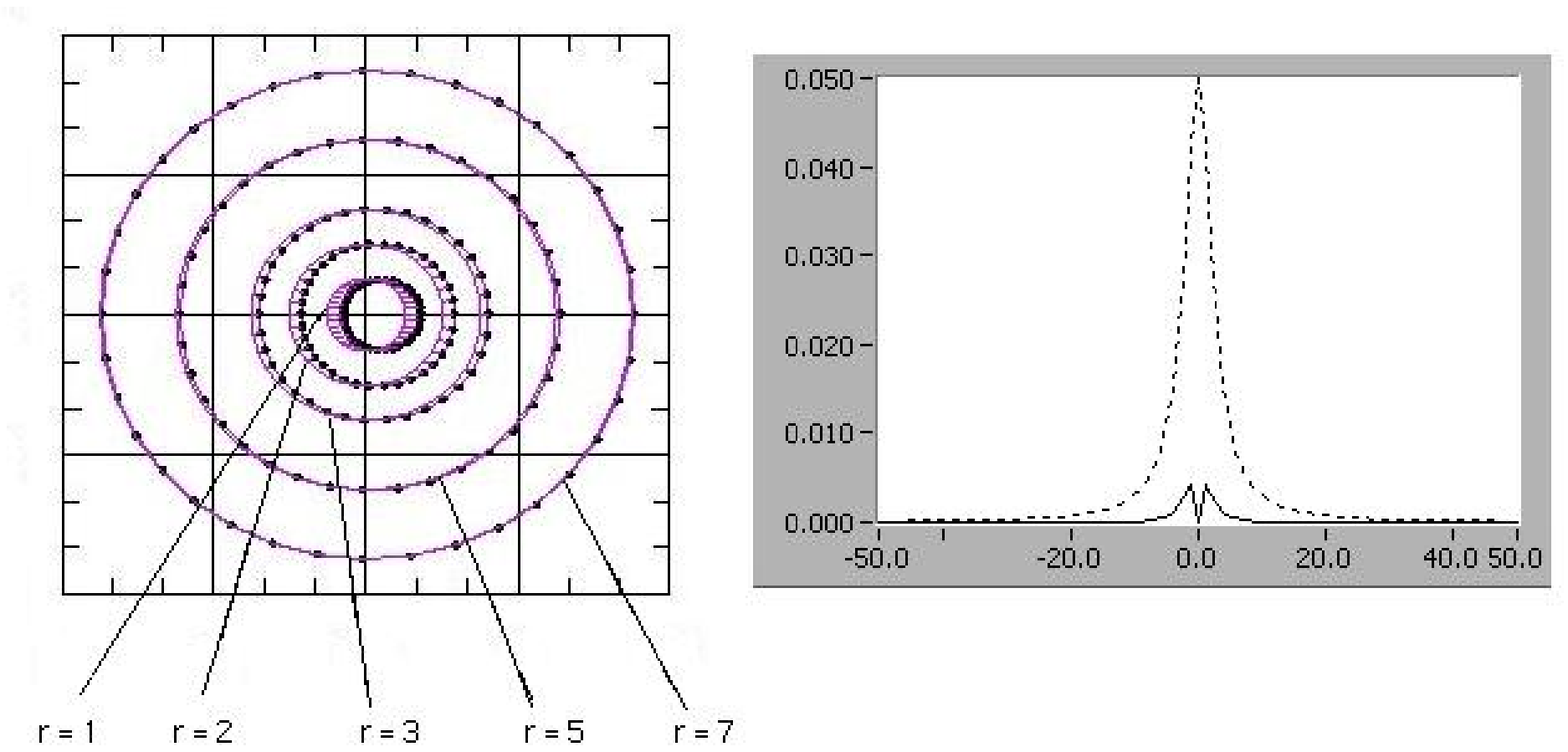}\\
FIG.4.1. n = 1, $F_1$ = 1, $\mu_1$ = 1
   \end{center}  
  \begin{center}  
    \includegraphics[height=60mm]{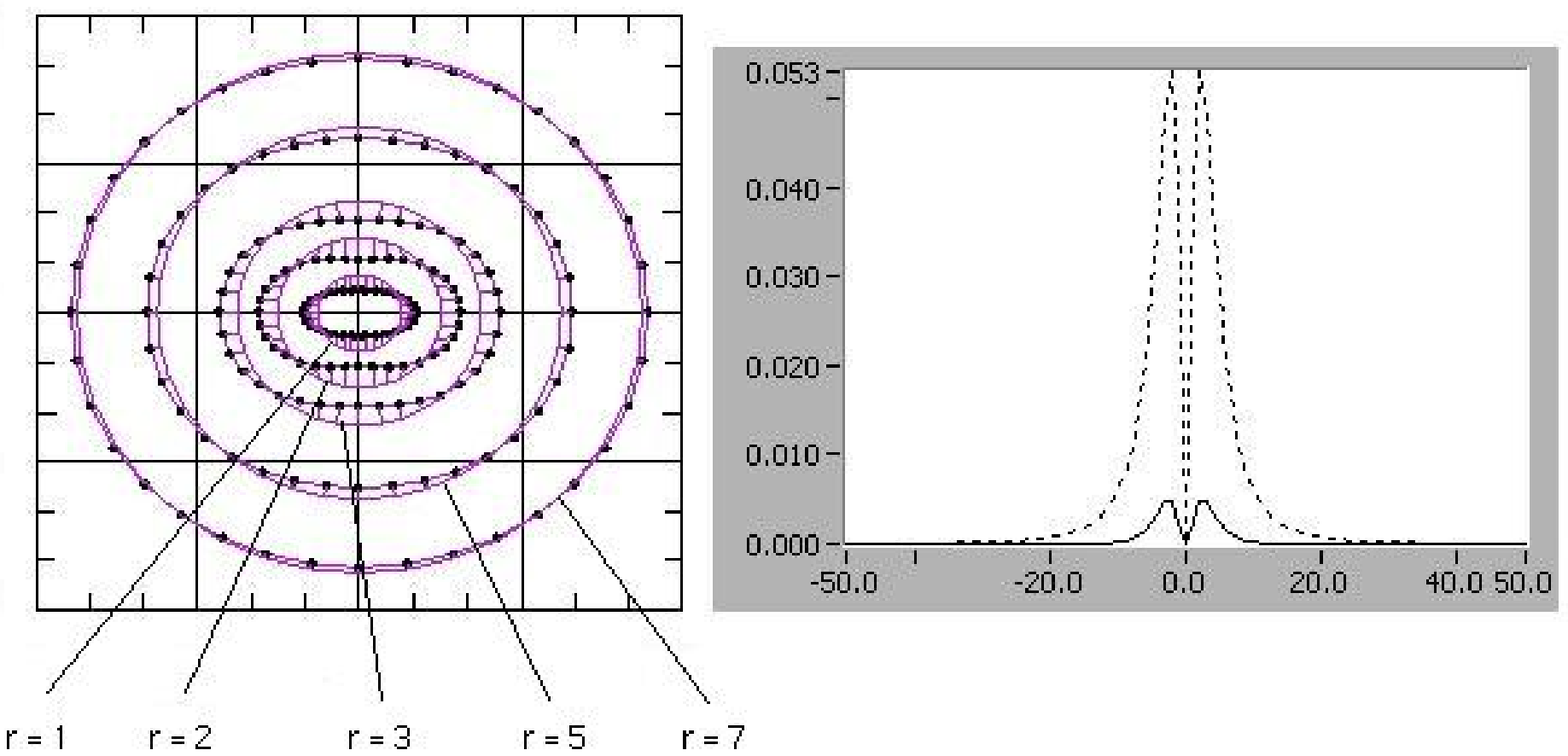}\\
     FIG.4.2. n = 2, $F_2$ = 0.5, $\mu_2$ = 1
  \end{center}    
    \begin{center}    
    \includegraphics[height=60mm]{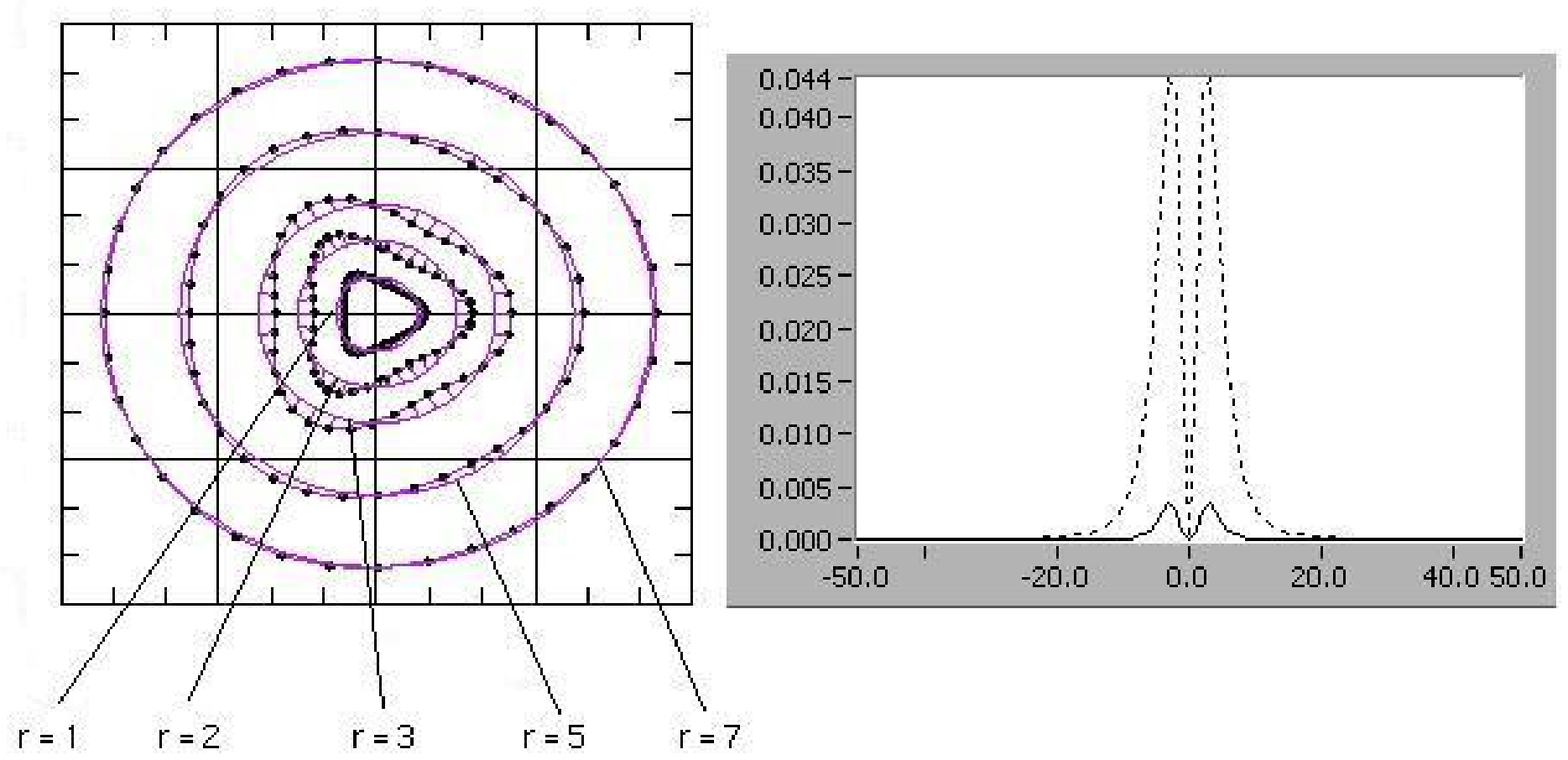}\\
    FIG.4.3. n = 3, $F_3$ = 0.33, $\mu_3$ = 1.3
      \end{center}  
 \begin{center} 
   \includegraphics[height=60mm]{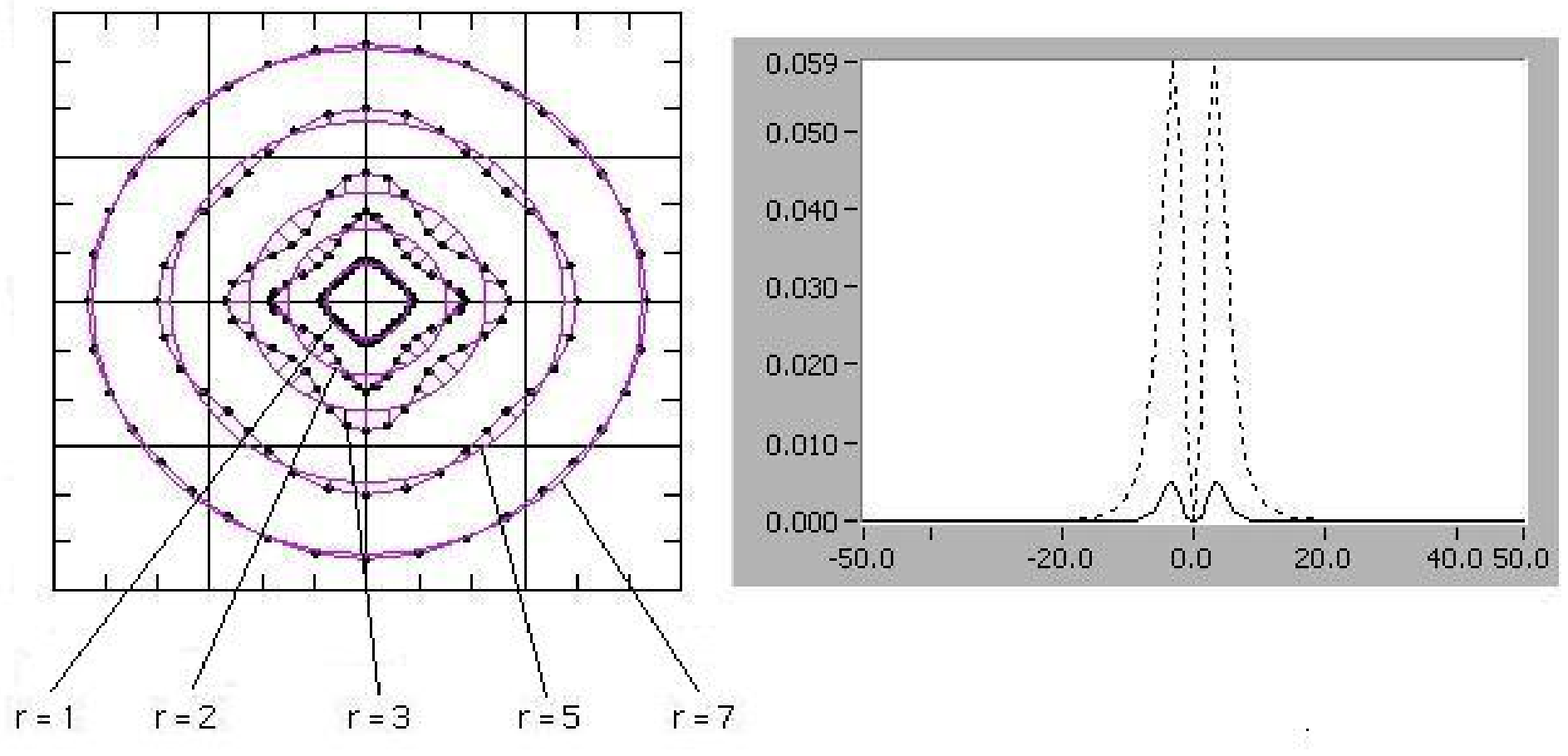}\\
    FIG.4.4. n = 4, $F_4$ = 0.25, $\mu_4$ = 1.5
  \end{center}    
 \begin{center}  
   \includegraphics[height=60mm]{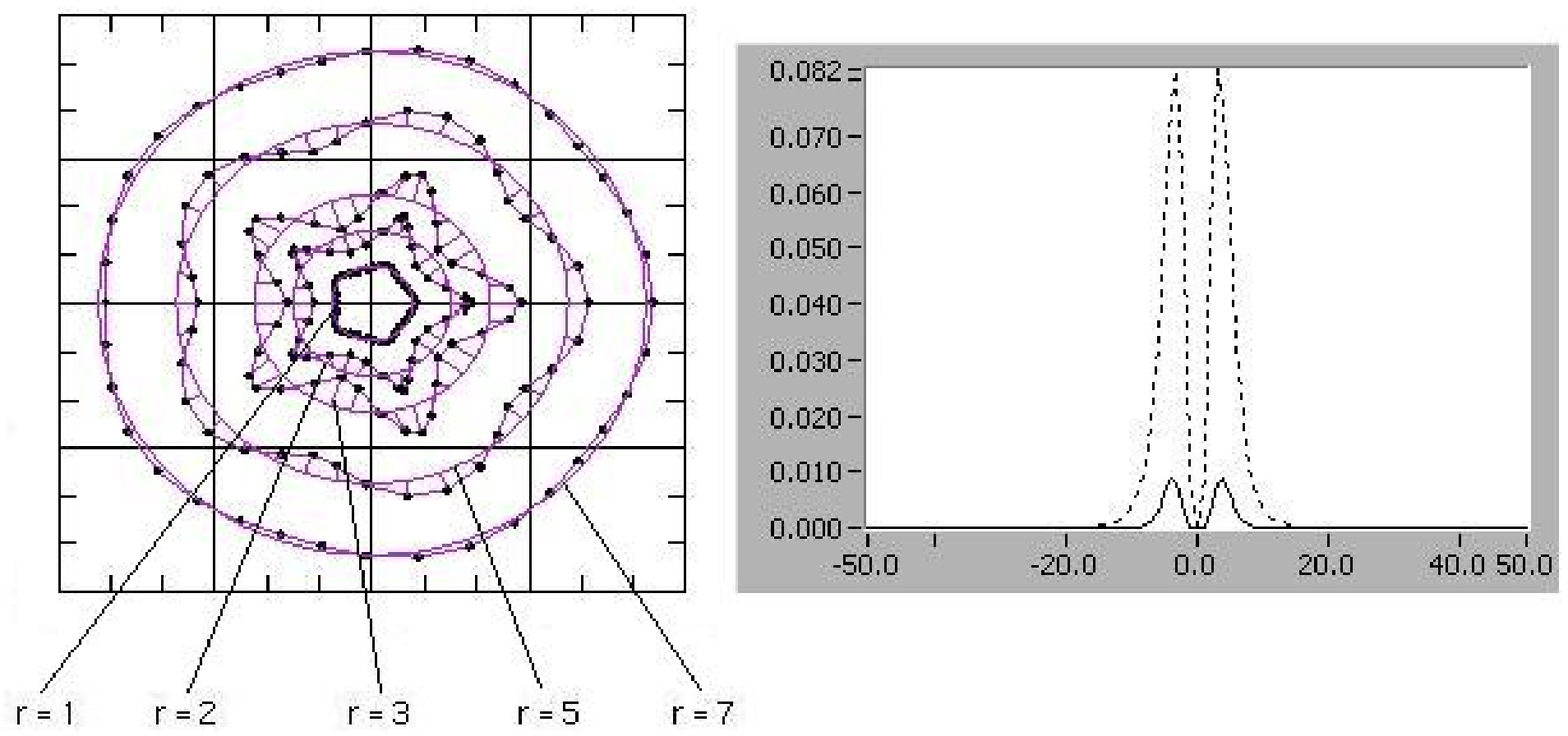}\\ 
   FIG.4.5. n = 5, $F_5$ = 0.2, $\mu_5$ = 1.7
 \end{center}  
 
$\nonumber\\
\nonumber\\$

\section{Example of the solution of the Cauchy problem  for the Navier - Stokes equations by the described iterative method with a particular applied force (N = 2)}\

We will consider an example of the solution of the Cauchy problem  for the Navier - Stokes equations for  N = 2 and with initial conditions:

\begin{equation}\label{eqn200}
\vec{u}(x,0)\; = \; \vec{u}^{0}(x)\; = \;0\;\;\;\;\;\;\;\; (x\in R^{2})
\end{equation}

Hence, and from formulas $(\ref{eqn39}), (\ref{eqn40})$ for arbitrary step j of the iterative process, it follows:

\begin{eqnarray}\label{eqn201}
u_{j1}(x_{1},x_{2},t)\;=\;\frac{1}{4\pi^2} \int_{-\infty}^{\infty}\int_{-\infty}^{\infty}\frac{\gamma_{2}^{2} }{(\gamma_{1}^{2}+\gamma_{2}^{2}) }   \int_{0}^{t} \sz{e}  ^{-\nu (\gamma_{1}^{2}+\gamma_{2}^{2}) (t-\tau)}\int_{-\infty}^{\infty}\int_{-\infty}^{\infty}\sz{e}  ^{i(\tilde x_{1}\gamma_{1}+\tilde x_{2}\gamma_{2})} f_{j1}(\tilde x_{1},\tilde x_{2},\tau)\,d \tilde x_{1}d \tilde  x_{2}d\tau\cdot
\nonumber\\
\nonumber\\
\nonumber\\
 \;\cdot\;\sz{e}  ^{-i(x_{1}\gamma_{1}+x_{2}\gamma_{2})}\,d\gamma_{1}d\gamma_{2}\;-
\quad\quad\quad\quad\quad\quad \quad\quad \quad\quad
\nonumber\\
\nonumber\\
\nonumber\\
- \;\frac{1}{4\pi^2} \int_{-\infty}^{\infty}\int_{-\infty}^{\infty}\frac{\gamma_{1} \gamma_{2} }{(\gamma_{1}^{2}+\gamma_{2}^{2}) }  \int_{0}^{t} \sz{e}  ^{-\nu (\gamma_{1}^{2}+\gamma_{2}^{2}) (t-\tau)}\int_{-\infty}^{\infty}\int_{-\infty}^{\infty}\sz{e}  ^{i(\tilde x_{1}\gamma_{1}+\tilde x_{2}\gamma_{2})} f_{j2}(\tilde x_{1},\tilde x_{2},\tau)\,d \tilde x_{1}d \tilde x_{2}d\tau\cdot
\nonumber\\
\nonumber\\
\nonumber\\
 \;\cdot\;\sz{e}  ^{-i(x_{1}\gamma_{1}+x_{2}\gamma_{2})}\,d\gamma_{1}d\gamma_{2}
\quad\quad\quad\quad\quad\quad \quad\quad \quad\quad
\nonumber\\
\quad\quad\quad\quad\quad\quad \quad\quad \quad\quad 
 \end{eqnarray}

\begin{eqnarray}\label{eqn202}
u_{j2}(x_{1},x_{2},t)\;=\;-\frac{1}{4\pi^2} \int_{-\infty}^{\infty}\int_{-\infty}^{\infty}\frac{\gamma_{1} \gamma_{2} }{(\gamma_{1}^{2}+\gamma_{2}^{2}) } \int_{0}^{t} \sz{e}  ^{-\nu (\gamma_{1}^{2}+\gamma_{2}^{2}) (t-\tau)}\int_{-\infty}^{\infty}\int_{-\infty}^{\infty}\sz{e}  ^{i(\tilde x_{1}\gamma_{1}+\tilde x_{2}\gamma_{2})} f_{j1}(\tilde x_{1},\tilde x_{2},\tau)\,d \tilde x_{1}d \tilde x_{2}d\tau\cdot
\nonumber\\
\nonumber\\
\nonumber\\
 \;\cdot\;\sz{e}  ^{-i(x_{1}\gamma_{1}+x_{2}\gamma_{2})}\,d\gamma_{1}d\gamma_{2}\;+
\quad\quad\quad\quad\quad\quad \quad\quad \quad\quad
\nonumber\\
\nonumber\\
\nonumber\\
+ \;\frac{1}{4\pi^2} \int_{-\infty}^{\infty}\int_{-\infty}^{\infty}\frac{\gamma_{1}^{2} }{(\gamma_{1}^{2}+\gamma_{2}^{2}) } \int_{0}^{t} \sz{e}  ^{-\nu (\gamma_{1}^{2}+\gamma_{2}^{2}) (t-\tau)}\int_{-\infty}^{\infty}\int_{-\infty}^{\infty}\sz{e}  ^{i(\tilde x_{1}\gamma_{1}+\tilde x_{2}\gamma_{2})} f_{j2}(\tilde x_{1},\tilde x_{2},\tau)\,d \tilde x_{1}d \tilde x_{2}d\tau \cdot
\nonumber\\
\nonumber\\
\nonumber\\
 \;\cdot\;\sz{e}  ^{-i(x_{1}\gamma_{1}+x_{2}\gamma_{2})}\,d\gamma_{1}d\gamma_{2}
\quad\quad\quad\quad\quad\quad \quad\quad \quad\quad
\nonumber\\
\quad\quad\quad\quad\quad\quad \quad\quad \quad\quad 
\end{eqnarray} 

We convert the Cartesian coordinates to the polar coordinates by formulas:

$x_{1}$ = $r\;\cdot\;$cos$\varphi\;$;$\;\;x_{2}$ = $r\;\cdot\;$sin$\varphi\;$;$\;\;\gamma_{1}$ = $\rho\;\cdot\;$cos$\psi\;$;$\;\;\gamma_{2}$ = $\rho\;\cdot\;$sin$\psi\;$;$\;\;\tilde x_{1}$ = $\tilde r\;\cdot\;$cos$\tilde \varphi\;$;$\;\;\tilde x_{2}$ = $\tilde r\;\cdot\;$sin$\tilde \varphi$;

and obtain from formulas $(\ref{eqn201}), (\ref{eqn202})$:

\begin{eqnarray}\label{eqn203}
u_{j1}(r,\varphi,t)\;=\;\frac{1}{4\pi^2} \int_{0}^{\infty}\int_{0}^{2\pi}{sin^{2}\psi } \int_{0}^{t} \sz{e}  ^{-\nu \rho^2 (t-\tau)}\int_{0}^{\infty}\int_{0}^{2\pi}\sz{e}  ^{i\tilde r\rho cos(\tilde \varphi-\psi)} f_{j1}(\tilde r,\tilde \varphi,\tau)\,\tilde r d \tilde r d \tilde  \varphi d\tau\cdot
\nonumber\\
\nonumber\\
\nonumber\\
 \;\cdot\;\sz{e}  ^{-ir\rho cos(\psi-\varphi)}\,\rho d\rho d\psi\;-
\quad\quad\quad\quad\quad\quad\quad
\nonumber\\
\nonumber\\
\nonumber\\
- \;\frac{1}{4\pi^2} \int_{0}^{\infty}\int_{0}^{2\pi}{sin\psi cos\psi }  \int_{0}^{t} \sz{e}  ^{-\nu \rho^2 (t-\tau)} \int_{0}^{\infty}\int_{0}^{2\pi}\sz{e}  ^{i\tilde r\rho cos(\tilde \varphi-\psi)} f_{j2}(\tilde r,\tilde \varphi,\tau)\,\tilde r d \tilde r d \tilde  \varphi d\tau\cdot
\nonumber\\
\nonumber\\
\nonumber\\
 \;\cdot\;
\sz{e}  ^{-ir\rho cos(\psi-\varphi)}\,\rho d\rho d\psi
\quad\quad\quad\quad\quad\quad\quad\quad\quad
\nonumber\\
\nonumber\\
\nonumber\\
\quad\quad\quad\quad\quad\quad\quad\quad\quad
\end{eqnarray}

\begin{eqnarray}\label{eqn204}
u_{j2}(r,\varphi,t)\;=\;-\frac{1}{4\pi^2} \int_{0}^{\infty}\int_{0}^{2\pi}{sin\psi cos\psi } \int_{0}^{t} \sz{e}  ^{-\nu \rho^2 (t-\tau)} \int_{0}^{\infty}\int_{0}^{2\pi}\sz{e}  ^{i\tilde r\rho cos(\tilde \varphi-\psi)} f_{j1}(\tilde r,\tilde \varphi,\tau)\,\tilde r d \tilde r d \tilde  \varphi d\tau\cdot
\nonumber\\
\nonumber\\
\nonumber\\
 \;\;\;\;\cdot\;\sz{e}  ^{-ir\rho cos(\psi-\varphi)}\,\rho d\rho d\psi\;+
\quad\quad\quad\quad\quad\quad \quad\quad \quad\quad
\nonumber\\
\nonumber\\
\nonumber\\
+ \;\frac{1}{4\pi^2} \int_{0}^{\infty}\int_{0}^{2\pi}cos^{2}\psi \int_{0}^{t} \sz{e}  ^{-\nu \rho^2 (t-\tau)}\int_{0}^{\infty}\int_{0}^{2\pi}\sz{e}  ^{i\tilde r\rho cos(\tilde \varphi-\psi)} f_{j2}(\tilde r,\tilde \varphi,\tau)\,\tilde r d \tilde r d \tilde  \varphi d\tau\cdot
\nonumber\\
\nonumber\\
\nonumber\\
 \;\;\;\;\cdot\;\sz{e}  ^{-ir\rho cos(\psi-\varphi)}\,\rho d\rho d\psi
\quad\quad\quad\quad\quad\quad \quad\quad \quad\quad
\nonumber\\
\quad\quad\quad\quad\quad\quad \quad\quad \quad\quad 
\end{eqnarray} 

We have the applied force $\vec{f}_{j}$ for arbitrary step j of the iterative process:

\begin{equation}\label{eqn205}
f_{j\tilde r}(\tilde r,\tilde \varphi,\tau) = f_{j\tilde r}(\tilde r,\tau)\sz{e}^{in_{j}\tilde \varphi}\;\;\;,\;\;\;f_{j\tilde \varphi}(\tilde r,\tilde \varphi,\tau) \equiv 0
\end{equation}

or

\begin{equation}\label{eqn206}
f_{j\tilde r}(\tilde r,\tilde \varphi,\tau)\equiv 0\;\;\;,\;\;\;f_{j\tilde \varphi}(\tilde r,\tilde \varphi,\tau) = f_{j\tilde \varphi}(\tilde r,\tau)\sz{e}^{in_{j}\tilde \varphi}
\end{equation}

where $\;\;f_{j\tilde r}(\tilde r,\tilde \varphi,\tau)\;\;\;,\;\;\;f_{j\tilde \varphi}(\tilde r,\tilde \varphi,\tau)\;\;\; - \;\;\;$  radial and tangential components of  the applied force. 

$n_{j}$ - separate circumferential mode, $n_{j}$ = 0,1,2,3,...

We take the radial and tangential components of  the applied force $(\ref{eqn205}), (\ref{eqn206})$ with condition $(\ref{eqn18})\;$. For the radial component of  the applied force we use De Moivre's formulas (\ref{A8}) and have:

\begin{eqnarray}\label{eqn207}
f_{j1}(\tilde r,\tilde \varphi,\tau) = f_{j\tilde r}(\tilde r,\tau)\sz{e}^{in_{j}\tilde \varphi}cos \tilde \varphi  = \frac{1}{2}f_{j\tilde r}(\tilde r,\tau)\bigl(\sz{e}^{i(n_{j}-1)\tilde \varphi} + \sz{e}^{i(n_{j}+1)\tilde \varphi}\bigr)
\nonumber\\
\nonumber\\
f_{j2}(\tilde r,\tilde \varphi,\tau) = f_{j\tilde r}(\tilde r,\tau)\sz{e}^{in_{j}\tilde \varphi}sin \tilde \varphi  = \frac{i}{2}f_{j\tilde r}(\tilde r,\tau)\bigl(\sz{e}^{i(n_{j}-1)\tilde \varphi} - \sz{e}^{i(n_{j}+1)\tilde \varphi}\bigr)
\end{eqnarray}

We put the applied force components $(\ref{eqn207})$ in formulas $(\ref{eqn203}), (\ref{eqn204})$ and find:

\begin{eqnarray}\label{eqn208}
u_{jr1}(r,\varphi,t)\;=\;\frac{1}{8\pi^2}  \int_{0}^{\infty}\int_{0}^{2\pi}{sin^{2}\psi } \biggl[ \int_{0}^{t} \sz{e}  ^{-\nu \rho^2 (t-\tau)}\int_{0}^{\infty} f_{j\tilde r}(\tilde r, \tau) \cdot
\quad\quad\quad\quad\quad\quad\quad\quad\quad\quad
\nonumber\\
\nonumber\\
\nonumber\\
\cdot\;\int_{0}^{2\pi}\sz{e}  ^{i\tilde r\rho cos(\tilde \varphi-\psi)} (\sz{e}^{i(n_{j}-1)\tilde \varphi} + \sz{e}^{i(n_{j}+1)\tilde \varphi})d \tilde  \varphi \tilde r d \tilde r  d\tau\biggr] \sz{e}  ^{-ir\rho cos(\psi-\varphi)}\,\rho d\rho d\psi\;-
\quad\quad\quad\quad\quad\quad\quad
\nonumber\\
\nonumber\\
\nonumber\\
-\;\; \frac{i}{8\pi^2} \int_{0}^{\infty}\int_{0}^{2\pi}{sin\psi cos\psi } \biggl[ \int_{0}^{t} \sz{e}  ^{-\nu \rho^2 (t-\tau)} \int_{0}^{\infty} f_{j\tilde r}(\tilde r, \tau) \cdot
\quad\quad\quad\quad\quad\quad\quad\quad\quad
\nonumber\\
\nonumber\\
\nonumber\\
 \;\cdot\;\int_{0}^{2\pi}\sz{e}  ^{i\tilde r\rho cos(\tilde \varphi-\psi)}(\sz{e}^{i(n_{j}-1)\tilde \varphi} - \sz{e}^{i(n_{j}+1)\tilde \varphi}) d \tilde  \varphi \tilde r d \tilde r  d\tau \biggr] \sz{e}  ^{-ir\rho cos(\psi-\varphi)}\,\rho d\rho d\psi 
\quad\quad\quad\quad\quad\quad\quad\quad\quad
\nonumber\\
\nonumber\\
\nonumber\\
\quad\quad\quad\quad\quad\quad\quad\quad\quad
\end{eqnarray}

\begin{eqnarray}\label{eqn209}
u_{jr2}(r,\varphi,t)\;=\;- \frac{1}{8\pi^2}  \int_{0}^{\infty}\int_{0}^{2\pi}{sin\psi cos\psi} \biggl[ \int_{0}^{t}\sz{e}  ^{-\nu \rho^2 (t-\tau)}\int_{0}^{\infty} f_{j\tilde r}(\tilde r, \tau) \cdot
\quad\quad\quad\quad\quad\quad\quad\quad\quad
\nonumber\\
\nonumber\\
\nonumber\\
\cdot\;\int_{0}^{2\pi}\sz{e}  ^{i\tilde r\rho cos(\tilde \varphi-\psi)} (\sz{e}^{i(n_{j}-1)\tilde \varphi} + \sz{e}^{i(n_{j}+1)\tilde \varphi})d \tilde  \varphi \tilde r d \tilde r  d\tau \biggr] \sz{e}  ^{-ir\rho cos(\psi-\varphi)}\,\rho d\rho d\psi\;+
\quad\quad\quad\quad\quad\quad\quad
\nonumber\\
\nonumber\\
\nonumber\\
+\;\; \frac{i}{8\pi^2} \int_{0}^{\infty}\int_{0}^{2\pi}{ cos^{2}\psi } \biggl[ \int_{0}^{t}  \sz{e}  ^{-\nu \rho^2 (t-\tau)}\int_{0}^{\infty} f_{j\tilde r}(\tilde r, \tau)  \cdot
\quad\quad\quad\quad\quad\quad\quad\quad\quad
\nonumber\\
\nonumber\\ 
\nonumber\\
 \;\cdot\;\int_{0}^{2\pi}\sz{e}  ^{i\tilde r\rho cos(\tilde \varphi-\psi)}(\sz{e}^{i(n_{j}-1)\tilde \varphi} - \sz{e}^{i(n_{j}+1)\tilde \varphi}) d \tilde  \varphi \tilde r d \tilde r  d\tau \biggr] \sz{e}  ^{-ir\rho cos(\psi-\varphi)}\,\rho d\rho d\psi  
\quad\quad\quad\quad\quad\quad\quad\quad\quad
\nonumber\\
\nonumber\\
\nonumber\\
\quad\quad\quad\quad\quad\quad\quad\quad\quad
\end{eqnarray}

Let us denote  internal integrals in $(\ref{eqn208}), (\ref{eqn209})$ as $\;\;I_{\underline{+}}(\tilde r,\rho,\psi)\;\;$: 

\begin{equation}\label{eqn210}
I_{\underline{+}}(\tilde r,\rho,\psi) = \int_{0}^{2\pi}\sz{e}  ^{i\tilde r\rho cos(\tilde \varphi-\psi)} (\sz{e}^{i(n_{j}-1)\tilde \varphi\;\;} \underline{+}\;\; \sz{e}^{i(n_{j}+1)\tilde \varphi})d \tilde  \varphi
\end{equation}

 We have two integrals here. Plus (+) is for the first part of each integral $(\ref{eqn208}), (\ref{eqn209})$ and minus (-) is for the second part.

We substitute $\tilde \theta$ for $\tilde \varphi$:$\;\;\;\tilde \theta\;$ = $\;\tilde \varphi\;$ - $\psi\;\;$ , d$\tilde \theta\;$ = d$\tilde \varphi\;$ and receive:

\begin{equation}\label{eqn211}
I_{\underline{+}}(\tilde r,\rho,\psi) = \sz{e}^{i(n_{j}-1)\psi}\int_{-\psi}^{2\pi-\psi}\sz{e}  ^{i\tilde r\rho cos\tilde \theta + i(n_{j}-1)\tilde \theta} d \tilde \theta \;\;\underline{+}\;\;\sz{e}^{i(n_{j}+1)\psi}\int_{-\psi}^{2\pi-\psi}\sz{e}  ^{i\tilde r\rho cos\tilde \theta + i(n_{j}+1)\tilde \theta} d \tilde \theta
\end{equation}

Then we use the Bessel function's integral representation $(\ref{A9})$ and have:

\begin{equation}\label{eqn212}
I_{\underline{+}}(\tilde r,\rho,\psi) = 2\pi i^{(n_{j}-1)}\sz{e}^{i(n_{j}-1)\psi}J_{n_{j}-1}(\tilde{r}\rho)\;\;\underline{+}\;\;2\pi i^{(n_{j}+1)}\sz{e}^{i(n_{j}+1)\psi}J_{n_{j}+1}(\tilde{r}\rho)
\end{equation}

Let us put $\;\;I_{\underline{+}}(\tilde r,\rho,\psi)\;\;$ from $(\ref{eqn212})$ in formulas $(\ref{eqn208}), (\ref{eqn209})\;$ , change order of integration and obtain:

\begin{eqnarray}\label{eqn213}
u_{jr1}(r,\varphi, t)\;=\;\frac{1}{8\pi^2}\int_{0}^{\infty}\int_{0}^{t}  \sz{e}  ^{-\nu \rho^2 (t-\tau)} \int_{0}^{\infty}f_{j\tilde r}(\tilde r, \tau) \cdot
\quad\quad\quad\quad\quad\quad\quad\quad\quad\quad\quad\quad\quad\quad\quad\quad\quad\quad
\nonumber\\
\nonumber\\
\cdot \int_{0}^{2\pi}\biggl[{sin^{2} \psi} I_{+}(\tilde r,\rho,\psi)-i\;{sin\psi cos\psi}I_{-}(\tilde r,\rho,\psi)\biggr]\;\sz{e}  ^{-ir\rho cos(\psi-\varphi)}d\psi \tilde r d \tilde r d\tau \rho d\rho 
\nonumber\\
\nonumber\\
\quad\quad\quad\quad\quad\quad\quad\quad\quad
\end{eqnarray}

\begin{eqnarray}\label{eqn214}
u_{jr2}(r,\varphi, t)\;=\;\frac{1}{8\pi^2}\int_{0}^{\infty}\int_{0}^{t}  \sz{e}  ^{-\nu \rho^2 (t-\tau)} \int_{0}^{\infty}f_{j\tilde r}(\tilde r, \tau) \cdot
\quad\quad\quad\quad\quad\quad\quad\quad\quad\quad\quad\quad\quad\quad\quad\quad\quad\quad
\nonumber\\
\nonumber\\
\cdot \int_{0}^{2\pi}\biggl[-{sin\psi}{cos\psi} I_{+}(\tilde r,\rho,\psi)+i\;{cos^{2}\psi}I_{-}(\tilde r,\rho,\psi)\biggr]\;\sz{e}  ^{-ir\rho cos(\psi-\varphi)}d\psi \tilde r d \tilde r d\tau \rho d\rho 
\nonumber\\
\nonumber\\
\quad\quad\quad\quad\quad\quad\quad\quad\quad
\end{eqnarray}

Then we group parts in brackets of formulas $(\ref{eqn213}), (\ref{eqn214})\;$ , use De Moivre's formulas $(\ref{A8})$ and the Bessel function's properties. And we get:

\begin{equation}\label{eqn215}
u_{jr1}(r,\varphi, t)\;=\;-\;\frac{n_{j}\;i^{n_{j}}}{2\pi}\int_{0}^{\infty}\int_{0}^{t}  \sz{e}  ^{-\nu \rho^2 (t-\tau)} \int_{0}^{\infty}f_{j\tilde r}(\tilde r, \tau)\int_{0}^{2\pi}{sin\psi}\;\sz{e}  ^{-ir\rho cos(\psi-\varphi)+i n_{j} \psi}\;d\psi \;J_{n_{j}}(\tilde r \rho)\;d \tilde r d\tau d\rho 
\quad\quad
\end{equation}

\begin{equation}\label{eqn216}
u_{jr2}(r,\varphi, t)\;=\;\frac{n_{j}\;i^{n_{j}}}{2\pi}\int_{0}^{\infty}\int_{0}^{t}  \sz{e}  ^{-\nu \rho^2 (t-\tau)} \int_{0}^{\infty}f_{j\tilde r}(\tilde r, \tau)\int_{0}^{2\pi}{cos\psi}\;\sz{e}  ^{-ir\rho cos(\psi-\varphi)+i n_{j} \psi}\;d\psi \;J_{n_{j}}(\tilde r \rho)\;d \tilde r d\tau d\rho 
\quad\quad
\end{equation}
\nonumber\\

We substitute $\theta$ for $\psi$:$\;\;\;\theta\;$ = $\;\psi\;$ - $\varphi\;\;$ , d$\theta\;$ = d$\psi\;$ in the internal integrals of formulas $(\ref{eqn215}), (\ref{eqn216})\;$, use De Moivre's formulas $(\ref{A8})$ and the Bessel function's integral representation $(\ref{A9})$ and have from formulas $(\ref{eqn215}), (\ref{eqn216})$:

\begin{equation}\label{eqn217}
u_{jr1}(r,\varphi, t)\;=\;\frac{n_{j}}{2}\;\sz{e}  ^{i n_{j} \varphi}\int_{0}^{\infty}\int_{0}^{t}  \sz{e}  ^{-\nu \rho^2 (t-\tau)} \biggl[\sz{e}  ^{i \varphi}J_{n_{j+1}}(r \rho)\;+\;\sz{e}  ^{-i \varphi}J_{n_{j-1}}(r \rho)\biggr]\int_{0}^{\infty}f_{j\tilde r}(\tilde r, \tau)\;J_{n_{j}}(\tilde r \rho)\;d \tilde r d\tau d\rho 
\quad\quad
\end{equation}

\begin{equation}\label{eqn218}
u_{jr2}(r,\varphi, t)\;=\;\frac{i\; n_{j}}{2}\;\sz{e}  ^{i n_{j} \varphi}\int_{0}^{\infty}\int_{0}^{t}  \sz{e}  ^{-\nu \rho^2 (t-\tau)} \biggl[\sz{e}  ^{i \varphi}J_{n_{j+1}}(r \rho)\;-\;\sz{e}  ^{-i \varphi}J_{n_{j-1}}(r \rho)\biggr]\int_{0}^{\infty}f_{j\tilde r}(\tilde r, \tau)\;J_{n_{j}}(\tilde r \rho)\;d \tilde r d\tau d\rho 
\quad\quad
\end{equation}
\nonumber\\

Let us denote:

\begin{equation}\label{eqn219}
R_{j,n_{j}-1,r}(r, t)\;=\;\int_{0}^{\infty}\int_{0}^{t}  \sz{e}  ^{-\nu \rho^2 (t-\tau)} J_{n_{j}-1}(r \rho)\int_{0}^{\infty}f_{j\tilde r}(\tilde r, \tau)\;J_{n_{j}}(\tilde r \rho)\;d \tilde r d\tau d\rho 
\quad\quad
\end{equation}

\begin{equation}\label{eqn220}
R_{j,n_{j}+1,r}(r, t)\;=\;\int_{0}^{\infty}\int_{0}^{t}  \sz{e}  ^{-\nu \rho^2 (t-\tau)} J_{n_{j}+1}(r \rho)\int_{0}^{\infty}f_{j\tilde r}(\tilde r, \tau)\;J_{n_{j}}(\tilde r \rho)\;d \tilde r d\tau d\rho 
\quad\quad
\end{equation}

Then we have from formulas $(\ref{eqn217}), (\ref{eqn218})$:

\begin{equation}\label{eqn221}
u_{jr1}(r,\varphi, t)\;=\;\frac{n_{j}}{2}\bigl[R_{j,n_{j}-1,r}(r, t)\;\sz{e}  ^{i (n_{j}-1) \varphi}\;+\;R_{j,n_{j}+1,r}(r, t)\;\sz{e}  ^{i (n_{j}+1) \varphi}\bigr]
\end{equation}

\begin{equation}\label{eqn222}
u_{jr2}(r,\varphi, t)\;=\;\frac{i\;n_{j}}{2}\bigl[R_{j,n_{j}-1,r}(r, t)\;\sz{e}  ^{i (n_{j}-1) \varphi}\;-\;R_{j,n_{j}+1,r}(r, t)\;\sz{e}  ^{i (n_{j}+1) \varphi}\bigr]
\end{equation}

Then if $n_{j}\;=\;0$ it follows from $(\ref{eqn219}), (\ref{eqn220}) , (\ref{eqn221}), (\ref{eqn222})$ that
$u_{jr1}(r,\varphi, t)\;=\;u_{jr2}(r,\varphi, t)\;=\;0\;$ and hence 
$u_{1}\;\;=\;\;u_{2}\;\;=\;\;0$.

In the equations bellow we will consider $n_{j}\;\geq\;1$.

Now we integrate the solution $(\ref{eqn203}), (\ref{eqn204})$ by the tangential component of the applied force $(\ref{eqn206})$ ( for $n_{j}\;\geq\;1$). Then we use De Moivre's formulas $(\ref{A8})$ and have:

\begin{eqnarray}\label{eqn223}
f_{j1}(\tilde r,\tilde \varphi,\tau) = - f_{j\tilde \varphi}(\tilde r, \tau)\sz{e}^{in_{j}\tilde \varphi}sin \tilde \varphi  = - \frac{i}{2}f_{j\tilde \varphi}(\tilde r, \tau) \bigl( \sz{e}^{i(n_{j}-1)\tilde \varphi} - \sz{e}^{i(n_{j}+1)\tilde \varphi} \bigr) 
\nonumber\\
\nonumber\\
f_{j2}(\tilde r,\tilde \varphi,\tau) \;\;= \;\; f_{j\tilde \varphi}(\tilde r, \tau)\sz{e}^{in_{j}\tilde \varphi}cos \tilde \varphi \;= \; \frac{1}{2}f_{j\tilde \varphi}(\tilde r, \tau)\bigl(\sz{e}^{i(n_{j}-1)\tilde \varphi} + \sz{e}^{i(n_{j}+1)\tilde \varphi}\bigr)
\nonumber\\
\nonumber\\
\end{eqnarray}

Hence formulas $(\ref{eqn223})$ are the components $f_{j1}$ and $f_{j2}$ from the tangential applied force $(\ref{eqn206})$, while formulas $(\ref{eqn207})$ are the components $f_{j1}$ and $f_{j2}$ from the radial applied force $(\ref{eqn205})$.

Let us put $(\ref{eqn223})$ in formulas $(\ref{eqn203}), (\ref{eqn204})$ and do the operations as we did in $(\ref{eqn208}) - (\ref{eqn222})\;  (\;n_{j}\;\geq\;1)$. We consider that $f_{j\tilde \varphi}(\tilde r, \tau)$ is restricted by condition $(\ref{eqn18})$ and get:

\begin{equation}\label{eqn224}
R_{j,n_{j}-1,\varphi}(r, t)\;=\;-\;\int_{0}^{\infty}\int_{0}^{t}  \sz{e}  ^{-\nu \rho^2 (t-\tau)}J_{n_{j}-1}(r \rho)\int_{0}^{\infty}\bigl(f_{j\tilde \varphi}(\tilde r, \tau)\cdot\tilde r\bigr)^{'}_{\tilde r}\;J_{n_{j}}(\tilde r \rho)\;d \tilde r d\tau d\rho 
\quad\quad
\end{equation}

\begin{equation}\label{eqn225}
R_{j,n_{j}+1,\varphi}(r, t)\;=\;-\;\int_{0}^{\infty}\int_{0}^{t}  \sz{e}  ^{-\nu \rho^2 (t-\tau)}J_{n_{j}+1}(r \rho)\int_{0}^{\infty}\bigl(f_{j\tilde \varphi}(\tilde r, \tau)\cdot\tilde r\bigr)^{'}_{\tilde r}\;J_{n_{j}}(\tilde r \rho)\;d \tilde r d\tau d\rho \;,
\quad\quad
\end{equation}

Here $\bigl(\bigr)^{'}_{\tilde r} \equiv \frac{\partial}{\partial \tilde r}$. Hence we have:

\begin{equation}\label{eqn226}
u_{j\varphi1}(r,\varphi, t)\;=\;-\;\frac{i}{2}\bigl[R_{j,n_{j}-1,\varphi}(r, t)\;\sz{e}  ^{i (n_{j}-1) \varphi}\;+\;R_{j,n_{j}+1,\varphi}(r, t)\;\sz{e}  ^{i (n_{j}+1) \varphi}\bigr]
\end{equation}

\begin{equation}\label{eqn227}
u_{j\varphi2}(r,\varphi, t)\;=\;\frac{1}{2}\bigl[R_{j,n_{j}-1,\varphi}(r, t)\;\sz{e}  ^{i (n_{j}-1) \varphi}\;-\;R_{j,n_{j}+1,\varphi}(r, t)\;\sz{e}  ^{i (n_{j}+1) \varphi}\bigr]
\end{equation}

We have obtain formulas $(\ref{eqn207})\;$ - $\;(\ref{eqn227})$ for an arbitrary step j of the iterative process and the applied forces $(\ref{eqn205})$ or $(\ref{eqn206})$. 

Now we investigate the first step (j = 1 , $n_{1} = n = 1,2,3,...$) of the iterative process with the particular radial applied force ${f}_{1}(x,t)\;\;$[look at $\;(\ref{eqn47})$]:

\begin{eqnarray}\label{eqn228}
f_{1\tilde r}(\tilde r,\tilde \varphi,\tau) = f_{1\tilde r}(\tilde r, \tau)\sz{e}^{in\tilde \varphi}\;\;\;,\;\;\;f_{1\tilde \varphi}(\tilde r,\tilde \varphi,\tau) \equiv 0
\nonumber\\
\nonumber\\
f_{1\tilde r}(\tilde r, \tau) = F_{n}\tilde r^{n+1}\sz{e}^{-\mu_{n}^2\tilde r^2}f_{1\tau}(\tau)
\quad\quad
\end{eqnarray}

$F_{n}, \mu_{n}$ - constants.   , $0 < F_{n}< \infty$ , $1 < \mu_{n}< \infty$.

Let us put the particular radial applied force $(\ref{eqn228})$ in formulas $(\ref{eqn219})$ , $(\ref{eqn220})$ and for the internal integral we have by using formula $ (\ref{A10}) , \cite{BE253}$:

\begin{equation}\label{eqn229}
I(\rho, \tau) = \int_{0}^{\infty}f_{1\tilde r}(\tilde r, \tau)\;J_{n}(\tilde r \rho)\;d \tilde r = F_{n} f_{1\tau}(\tau) \int_{0}^{\infty} \tilde r^{n+1}\sz{e}^{-\mu_{n}^2\tilde r^2}J_{n}(\tilde r \rho)\;d \tilde r = \frac{F_{n} f_{1\tau}(\tau)\rho^n}{(2\mu_{n}^2)^{n+1}}\sz{e}^{-\frac{\rho^2}{4\mu_{n}^2}}
\end{equation}

Now we put $I(\rho, \tau)$ from formula $(\ref{eqn229})$ in formulas $(\ref{eqn219})$ , $(\ref{eqn220})$ , change the order of integration and have by using formula $ (\ref{A11}) , \cite{BE253}$:

\begin{eqnarray}\label{eqn230}
R_{1,n-1,r}(r, t)\;=\;F_{n}  \int_{0}^{t}f_{1\tau}(\tau) \int_{0}^{\infty}  \sz{e}  ^{-\bigl [\nu  (t-\tau)+\frac{1}{4\mu_{n}^2}\bigr ]\rho^2}\frac{\rho^n}{(2\mu_{n}^2)^{n+1}} J_{n-1}(r \rho)  d\rho d\tau \;=
\nonumber\\
\nonumber\\
=\;\frac{F_{n} r^{n-1}}{2\mu_{n}^2}\int_{0}^{t}\frac{f_{1\tau}(\tau)\Phi\bigl(n+2, n+2;\frac{-\mu_{n}^2 r^2}{\bigl[4\mu_{n}^2\nu (t-\tau)+ 1\bigr]}\bigr)}{[4\mu_{n}^2\nu (t-\tau)+ 1]^n}d\tau
\quad\quad\quad\quad\quad\quad
\end{eqnarray}

\begin{eqnarray}\label{eqn231}
R_{1,n+1,r}(r, t)\;=\;F_{n} \int_{0}^{t}f_{1\tau}(\tau) \int_{0}^{\infty}  \sz{e}  ^{-\bigl [\nu  (t-\tau)+\frac{1}{4\mu_{n}^2}\bigr ]\rho^2}\frac{\rho^n}{(2\mu_{n}^2)^{n+1}} J_{n+1}(r \rho)  d\rho d\tau\;=
\nonumber\\
\nonumber\\
=\;\frac{F_{n} r^{n+1}}{2(n+1)}\int_{0}^{t}\frac{f_{1\tau}(\tau)\Phi\bigl(n+1, n+2;\frac{-\mu_{n}^2 r^2}{\bigl[4\mu_{n}^2\nu (t-\tau)+ 1\bigr]}\bigr)}{[4\mu_{n}^2\nu (t-\tau)+ 1]^{n+1}}d\tau
\quad\quad\quad\quad\quad\quad
\end{eqnarray}

Here $\Phi(a, c; x)$ is a confluent hypergeometric function $\cite{BE153}$.

We substitute y for $\tau$: y = $\frac{1}{[4\mu_{n}^2\nu (t-\tau)+ 1]}$, dy = $\frac{4\mu_{n}^2\nu}{[4\mu_{n}^2\nu (t-\tau)+ 1]^2} d\tau$ and receive:

\begin{eqnarray}\label{eqn232}
R_{1,n-1,r}(r, t)\;=\;\frac{F_{n} r^{n-1}}{8\mu_{n}^4\nu}\int_{\frac{1}{[4\mu_{n}^2\nu t+1]}}^{1}f_{1\tau}(y)\cdot y^{n-2}\cdot \Phi\bigl(n+2, n+2;-\mu_{n}^2 r^2y\bigr)dy
\quad\quad\quad\quad\quad\quad
\end{eqnarray}

\begin{eqnarray}\label{eqn233}
R_{1,n+1,r}(r, t)\;=\frac{F_{n} r^{n+1}}{8\mu_{n}^2\nu(n+1)}\int_{\frac{1}{[4\mu_{n}^2\nu t+1]}}^{1}f_{1\tau}(y)\cdot y^{n-1}\cdot\Phi\bigl(n+1, n+2;-\mu_{n}^2 r^2y\bigr)dy
\quad\quad\quad\quad\quad\quad
\end{eqnarray}

Let us denote $f_{1\tau}(y) = y^2$ and get:

\begin{eqnarray}\label{eqn234}
R_{1,n-1,r}(r, t)\;=\;\frac{F_{n} r^{n-1}}{8\mu_{n}^4\nu}\int_{\frac{1}{[4\mu_{n}^2\nu t+1]}}^{1}y^{n}\cdot \Phi\bigl(n+2, n+2;-\mu_{n}^2 r^2y\bigr)dy
\quad\quad\quad\quad\quad\quad
\end{eqnarray}

\begin{eqnarray}\label{eqn235}
R_{1,n+1,r}(r, t)\;=\frac{F_{n} r^{n+1}}{8\mu_{n}^2\nu(n+1)}\int_{\frac{1}{[4\mu_{n}^2\nu t+1]}}^{1} y^{n+1}\cdot\Phi\bigl(n+1, n+2;-\mu_{n}^2 r^2y\bigr)dy
\quad\quad\quad\quad\quad\quad
\end{eqnarray}

We use formula $(\ref{A12})$ for integrand in the integral $(\ref{eqn234})$ and formula $(\ref{A13})$ for integrand in the integral $(\ref{eqn235})$ $\cite{BE153}$, integrate and then we have:

\begin{eqnarray}\label{eqn236}
R_{1,n-1,r}(r, t)\;=\;\frac{F_{n} r^{n-1}}{8\mu_{n}^4\nu(n+1)}\biggl[\Phi\bigl(n+1, n+2;-\mu_{n}^2 r^2\bigr) - \frac{\Phi\bigl(n+1, n+2;-\frac{\mu_{n}^2 r^2}{(4\mu_{n}^2\nu t+1)}\bigr)}{(4\mu_{n}^2\nu t+1)^{n+1}}\biggr]
\nonumber\\
\nonumber\\
\quad\quad\quad\quad\quad\quad
\end{eqnarray}

\begin{eqnarray}\label{eqn237}
R_{1,n+1,r}(r, t)\;=\;\frac{F_{n} r^{n+1}}{8\mu_{n}^2\nu(n+1)(n+2)}\biggl[\Phi\bigl(n+1, n+3;-\mu_{n}^2 r^2\bigr) - \frac{\Phi\bigl(n+1, n+3;-\frac{\mu_{n}^2 r^2}{(4\mu_{n}^2\nu t+1)}\bigr)}{(4\mu_{n}^2\nu t+1)^{n+2}}\biggr]
\nonumber\\
\nonumber\\
\quad\quad\quad\quad\quad\quad
\end{eqnarray}

Hence, and from formulas $(\ref{eqn221}), (\ref{eqn222})$ it follows for j = 1:

\begin{equation}\label{eqn238}
u_{1r1}(r,\varphi, t)\;=\;\frac{n}{2}\bigl[R_{1,n-1,r}(r, t)\;\sz{e}  ^{i (n-1) \varphi}\;+\;R_{1,n+1,r}(r, t)\;\sz{e}  ^{i (n+1) \varphi}\bigr]
\end{equation}

\begin{equation}\label{eqn239}
u_{1r2}(r,\varphi, t)\;=\;\frac{i\;n}{2}\bigl[R_{1,n-1,r}(r, t)\;\sz{e}  ^{i (n-1) \varphi}\;-\;R_{1,n+1,r}(r, t)\;\sz{e}  ^{i (n+1) \varphi}\bigr]
\end{equation}

We continue operations and get:

\begin{equation}\label{eqn240}
u_{1 r}(r,\varphi, t)\;=\;\frac{n}{2}\bigl[R_{1,n-1,r}(r, t)\;+\;R_{1,n+1,r}(r, t)\bigr]\;\sz{e}  ^{i n \varphi}\;
\end{equation}

\begin{equation}\label{eqn241}
u_{1\varphi}(r,\varphi, t)\;=\;\frac{i\;n}{2}\bigl[R_{1,n-1,r}(r, t)\;-\;R_{1,n+1,r}(r, t)\bigr]\;\sz{e}^{i n \varphi}\;
\end{equation}

$u_{1 r}(r,\varphi, t) , \;u_{1\varphi}(r,\varphi, t)$ are the radial and tangential components of the velocity $\vec u_{1}$.

We use the properties of the  confluent hypergeometric function $\Phi(a, c; x)$ and have from formulas $(\ref{eqn236})-(\ref{eqn241})$:

\begin{equation}
\begin{array}{ll} 
lim \;\;u_{1r1}(r,\varphi,t)= 0;\;\;\;   lim \;\;u_{1r2}(r,\varphi,t)= 0;\\
t \rightarrow 0; \;\;\;\;\;\;\;\;\;\;\;\;\;\;\;\;\;\;\;\;\;\;\;\;\;\; t \rightarrow 0;\\
lim \;\;u_{1r}(r,\varphi,t)= 0;\;\;\;   lim \;\;u_{1\varphi}(r,\varphi,t)= 0;\\
t \rightarrow 0; \;\;\;\;\;\;\;\;\;\;\;\;\;\;\;\;\;\;\;\;\;\;\;\;\;\; t \rightarrow 0;\\
\end{array}
\end{equation}

Hence we have velocity $\vec u_{1}$ satisfies the initial conditions $(\ref{eqn200})$.

Then we use the asymptotic properties of the confluent hypergeometric function $\Phi(a, c; x)\; \cite{BE153}$ and from formulas $(\ref{eqn236})-(\ref{eqn241})$ we have velocity $\vec u_{1}$ satisfies conditions  $(\ref{eqn16})$ ( for $r \;\rightarrow\; \infty )$.

Let us continue investigation for the second step (j = 2) of the iterative process. 

Find $\vec{f}_{2}^{*}(r,\varphi, t) = \{f_{21}^{*}, f_{22}^{*}\}$ - the first correction of the particular radial applied force ${f}_{1}(x,t)$ $(\ref{eqn228})$.

We have for $\vec{f}_{2}^{*}$ from formula $(\ref{eqn48})$:

\begin{equation}\label{eqn242}
f_{21}^{*} =  u_{1r1}\;\frac{\partial u_{1r1}}{\partial x_{1}}\;+\;u_{1r2}\;\frac{\partial u_{1r1}}{\partial x_{2}}
\end{equation}

\begin{equation}\label{eqn243}
f_{22}^{*} =  u_{1r1}\;\frac{\partial u_{1r2}}{\partial x_{1}}\;+\;u_{1r2}\;\frac{\partial u_{1r2}}{\partial x_{2}}
\end{equation}

where $u_{1r1},\; u_{1r2}$ are the components of $\vec{u_{1}}$ and were taken from formulas $(\ref{eqn238}) , (\ref{eqn239})$.

We have here:

\begin{eqnarray}\label{eqn244}
\frac{\partial u_{1r1}(r,\varphi, t)}{\partial x_{1}} \; = \; \frac{\partial u_{1r1}(r,\varphi, t)}{\partial r}\; \frac{\partial r}{\partial x_{1}} \; + \; \frac{\partial u_{1r1}(r,\varphi, t)}{\partial \varphi}\; \frac{\partial \varphi}{\partial x_{1}}
\nonumber\\
\nonumber\\
\frac{\partial u_{1r1}(r,\varphi, t)}{\partial x_{2}}\;  = \; \frac{\partial u_{1r1}(r,\varphi, t)}{\partial r}\; \frac{\partial r}{\partial x_{2}} \; + \; \frac{\partial u_{1r1}(r,\varphi, t)}{\partial \varphi}\; \frac{\partial \varphi}{\partial x_{2}}
\nonumber\\
\nonumber\\
\frac{\partial u_{1r2}(r,\varphi, t)}{\partial x_{1}}\;  = \; \frac{\partial u_{1r2}(r,\varphi, t)}{\partial r}\; \frac{\partial r}{\partial x_{1}} \; + \; \frac{\partial u_{1r2}(r,\varphi, t)}{\partial \varphi}\; \frac{\partial \varphi}{\partial x_{1}}
\nonumber\\
\nonumber\\
\frac{\partial u_{1r2}(r,\varphi, t)}{\partial x_{2}}\;  = \; \frac{\partial u_{1r2}(r,\varphi, t)}{\partial r}\; \frac{\partial r}{\partial x_{2}} \; + \; \frac{\partial u_{1r2}(r,\varphi, t)}{\partial \varphi}\; \frac{\partial \varphi}{\partial x_{2}}
\nonumber\\
\nonumber\\
\frac{\partial r}{\partial x_{1}} = cos \varphi \; , \; \frac{\partial \varphi}{\partial x_{1}} = \;-\; \frac{sin \varphi}{r}\;,\; \frac{\partial r}{\partial x_{2}} = sin \varphi \; , \;\frac{\partial \varphi}{\partial x_{2}} =  \frac{cos \varphi}{r}
\nonumber\\
\nonumber\\
\end{eqnarray}

Then we use formulas $(\ref{eqn238}) , (\ref{eqn239})$ for $u_{1r1}(r,\varphi, t),\; u_{1r2}(r,\varphi, t)$ and have from $(\ref{eqn244})$:

\begin{eqnarray}\label{eqn245}
\frac{\partial u_{1r1}(r,\varphi, t)}{\partial x_{1}} \; = \; \frac{n}{2}\biggl\{\bigl[R_{1,n-1,r}^{'}(r, t)\;\sz{e}  ^{i (n-1) \varphi}\;+\;R_{1,n+1,r}^{'}(r, t)\;\sz{e}  ^{i (n+1) \varphi}\bigr]\;cos \varphi\;+
\nonumber\\
+\;i\;\bigl[(n-1)R_{1,n-1,r}(r, t)\;\sz{e}  ^{i (n-1) \varphi}\;+\;(n+1)R_{1,n+1,r}(r, t)\;\sz{e}  ^{i (n+1) \varphi}\bigr]\bigl(-\frac{sin \varphi}{r} \bigr ) \biggr \}
\nonumber\\
\nonumber\\
\frac{\partial u_{1r1}(r,\varphi, t)}{\partial x_{2}}\;  = \; \frac{n}{2}\biggl\{\bigl[R_{1,n-1,r}^{'}(r, t)\;\sz{e}  ^{i (n-1) \varphi}\;+\;R_{1,n+1,r}^{'}(r, t)\;\sz{e}  ^{i (n+1) \varphi}\bigr]\;sin \varphi\;+
\nonumber\\
+\;i\;\bigl[(n-1)R_{1,n-1,r}(r, t)\;\sz{e}  ^{i (n-1) \varphi}\;+\;(n+1)R_{1,n+1,r}(r, t)\;\sz{e}  ^{i (n+1) \varphi}\bigr]\frac{cos \varphi}{r}  \biggr \}
\nonumber\\
\nonumber\\
\frac{\partial u_{1r2}(r,\varphi, t)}{\partial x_{1}}\;  = \; \frac{i n}{2}\biggl\{\bigl[R_{1,n-1,r}^{'}(r, t)\;\sz{e}  ^{i (n-1) \varphi}\;-\;R_{1,n+1,r}^{'}(r, t)\;\sz{e}  ^{i (n+1) \varphi}\bigr]\;cos \varphi\;+
\nonumber\\
+\;i\;\bigl[(n-1)R_{1,n-1,r}(r, t)\;\sz{e}  ^{i (n-1) \varphi}\;-\;(n+1)R_{1,n+1,r}(r, t)\;\sz{e}  ^{i (n+1) \varphi}\bigr]\bigl(-\frac{sin \varphi}{r} \bigr ) \biggr \}
\nonumber\\
\nonumber\\
\frac{\partial u_{1r2}(r,\varphi, t)}{\partial x_{2}}\;  = \; \frac{i n}{2}\biggl\{\bigl[R_{1,n-1,r}^{'}(r, t)\;\sz{e}  ^{i (n-1) \varphi}\;-\;R_{1,n+1,r}^{'}(r, t)\;\sz{e}  ^{i (n+1) \varphi}\bigr]\;sin \varphi\;+
\nonumber\\
+\;i\;\bigl[(n-1)R_{1,n-1,r}(r, t)\;\sz{e}  ^{i (n-1) \varphi}\;-\;(n+1)R_{1,n+1,r}(r, t)\;\sz{e}  ^{i (n+1) \varphi}\bigr]\frac{cos \varphi}{r}  \biggr \}
\nonumber\\
\nonumber\\
\end{eqnarray}

where 

\begin{eqnarray}\label{eqn246}
R_{1,n-1,r}^{'}(r, t)\;=\;\frac{\partial R_{1,n-1,r}(r, t)}{\partial r}\;=\;\frac{F_{n}(n-1) r^{n-2}}{8\mu_{n}^4\nu(n+1)}\biggl[\Phi\bigl(n+1, n+2;-\mu_{n}^2 r^2\bigr) - \frac{\Phi\bigl(n+1, n+2;-\frac{\mu_{n}^2 r^2}{(4\mu_{n}^2\nu t+1)}\bigr)}{(4\mu_{n}^2\nu t+1)^{n+1}}\biggr] -
\nonumber\\
\nonumber\\
\quad\quad\quad\quad\quad
- \frac{F_{n} r^{n}}{4\mu_{n}^2\nu(n+2)}\biggl[\Phi\bigl(n+2, n+3;-\mu_{n}^2 r^2\bigr) - \frac{\Phi\bigl(n+2, n+3;-\frac{\mu_{n}^2 r^2}{(4\mu_{n}^2\nu t+1)}\bigr)}{(4\mu_{n}^2\nu t+1)^{n+2}}\biggr]
\quad\quad\quad\quad\quad\quad\quad\quad\quad\quad
\nonumber\\
\nonumber\\
\end{eqnarray}
\begin{eqnarray}\label{eqn247}
R_{1,n+1,r}^{'}(r, t)\;=\;\frac{\partial R_{1,n+1,r}(r, t)}{\partial r}\;=\;\frac{F_{n} r^{n}}{8\mu_{n}^2\nu(n+2)}\biggl[\Phi\bigl(n+1, n+3;-\mu_{n}^2 r^2\bigr) - \frac{\Phi\bigl(n+1, n+3;-\frac{\mu_{n}^2 r^2}{(4\mu_{n}^2\nu t+1)}\bigr)}{(4\mu_{n}^2\nu t+1)^{n+2}}\biggr] -
\nonumber\\
\nonumber\\
- \frac{F_{n} r^{n+2}}{4\nu(n+2)(n+3)}\biggl[\Phi\bigl(n+2, n+4;-\mu_{n}^2 r^2\bigr) - \frac{\Phi\bigl(n+2, n+4;-\frac{\mu_{n}^2 r^2}{(4\mu_{n}^2\nu t+1)}\bigr)}{(4\mu_{n}^2\nu t+1)^{n+3}}\biggr]
\quad\quad\quad\quad\quad\quad\quad\quad
\nonumber\\
\nonumber\\
\end{eqnarray}

Let us put $u_{1r1},\; u_{1r2}, \;\frac{\partial u_{1r1}}{\partial x_{1}}, \;\frac{\partial u_{1r1}}{\partial x_{2}}, \;\frac{\partial u_{1r2}}{\partial x_{1}},\; \frac{\partial u_{1r2}}{\partial x_{2}}$ from formulas $(\ref{eqn238})  , (\ref{eqn239}) , (\ref{eqn245})$ in formulas $(\ref{eqn242}), (\ref{eqn243})$ for $f_{21}^{*} \;,\; f_{22}^{*}$. 
\\
\\
After compliting appropriate operations we have:

\begin{equation}\label{eqn248}
f_{21}^{*} (r,\varphi, t) =  \frac{n^{2}}{2^{2}}\bigl [ T_{2,2n-1,r}(r, t)\;\sz{e}  ^{i (2n-1) \varphi} + T_{2,2n+1,r}(r, t)\;\sz{e}  ^{i (2n+1) \varphi} \bigr ]
\end{equation}

\begin{equation}\label{eqn249}
f_{22}^{*} (r,\varphi, t) =  \frac{i n^{2}}{2^{2}}\bigl [ T_{2,2n-1,r}(r, t)\;\sz{e}  ^{i (2n-1) \varphi} - T_{2,2n+1,r}(r, t)\;\sz{e}  ^{i (2n+1) \varphi} \bigr ]
\end{equation}

where

\begin{eqnarray}\label{eqn250}
T_{2,2n-1,r}(r, t) = \bigl [ R_{1,n-1,r}(r, t) + R_{1,n+1,r}(r, t) \bigr ] R_{1,n-1,r}^{'}(r, t) - 
\nonumber\\
\nonumber\\
- \frac{(n-1)}{r} R_{1,n-1,r}(r, t)\bigl [ R_{1,n-1,r}(r, t) - R_{1,n+1,r}(r, t) \bigr ]
\nonumber\\
\nonumber\\
T_{2,2n+1,r}(r, t) = \bigl [ R_{1,n-1,r}(r, t) + R_{1,n+1,r}(r, t) \bigr ] R_{1,n+1,r}^{'}(r, t) - 
\nonumber\\
\nonumber\\
- \frac{(n+1)}{r} R_{1,n+1,r}(r, t)\bigl [ R_{1,n-1,r}(r, t) - R_{1,n+1,r}(r, t) \bigr ]
\nonumber\\
\nonumber\\
\quad\quad
\end{eqnarray}

For radial $f_{2r}^{*}$ and tangential $f_{2\varphi}^{*}$ components of the first correction $\vec{f_{2}^{*}}(r, \varphi, t)$ of the particular radial applied force we have:

\begin{equation}\label{eqn251}
f_{2r}^{*} (r,\varphi, t) =  \frac{n^{2}}{2^{2}}\bigl [ T_{2,2n-1,r}(r, t) + T_{2,2n+1,r}(r, t) \bigr ]\; \sz{e}  ^{i\; 2n \varphi}\;=\;\frac{n^{2}}{2^{2}}T_{2,2n,r}(r, t)\; \sz{e}  ^{i\; 2n \varphi}\;
\end{equation}

\begin{equation}\label{eqn252}
f_{2\varphi}^{*} (r,\varphi, t) =  \frac{i n^{2}}{2^{2}}\bigl [ T_{2,2n-1,r}(r, t) - T_{2,2n+1,r}(r, t) \bigr ]\; \sz{e}  ^{i\; 2n \varphi}\;=\;\frac{i n^{2}}{2^{2}}T_{2,2n,\varphi}(r, t)\; \sz{e}  ^{i\; 2n \varphi}\;
\end{equation}

where ( see $(\ref{eqn250})$)

\begin{eqnarray}\label{eqn253}
T_{2,2n,r}(r, t) = \bigl [ R_{1,n-1,r}(r, t) + R_{1,n+1,r}(r, t) \bigr ] \bigl [R_{1,n-1,r}^{'}(r, t) + R_{1,n+1,r}^{'}(r, t)\bigr ] - 
\nonumber\\
\nonumber\\
- \frac{1}{r}\bigl [(n-1)R_{1,n-1,r}(r, t) + (n+1)R_{1,n+1,r}(r, t) \bigr ]\bigl [ R_{1,n-1,r}(r, t) - R_{1,n+1,r}(r, t) \bigr ]
\nonumber\\
\nonumber\\
T_{2,2n,\varphi}(r, t) = \bigl [ R_{1,n-1,r}(r, t) + R_{1,n+1,r}(r, t) \bigr ] \bigl [R_{1,n-1,r}^{'}(r, t) - R_{1,n+1,r}^{'}(r, t)\bigr ] - 
\nonumber\\
\nonumber\\
- \frac{1}{r}\bigl [(n-1)R_{1,n-1,r}(r, t) - (n+1)R_{1,n+1,r}(r, t) \bigr ]\bigl [ R_{1,n-1,r}(r, t) - R_{1,n+1,r}(r, t) \bigr ]
\nonumber\\
\nonumber\\
\quad\quad
\end{eqnarray}

We use formulas $(\ref{eqn236}) , (\ref{eqn237})$ for $R_{1,n-1,r}(r, t)\;,\;R_{1,n+1,r}(r, t)$ and $(\ref{eqn246}) , (\ref{eqn247})$ for  $R_{1,n-1,r}^{'}(r, t)\;$ , $\;R_{1,n+1,r}^{'}(r, t)$ then do the appropriate operations for $T_{2,2n-1,r}(r, t)\;,\;T_{2,2n+1,r}(r, t)$, using formula $(\ref{A14})$, and get:

\begin{eqnarray}\label{eqn254}
T_{2,2n,r}(r, t) = \frac{- F_{n}^2\cdot r^{2n-1}}{16\mu_{n}^6\nu^2(n+1)^2(n+2)}\biggl[\Phi\bigl(n+1, n+2;-\mu_{n}^2 r^2\bigr) - \frac{\Phi\bigl(n+1, n+2;-\frac{\mu_{n}^2 r^2}{(4\mu_{n}^2\nu t+1)}\bigr)}{(4\mu_{n}^2\nu t+1)^{n+1}}\biggr]\cdot
\nonumber\\
\nonumber\\
\cdot\biggl[\Phi\bigl(n+1, n+3;-\mu_{n}^2 r^2\bigr) - \frac{\Phi\bigl(n+1, n+3;-\frac{\mu_{n}^2 r^2}{(4\mu_{n}^2\nu t+1)}\bigr)}{(4\mu_{n}^2\nu t+1)^{n+2}}\biggr]
\quad\quad\quad\quad\quad\quad\quad\quad\quad
\nonumber\\
\nonumber\\
\quad\quad\quad\quad\quad\quad\quad\quad
\end{eqnarray}

\begin{eqnarray}\label{eqn255}
T_{2,2n,\varphi}(r, t) = \frac{- F_{n}^2\cdot r^{2n-1}}{16\mu_{n}^6\nu^2(n+1)(n+2)}\biggl[\Phi\bigl(n, n+2;-\mu_{n}^2 r^2\bigr) - \frac{\Phi\bigl(n, n+2;-\frac{\mu_{n}^2 r^2}{(4\mu_{n}^2\nu t+1)}\bigr)}{(4\mu_{n}^2\nu t+1)^{n+1}}\biggr]\cdot
\nonumber\\
\nonumber\\
\cdot\biggl[\Phi\bigl(n+2, n+3;-\mu_{n}^2 r^2\bigr) - \frac{\Phi\bigl(n+2, n+3;-\frac{\mu_{n}^2 r^2}{(4\mu_{n}^2\nu t+1)}\bigr)}{(4\mu_{n}^2\nu t+1)^{n+2}}\biggr] +
\quad\quad\quad\quad\quad\quad\quad
\nonumber\\
\nonumber\\
+ \frac{F_{n}^2 \cdot n \cdot r^{2n-1}}{16\mu_{n}^6\nu^2(n+1)^2(n+2)}\biggl[\Phi\bigl(n+1, n+2;-\mu_{n}^2 r^2\bigr) - \frac{\Phi\bigl(n+1, n+2;-\frac{\mu_{n}^2 r^2}{(4\mu_{n}^2\nu t+1)}\bigr)}{(4\mu_{n}^2\nu t+1)^{n+1}}\biggr]\cdot
\nonumber\\
\nonumber\\
\cdot\biggl[\Phi\bigl(n+1, n+3;-\mu_{n}^2 r^2\bigr) - \frac{\Phi\bigl(n+1, n+3;-\frac{\mu_{n}^2 r^2}{(4\mu_{n}^2\nu t+1)}\bigr)}{(4\mu_{n}^2\nu t+1)^{n+2}}\biggr]
\quad\quad\quad\quad\quad\quad\quad\quad\quad
\nonumber\\
\nonumber\\
\quad\quad\quad\quad\quad\quad\quad\quad
\end{eqnarray}

By comparing particular radial applied force $\vec{f_{1}}\;$ from $\;(\ref{eqn228})$ with the first correction $\vec{f_{2}^{*}}\;$ from $\;((\ref{eqn251})- (\ref{eqn255}))$ of this particular radial applied force  we have:

\begin{equation}\label{eqn256}
\mid\vec{f_{2}^{*}}\mid\; <<\; \mid\vec{f_{1}}\mid
\end{equation}

with condition 

\begin{equation}\label{eqn257}
F_{n}\;\leq\;\frac{1}{n}
\end{equation}

After the first step of the iterative process (j = 1) we had the velocity $\vec{u_{1}}\;$ - formulas $\; (\ref{eqn240})\;,\;(\ref{eqn241})$.
\\  
Now we will calculate $\vec{u_{2}^{*}}\;$ - the first correction of the velocity $\vec{u_{1}} $. Solution of this problem has two stages. On the first stage we find the part of the first correction $\vec{u_{2r}^{*}}$, corresponding to the first correction $f_{2r}^{*}\;$ from formula $\;(\ref{eqn251})$ of the applied force: 

\begin{equation}\label{eqn258}
f_{2r}^{*}(r,\varphi,t) = \;\frac{n^{2}}{2^{2}}T_{2,2n,r}(r, t)\; \sz{e}  ^{i\; 2n \varphi}\;\;\;,\;\;\;f_{2\varphi}^{*}(r,\varphi,t) \equiv 0
\end{equation}

On the second stage we calculate the other part of the first correction $\vec{u_{2\varphi}^{*}}\;$, corresponding to  the first correction $f_{2\varphi}^{*}\;$ from formula $\;(\ref{eqn252})$ of the applied force:

\begin{equation}\label{eqn259}
f_{2r}^{*}(r,\varphi,t) \equiv 0\;\;\;,\;\;\; f_{2\varphi}^{*}(r,\varphi,t) = \;\frac{i n^{2}}{2^{2}}T_{2,2n,\varphi}(r, t)\; \sz{e}  ^{i\; 2n\varphi}
\end{equation}

In other words 

\begin{eqnarray}\label{eqn260}
\vec{u_{2}^{*}} = \vec{u_{2r}^{*}} + \vec{u_{2\varphi}^{*}}, \;\;\;\vec{u_{2r}^{*}} = \{u_{2r1}^{*}, u_{2r2}^{*}\}, \;\;\;\vec{u_{2\varphi}^{*}} = \{u_{2\varphi1}^{*}, u_{2\varphi2}^{*}\}.
\end{eqnarray}

First stage: After completing appropriate operations we have from formula $(\ref{eqn254})$:

\begin{eqnarray}\label{eqn261}
T_{2,2n,r}(r, t) = \frac{- F_{n}^2\cdot r^{2n-1}}{16\mu_{n}^6\nu^2(n+1)^2(n+2)}\biggl[\Phi\bigl(n+1, n+2;-\mu_{n}^2 r^2\bigr)\cdot \Phi\bigl(n+1, n+3;-\mu_{n}^2 r^2\bigr)- 
\nonumber\\
\nonumber\\
- \frac{1}{(4\mu_{n}^2\nu t+1)^{n+2}}\cdot\Phi\bigl(n+1, n+2;-\mu_{n}^2 r^2\bigr)\cdot\Phi\bigl(n+1, n+3;-\frac{\mu_{n}^2 r^2}{(4\mu_{n}^2\nu t+1)}\bigr) -
\quad\quad\quad\quad\quad
\nonumber\\
\nonumber\\
- \frac{1}{(4\mu_{n}^2\nu t+1)^{n+1}}\cdot\Phi\bigl(n+1, n+2;-\frac{\mu_{n}^2 r^2}{(4\mu_{n}^2\nu t+1)}\bigr)\cdot\Phi\bigl(n+1, n+3;-\mu_{n}^2 r^2\bigr) +
\quad\quad\quad\quad\quad
\nonumber\\
\nonumber\\
+ \frac{1}{(4\mu_{n}^2\nu t+1)^{2n+3}}\cdot \Phi\bigl(n+1, n+2;-\frac{\mu_{n}^2 r^2}{(4\mu_{n}^2\nu t+1)}\bigr)\cdot\Phi\bigl(n+1, n+3;-\frac{\mu_{n}^2 r^2}{(4\mu_{n}^2\nu t+1)}\bigr)\biggr]
\quad\quad
\nonumber\\
\nonumber\\
\quad\quad\quad\quad\quad
\end{eqnarray}

We take formulas $(\ref{eqn221})\; ,\; (\ref{eqn222})$ for components $u_{jr1}(r, \varphi, t),\; u_{jr2}(r, \varphi, t)$ and formulas $(\ref{eqn219}) , (\ref{eqn220})$ for $R_{j,n_{j}-1,r}(r, t),\; R_{j,n_{j}+1,r}(r, t)$ for j = 2 and then formulas $(\ref{eqn258}) , (\ref{eqn261})$ for $f_{2r}^{*} (r,\varphi, t)\;,\;T_{2,2n,r}(r, t)$ . We do appropriate operations and have:

\begin{equation}\label{eqn264}
u_{2r1}^{*}(r,\varphi, t)\;=\;n\bigl[R_{2,2n-1,r}(r, t)\;\sz{e}  ^{i (2n-1) \varphi}\;+\;R_{2,2n+1,r}(r, t)\;\sz{e}  ^{i (2n+1) \varphi}\bigr]
\end{equation}

\begin{equation}\label{eqn265}
u_{2r2}^{*}(r,\varphi, t)\;=\;i\;n\bigl[R_{2,2n-1,r}(r, t)\;\sz{e}  ^{i (2n-1) \varphi}\;-\;R_{2,2n+1,r}(r, t)\;\sz{e}  ^{i (2n+1) \varphi}\bigr]
\end{equation}

After changing the order of integration we receive:  

\begin{equation}\label{eqn266}
R_{2,2n-1,r}(r, t)\;=\;\frac{n^{2}}{2^{2}}\int_{0}^{t}\int_{0}^{\infty} T_{2,2n,r}(\tilde r, \tau)\int_{0}^{\infty} \sz{e}  ^{-\nu \rho^2 (t-\tau)} J_{2n-1}(r \rho)\;J_{2n}(\tilde r \rho)\;d\rho d \tilde r d\tau  
\quad\quad
\end{equation}

\begin{equation}\label{eqn267}
R_{2,2n+1,r}(r, t)\;=\;\frac{n^{2}}{2^{2}}\int_{0}^{t}\int_{0}^{\infty} T_{2,2n,r}(\tilde r, \tau) \int_{0}^{\infty}\sz{e}  ^{-\nu \rho^2 (t-\tau)} J_{2n+1}(r \rho)\;J_{2n}(\tilde r \rho)\;d\rho d \tilde r d\tau  
\quad\quad
\end{equation}


Second stage: After operations with $T_{2,2n,\varphi}(r, t)$ from formula $(\ref{eqn255})$ and using formula  $(\ref{A14})$, we obtain:

\begin{eqnarray}\label{eqn278}
\bigl(T_{2,2n,\varphi}( r, t)\cdot r\bigr)^{'}_{r} = \frac{\partial \bigl(T_{2,2n,\varphi}( r, t)\cdot r\bigr)}{\partial r} = T_{2,2n,\varphi, \varphi}(r, t) + T_{2,2n,\varphi, r}(r, t)
\nonumber\\
\nonumber\\
\end{eqnarray}

We denote here

\begin{eqnarray}\label{eqn279}
T_{2,2n,\varphi, \varphi}(r, t) = \frac{- F_{n}^2\cdot n \cdot r^{2n+1}}{8\mu_{n}^4\nu^2(n+1)(n+2)^2}\biggl[\Phi\bigl(n+1, n+3;-\mu_{n}^2 r^2\bigr) - \frac{\Phi\bigl(n+1, n+3;-\frac{\mu_{n}^2 r^2}{(4\mu_{n}^2\nu t+1)}\bigr)}{(4\mu_{n}^2\nu t+1)^{n+2}}\biggr]\cdot
\nonumber\\
\nonumber\\
\cdot\biggl[\Phi\bigl(n+2, n+3;-\mu_{n}^2 r^2\bigr) - \frac{\Phi\bigl(n+2, n+3;-\frac{\mu_{n}^2 r^2}{(4\mu_{n}^2\nu t+1)}\bigr)}{(4\mu_{n}^2\nu t+1)^{n+2}}\biggr] +
\quad\quad\quad\quad\quad\quad
\nonumber\\
\nonumber\\
+\frac{F_{n}^2\cdot r^{2n+1}}{8\mu_{n}^4\nu^2(n+1)(n+3)}\biggl[\Phi\bigl(n, n+2;-\mu_{n}^2 r^2\bigr) - \frac{\Phi\bigl(n, n+2;-\frac{\mu_{n}^2 r^2}{(4\mu_{n}^2\nu t+1)}\bigr)}{(4\mu_{n}^2\nu t+1)^{n+1}}\biggr]\cdot
\quad\quad\quad\quad
\nonumber\\
\nonumber\\
\cdot\biggl[\Phi\bigl(n+3, n+4;-\mu_{n}^2 r^2\bigr) - \frac{\Phi\bigl(n+3, n+4;-\frac{\mu_{n}^2 r^2}{(4\mu_{n}^2\nu t+1)}\bigr)}{(4\mu_{n}^2\nu t+1)^{n+3}}\biggr]
\quad\quad\quad\quad\quad\quad
\nonumber\\
\nonumber\\
\end{eqnarray}

and

\begin{eqnarray}\label{eqn280}
T_{2,2n,\varphi, r}(r, t) =  \frac{- F_{n}^2 \cdot n \cdot r^{2n-1}}{8\mu_{n}^6\nu^2(n+1)^2(n+2)}\biggl[\Phi\bigl(n+1, n+2;-\mu_{n}^2 r^2\bigr) - \frac{\Phi\bigl(n+1, n+2;-\frac{\mu_{n}^2 r^2}{(4\mu_{n}^2\nu t+1)}\bigr)}{(4\mu_{n}^2\nu t+1)^{n+1}}\biggr]\cdot
\quad
\nonumber\\
\nonumber\\
\cdot\biggl[\Phi\bigl(n+1, n+3;-\mu_{n}^2 r^2\bigr) - \frac{\Phi\bigl(n+1, n+3;-\frac{\mu_{n}^2 r^2}{(4\mu_{n}^2\nu t+1)}\bigr)}{(4\mu_{n}^2\nu t+1)^{n+2}}\biggr] = 2 n \cdot T_{2,2n, r}(r, t)
\quad\quad\quad\quad
\nonumber\\
\nonumber\\
\end{eqnarray}

$T_{2,2n, r}(r, t)$ is taken  from formulas $(\ref{eqn254})$ and $(\ref{eqn261})$. 

Then we do several operations and have from formula $(\ref{eqn279})$:

\begin{eqnarray}\label{eqn281}
T_{2,2n,\varphi, \varphi}(r, t) = \frac{- F_{n}^2\cdot n \cdot r^{2n+1}}{8\mu_{n}^4\nu^2(n+1)(n+2)^2}\biggl[\Phi\bigl(n+1, n+3;-\mu_{n}^2 r^2\bigr) \cdot\Phi\bigl(n+2, n+3;-\mu_{n}^2 r^2\bigr) - 
\nonumber\\
\nonumber\\
- \frac{1}{(4\mu_{n}^2\nu t+1)^{n+2}}\cdot\Phi\bigl(n+1, n+3;-\mu_{n}^2 r^2\bigr)\cdot\Phi\bigl(n+2, n+3;-\frac{\mu_{n}^2 r^2}{(4\mu_{n}^2\nu t+1)}\bigr) -
\quad\quad
\nonumber\\
\nonumber\\
- \frac{1}{(4\mu_{n}^2\nu t+1)^{n+2}}\cdot\Phi\bigl(n+1, n+3;-\frac{\mu_{n}^2 r^2}{(4\mu_{n}^2\nu t+1)}\bigr)\cdot \Phi\bigl(n+2, n+3;-\mu_{n}^2 r^2\bigr) + 
\quad\quad
\nonumber\\
\nonumber\\
+ \frac{1}{(4\mu_{n}^2\nu t+1)^{2n+4}}\cdot\Phi\bigl(n+1, n+3;-\frac{\mu_{n}^2 r^2}{(4\mu_{n}^2\nu t+1)}\bigr)\cdot\Phi\bigl(n+2, n+3;-\frac{\mu_{n}^2 r^2}{(4\mu_{n}^2\nu t+1)}\bigr)\biggr ] + 
\nonumber\\
\nonumber\\
+ \frac{F_{n}^2\cdot r^{2n+1}}{8\mu_{n}^4\nu^2(n+1)(n+3)}\biggl[\Phi\bigl(n, n+2;-\mu_{n}^2 r^2\bigr) \cdot\Phi\bigl(n+3, n+4;-\mu_{n}^2 r^2\bigr) - 
\quad\quad\quad\quad
\nonumber\\
\nonumber\\
- \frac{1}{(4\mu_{n}^2\nu t+1)^{n+3}}\cdot\Phi\bigl(n, n+2;-\mu_{n}^2 r^2\bigr)\cdot\Phi\bigl(n+3, n+4;-\frac{\mu_{n}^2 r^2}{(4\mu_{n}^2\nu t+1)}\bigr) -
\quad\quad
\quad\quad
\nonumber\\
\nonumber\\
- \frac{1}{(4\mu_{n}^2\nu t+1)^{n+1}}\cdot\Phi\bigl(n, n+2;-\frac{\mu_{n}^2 r^2}{(4\mu_{n}^2\nu t+1)}\bigr)\cdot \Phi\bigl(n+3, n+4;-\mu_{n}^2 r^2\bigr) + 
\quad\quad\quad\quad
\nonumber\\
\nonumber\\
+ \frac{1}{(4\mu_{n}^2\nu t+1)^{2n+4}}\cdot\Phi\bigl(n, n+2;-\frac{\mu_{n}^2 r^2}{(4\mu_{n}^2\nu t+1)}\bigr)\cdot\Phi\bigl(n+3, n+4;-\frac{\mu_{n}^2 r^2}{(4\mu_{n}^2\nu t+1)}\bigr)\biggr ]
\quad\quad
\nonumber\\
\nonumber\\
\end{eqnarray}

We transform formulas $(\ref{eqn226}), (\ref{eqn227})$ for components $u_{j\varphi1}(r, \varphi),\; u_{j\varphi2}(r, \varphi)$ and formulas $(\ref{eqn224}), (\ref{eqn225})\;\;$ for $R_{j,n_{j}-1,\varphi}(r),\; R_{j,n_{j}+1,\varphi}(r)$ for j = 2 and then use formula $(\ref{eqn259})$ for $f_{2\varphi}^{*} (r,\varphi, t)$. We do several operations and have:

\begin{equation}\label{eqn284}
u_{2\varphi1}^{*}(r,\varphi, t)\;=\;- \;\frac{i}{2}\bigl[R_{2,2n-1,\varphi}(r, t)\;\sz{e}  ^{i (2n-1) \varphi}\;+\;R_{2,2n+1,\varphi}(r, t)\;\sz{e}  ^{i (2n+1) \varphi}\bigr]
\end{equation}

\begin{equation}\label{eqn285}
u_{2\varphi2}^{*}(r,\varphi, t)\;=\;\frac{1}{2}\bigl[R_{2,2n-1,\varphi}(r, t)\;\sz{e}  ^{i (2n-1) \varphi}\;-\;R_{2,2n+1,\varphi}(r, t)\;\sz{e}  ^{i (2n+1) \varphi}\bigr]
\end{equation}

After changing the order of integration we obtain:  

\begin{equation}\label{eqn286}
R_{2,2n-1,\varphi}(r, t) = -\;\frac{in^{2}}{2^{2}}\int_{0}^{t}\int_{0}^{\infty} \bigl(T_{2,2n,\varphi}(\tilde r, \tau)\cdot\tilde r\bigr)^{'}_{\tilde r} \int_{0}^{\infty}\sz{e}  ^{-\nu \rho^2 (t-\tau)}J_{2n-1}(r \rho)\;J_{2n}(\tilde r \rho)\;d\rho d \tilde r d\tau  
\end{equation}

\begin{equation}\label{eqn287}
R_{2,2n+1,\varphi}(r, t) = -\;\frac{in^{2}}{2^{2}}\int_{0}^{t}\int_{0}^{\infty} \bigl(T_{2,2n,\varphi}(\tilde r, \tau)\cdot\tilde r\bigr)^{'}_{\tilde r} \int_{0}^{\infty}\sz{e}  ^{-\nu \rho^2 (t-\tau)}J_{2n+1}(r \rho)\;J_{2n}(\tilde r \rho)\;d\rho d \tilde r d\tau  
\end{equation}

Then we take $\bigl(T_{2,2n,\varphi}(\tilde r, \tau)\cdot\tilde r\bigr)^{'}_{\tilde r}$ from formula $(\ref{eqn278})$  and with use of formula $(\ref{eqn280})$ put it in formulas $(\ref{eqn286}), (\ref{eqn287})$ and have:

\begin{eqnarray}\label{eqn288}
R_{2,2n-1,\varphi}(r, t) = -\;\frac{in^{2}}{2^{2}}\int_{0}^{t}\int_{0}^{\infty} T_{2,2n,\varphi, \varphi}(\tilde r, \tau) \int_{0}^{\infty}\sz{e}  ^{-\nu \rho^2 (t-\tau)}J_{2n-1}(r \rho)\cdot J_{2n}(\tilde r \rho)\;d\rho d \tilde r d\tau  -
\quad\quad\quad
\nonumber\\
\nonumber\\
-\;\frac{in^{3}}{2}\int_{0}^{t}\int_{0}^{\infty} T_{2,2n,r}(\tilde r, \tau) \int_{0}^{\infty}\sz{e}  ^{-\nu \rho^2 (t-\tau)}J_{2n-1}(r \rho)\cdot J_{2n}(\tilde r \rho)\;d\rho d \tilde r d\tau  = R_{2,2n-1,\varphi, \varphi}(r, t) - 2niR_{2,2n-1,r}(r, t)
\nonumber\\
\nonumber\\
\end{eqnarray}

\begin{eqnarray}\label{eqn289}
R_{2,2n+1,\varphi}(r, t) = -\;\frac{in^{2}}{2^{2}}\int_{0}^{t}\int_{0}^{\infty} T_{2,2n,\varphi, \varphi}(\tilde r, \tau) \int_{0}^{\infty}\sz{e}  ^{-\nu \rho^2 (t-\tau)}J_{2n+1}(r \rho)\cdot J_{2n}(\tilde r \rho)\;d\rho d \tilde r d\tau  -
\quad\quad\quad
\nonumber\\
\nonumber\\
-\;\frac{in^{3}}{2}\int_{0}^{t}\int_{0}^{\infty} T_{2,2n,r}(\tilde r, \tau) \int_{0}^{\infty}\sz{e}  ^{-\nu \rho^2 (t-\tau)}J_{2n+1}(r \rho)\cdot J_{2n}(\tilde r \rho)\;d\rho d \tilde r d\tau  = R_{2,2n+1,\varphi, \varphi}(r, t) - 2niR_{2,2n+1,r}(r, t)
\nonumber\\
\nonumber\\
\end{eqnarray}

where 

\begin{eqnarray}\label{eqn290}
R_{2,2n-1,\varphi,\varphi}(r, t) = -\;\frac{in^{2}}{2^{2}}\int_{0}^{t}\int_{0}^{\infty} T_{2,2n,\varphi, \varphi}(\tilde r, \tau) \int_{0}^{\infty}\sz{e}  ^{-\nu \rho^2 (t-\tau)}J_{2n-1}(r \rho)\cdot J_{2n}(\tilde r \rho)\;d\rho d \tilde r d\tau  
\end{eqnarray}

\begin{eqnarray}\label{eqn291}
R_{2,2n+1,\varphi,\varphi}(r, t) = -\;\frac{in^{2}}{2^{2}}\int_{0}^{t}\int_{0}^{\infty} T_{2,2n,\varphi, \varphi}(\tilde r, \tau) \int_{0}^{\infty}\sz{e}  ^{-\nu \rho^2 (t-\tau)}J_{2n+1}(r \rho)\cdot J_{2n}(\tilde r \rho)\;d\rho d \tilde r d\tau  
\end{eqnarray}

and $R_{2,2n-1,r}(r, t)\;,\; R_{2,2n+1,r}(r, t)$ we take from formulas $(\ref{eqn266}), (\ref{eqn267})$. 


Then we use formulas $(\ref{eqn260}), (\ref{eqn264}), (\ref{eqn265}), (\ref{eqn284}), (\ref{eqn285})$ and have:

\begin{eqnarray}\label{eqn302}
u_{21}^{*}(r,\varphi, t) = u_{2r1}^{*}(r,\varphi, t) + u_{2\varphi1}^{*}(r,\varphi, t)\;=
\quad\quad\quad\quad\quad\quad\quad\quad\quad\quad\quad\quad\quad\quad\quad
\nonumber\\
\nonumber\\
=\;\bigl [ n R_{2,2n-1,r}(r, t) - \frac{i}{2} R_{2,2n-1,\varphi}(r, t)\bigr ]\;\sz{e}  ^{i (2n-1) \varphi} + \bigl [ n R_{2,2n+1,r}(r, t) - \frac{i}{2} R_{2,2n+1,\varphi}(r, t)\bigr ]\;\sz{e}  ^{i (2n+1) \varphi}
\nonumber\\
\nonumber\\
\end{eqnarray}

\begin{eqnarray}\label{eqn303}
u_{22}^{*}(r,\varphi, t) = u_{2r2}^{*}(r,\varphi, t) + u_{2\varphi2}^{*}(r,\varphi, t)\;=
\quad\quad\quad\quad\quad\quad\quad\quad\quad\quad\quad\quad\quad\quad\quad
\nonumber\\
\nonumber\\
=\;i\;\bigl [ n R_{2,2n-1,r}(r, t) - \frac{i}{2} R_{2,2n-1,\varphi}(r, t)\bigr ]\;\sz{e}  ^{i (2n-1) \varphi} - i\;\bigl [ n R_{2,2n+1,r}(r, t) - \frac{i}{2} R_{2,2n+1,\varphi}(r, t)\bigr ]\;\sz{e}  ^{i (2n+1) \varphi}
\nonumber\\
\nonumber\\
\end{eqnarray}

We obtain by performing appropriate transformations:

\begin{eqnarray}\label{eqn304}
u_{2 r}^{*}(r,\varphi, t)\;=\;\bigl \{\bigl [ n R_{2,2n-1,r}(r, t) - \frac{i}{2} R_{2,2n-1,\varphi}(r, t)\bigr ] + \bigl [ n R_{2,2n+1,r}(r, t) - \frac{i}{2} R_{2,2n+1,\varphi}(r, t)\bigr ]\bigr \}\;\sz{e}  ^{i 2n \varphi}
\nonumber\\
\nonumber\\
\quad\quad\quad\quad
\end{eqnarray}

\begin{eqnarray}\label{eqn305}
u_{2\varphi}^{*}(r,\varphi, t)\;=\;i\;\bigl \{\bigl [ n R_{2,2n-1,r}(r, t) - \frac{i}{2} R_{2,2n-1,\varphi}(r, t)\bigr ] - \bigl [ n R_{2,2n+1,r}(r, t) - \frac{i}{2} R_{2,2n+1,\varphi}(r, t)\bigr ]\bigr \}\;\sz{e}  ^{i 2n \varphi}
\nonumber\\
\nonumber\\
\quad\quad\quad\quad
\end{eqnarray}

Here $u_{2 r}^{*}(r,\varphi, t) ,\; u_{2\varphi}^{*}(r,\varphi, t)$ are the radial and tangential components of the first correction $\vec u_{2}^{*}$ of the velocity $\vec u_{1}$. $R_{2,2n-1,r}(r, t),\; R_{2,2n+1,r}(r, t)$ are taken from formulas $(\ref{eqn266}), (\ref{eqn267})$.

\begin{eqnarray}\label{eqn306}
R_{2,2n-1,\varphi}(r, t)\;=\;R_{2,2n-1,\varphi, \varphi}(r, t) - 2niR_{2,2n-1,r}(r, t)
\quad\quad\quad\quad
\nonumber\\
\nonumber\\
R_{2,2n+1,\varphi}(r, t)\;=\;R_{2,2n+1,\varphi, \varphi}(r, t) - 2niR_{2,2n+1,r}(r, t)
\quad\quad\quad\quad
\nonumber\\
\nonumber\\
\end{eqnarray}

and $R_{2,2n-1,\varphi, \varphi}(r, t),\; R_{2,2n+1, \varphi, \varphi}(r, t)$ are taken from formulas $(\ref{eqn290}), (\ref{eqn291})$.

Then we do appropriate operations and have from formulas $(\ref{eqn304}), (\ref{eqn305})$:

\begin{eqnarray}\label{eqn304a}
u_{2 r}^{*}(r,\varphi, t)\;=\;- \frac{i}{2}\bigl [  R_{2,2n-1,\varphi,\varphi}(r, t) + R_{2,2n+1,\varphi,\varphi}(r, t)\bigr ]\;\sz{e}  ^{i 2n \varphi}
\nonumber\\
\quad\quad\quad\quad
\end{eqnarray}

\begin{eqnarray}\label{eqn305a}
u_{2\varphi}^{*}(r,\varphi, t)\;=\;\frac{1}{2}\bigl [  R_{2,2n-1,\varphi,\varphi}(r, t) - R_{2,2n+1,\varphi,\varphi}(r, t)\bigr ]\;\sz{e}  ^{i 2n \varphi}
\nonumber\\
\quad\quad\quad\quad
\end{eqnarray}

From formulas $(\ref{eqn302}), (\ref{eqn303})$ with properties of $R_{2,2n-1,r}(r, t),\; R_{2,2n+1,r}(r, t),\; R_{2,2n-1,\varphi}(r, t),\; R_{2,2n+1,\varphi}(r, t)$ it follows:

\begin{eqnarray}\label{eqn307}
  \begin{array}{ll} 								
lim \;\;u_{21}^{*}(r,\varphi,t)= 0;\;\;\;   lim \;\;u_{22}^{*}(r,\varphi,t)= 0;\\
t \rightarrow 0; \;\;\;\;\;\;\;\;\;\;\;\;\;\;\;\;\;\;\;\;\;\;\;\; t \rightarrow 0;
\end{array}
\end{eqnarray}

\begin{eqnarray}\label{eqn308}
  \begin{array}{ll} 								
lim \;\;u_{2r}^{*}(r,\varphi,t)= 0;\;\;\;   lim \;\;u_{2\varphi}^{*}(r,\varphi,t)= 0;\\
t \rightarrow 0; \;\;\;\;\;\;\;\;\;\;\;\;\;\;\;\;\;\;\;\;\;\;\;\; t \rightarrow 0;
\end{array}
\end{eqnarray}

In other words the velocity $\vec u_{2} = \vec u_{1} - \vec u_{2}^{*}$ [look at $(\ref{eqn50})$] satisfies the initial conditions $(\ref{eqn200})$.
We use the asymptotic properties of the  confluent hypergeometric function $\Phi(a, c; x)$ and have from formulas $(\ref{eqn266}),(\ref{eqn267}), (\ref{eqn290}),$ $(\ref{eqn291})$: the first correction $\vec u_{2}^{*}$ and therefore the velocity $\vec u_{2}$ satisfies conditions  $\;(\ref{eqn16})\;$ ( for $r \;\rightarrow\; \infty $).

By comparing the solution  $\vec u_{1}$ from $(\ref{eqn238}), (\ref{eqn239})$ or $(\ref{eqn240}), (\ref{eqn241})$ of the first step of iterative process  with the first correction $\vec u_{2}^{*}$ from $(\ref{eqn302}), (\ref{eqn303})$ or $(\ref{eqn304a}), (\ref{eqn305a})$, which is received on the second step of iterative process, we see that

\begin{equation}\label{eqn309}
\mid\vec{u_{2}^{*}}\mid\; <<\; \mid\vec{u_{1}}\mid
\end{equation}

with conditions 

\begin{eqnarray}\label{eqn310}
F_{n}\;\leq\;\frac{1}{n} 
\nonumber\\
\nonumber\\
\end{eqnarray}

By continuing this iterative process we can obtain next parts $\;\vec{u_{3}^{*}}\;, \vec{u_{4}^{*}}...$, of the converging series for $\vec{u}$. For arbitrary step j  of the iterative process we have by using formula $(\ref{eqn57})$:

\begin{equation}\label{eqn311}
\vec{u}_{j}\;=\;\vec{u}_{1}\;-\;\sum_{l=2}^{j} \vec{u}_{l}^{*}\
\end{equation}

and then: 

\begin{equation}\label{eqn312}
  \begin{array}{ll} 								
lim \;\;\vec{u_{j}}= \vec{u}\\   
j \rightarrow \infty 
\end{array}
\end{equation}

where $\vec{u}$ is the solution of the problem $(\ref{eqn1}) - (\ref{eqn6})$.
\nonumber\\
\nonumber\\
Below we provide numerical analysis of these results for the following values of problem's parameters:
 
Circumferential modes n = 1, 2, 3, 4, 5.

$\mu_n$ = 1 (n = 1, 2, 3, 4, 5).

Results were obtained for functions 

$\vec u_{1} - (\ref{eqn238}), (\ref{eqn239})$ or $(\ref{eqn240}), (\ref{eqn241})$ with calculations of the confluent hypergeometric functions $\cite{BE153}$; 

$\vec u_{2}^{*} - (\ref{eqn302}), (\ref{eqn303})$ or $(\ref{eqn304a}), (\ref{eqn305a})$ by using numerical integration of the triple integrals $(\ref{eqn290}), (\ref{eqn291})$. Each of those integrals is computed as an iterated integral.

Let us consider first the calculation of the inner integrals from $(\ref{eqn290}), (\ref{eqn291})$:

\begin{eqnarray}\label{eqn290a}
I_{\underline{+}}(r, \tilde r, t, \tau) = \int_{0}^{\infty}\sz{e}  ^{-\nu \rho^2 (t-\tau)}J_{2n_{\underline{+}}1}(r \rho)\cdot J_{2n}(\tilde r \rho)\;d\rho 
\end{eqnarray}

For condition t $< \tau$ the integrand is diminishing fast enough. It is easy to find upper limit of integration, so we can substitute integral $(\ref{eqn290a})$ for 

\begin{eqnarray}\label{eqn290b}
I_{\underline{+}}(r, \tilde r, t, \tau) = \int_{0}^{A_1}\sz{e}  ^{-\nu \rho^2 (t-\tau)}J_{2n_{\underline{+}}1}(r \rho)\cdot J_{2n}(\tilde r \rho)\;d\rho, 
\end{eqnarray}

where ($0 < A_1 = 200 < \infty$), and hence we are integrating over the finite interval. For additional check let us increase $A_1$ in 1.5 times and change the number of steps of integration $n_1($ from 4001 to 6001 $)$. We see that the difference in result values of integral $(\ref{eqn290b})$ is within the range of required precision $\epsilon_1(10^{-14})$.

For condition t $= \tau$ the integral $(\ref{eqn290a})$ is in fact an integral of Weber and Schafheitlin  and it is possible to calculate it analytically $(\ref{A17})\; \cite{gW44}$.

Let us now consider the calculation of middle integrals 

\begin{eqnarray}\label{eqn290c}
\tilde I_{\underline{+}}(r, t, \tau) = \int_{0}^{\infty} T_{2,2n,\varphi, \varphi}(\tilde r, \tau) I_{\underline{+}}(r, \tilde r, t, \tau) d \tilde r  
\end{eqnarray}

We use asymptotical properties of  confluent hypergeometric functions $\Phi(a, c; x)$ $ \;\;\cite{BE153}$ and we have:

\begin{equation}
  \begin{array}{ll} 
	T_{2,2n,\varphi, \varphi}(\tilde r, \tau)\rightarrow (1/ \tilde r^{2n + 5})\\
				\tilde r \rightarrow \infty
					\end{array}     
\end{equation}

Hence, we substitute integral $(\ref{eqn290c})$ for 

\begin{eqnarray}\label{eqn290d}
\tilde I_{\underline{+}}(r, t, \tau) = \int_{0}^{A_2} T_{2,2n,\varphi, \varphi}(\tilde r, \tau) I_{\underline{+}}(r, \tilde r, t, \tau) d \tilde r  
\end{eqnarray}

where ($0 < A_2 = 20 < \infty$) and integration is really over the finite interval. For additional check let us increase value $A_2$ in 1.5 times and change the number of integration steps $n_2$(from 201 to 301). We see that the difference in result values of integral $(\ref{eqn290d})$ is within the range of required precision $\epsilon_2(10^{-11})$.

Confluent hypergeometric functions $\Phi(a, c; x)$ were computed with precision $\epsilon(10^{-15})$.

The outer integrals in $(\ref{eqn290}), (\ref{eqn291})$ are the integrals over finite interval (0, t = 10). These integrals are computed with precision $\epsilon_3(10^{-5})$ and the number of steps of integration $n_3$= 101. For additional check let us change the number of integration steps $n_3$(from 101 to 201), and we see that the difference in result integral values is within the required precision $\epsilon_3(10^{-5})$. 
All integrals were computed by Simpson's method and $\epsilon_1(10^{-14}) < \epsilon_2(10^{-11}) < \epsilon_3(10^{-5})$.

$\vec u_{2} = \vec u_{1} - \vec u_{2}^{*}$ and is shown in FIG. 5.1.1 - 5.1.15.
The vector field $\vec u_{2}$ at distances r = 1, 2, 3, 5, 7 is represented by the dotted curves in top diagrams.
The comparison of $\mid\vec u_{1}\mid$ (dashed plots) and $\mid\vec u_{2}^{*}\mid$ (solid plots) in plane $\varphi$ = [0, $\pi$], at distances $\;\;\;$        0 $\leq$ r $\leq$ 50 is represented in bottom diagrams. This comparison shows $\mid\vec{u_{2}^{*}}\mid\; <<\; \mid\vec{u_{1}}\mid$ 
and is corresponding to the conclusion $(\ref{eqn484})$.
$\nonumber\\$
$\nonumber\\$
   \begin{center}   
     \includegraphics[height=80mm]{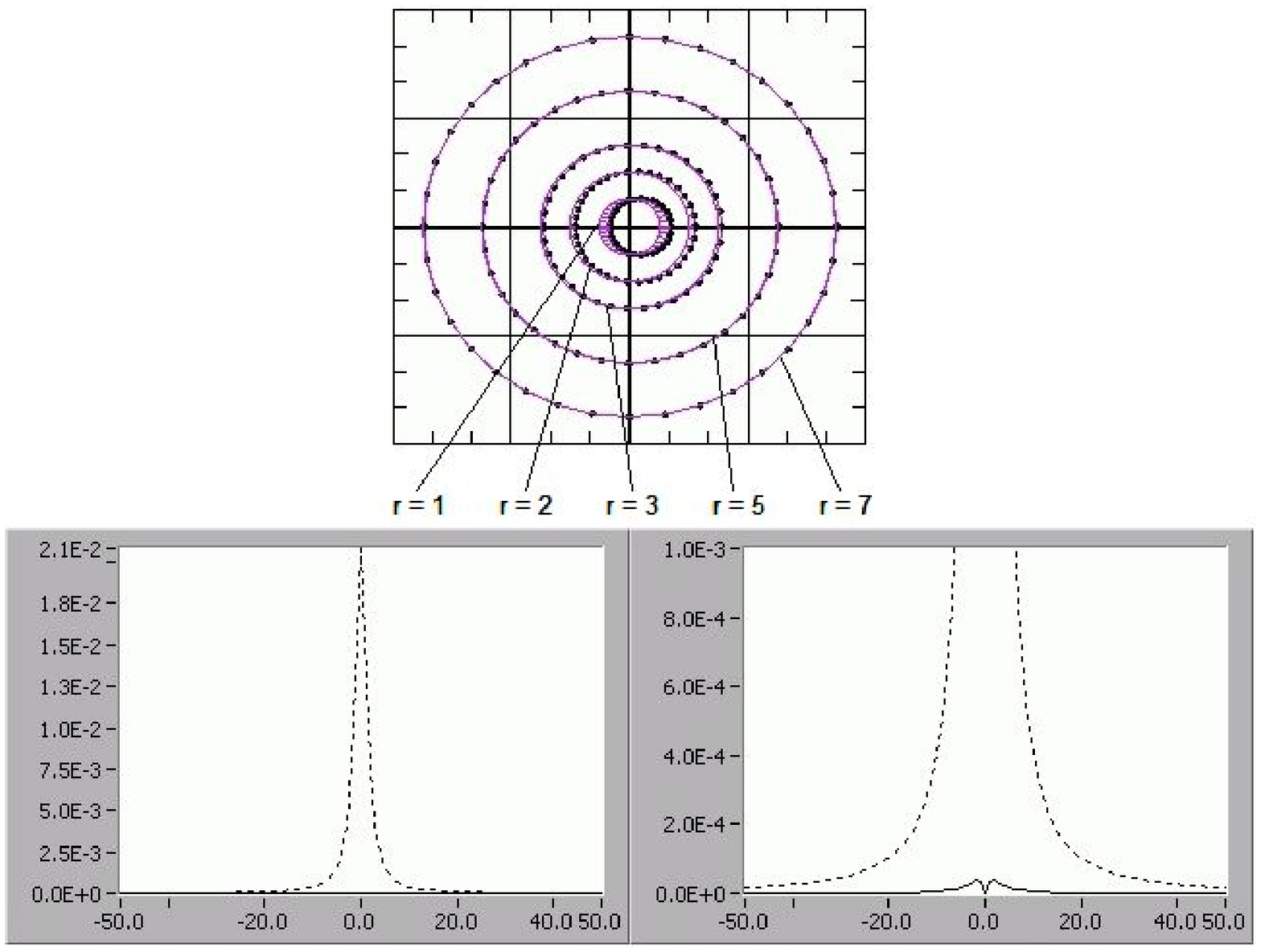}\\
     FIG.5.1.1. n = 1, $F_1$ = 1, $\nu$ = 1.5
   \end{center}
  \begin{center}   
     \includegraphics[height=80mm]{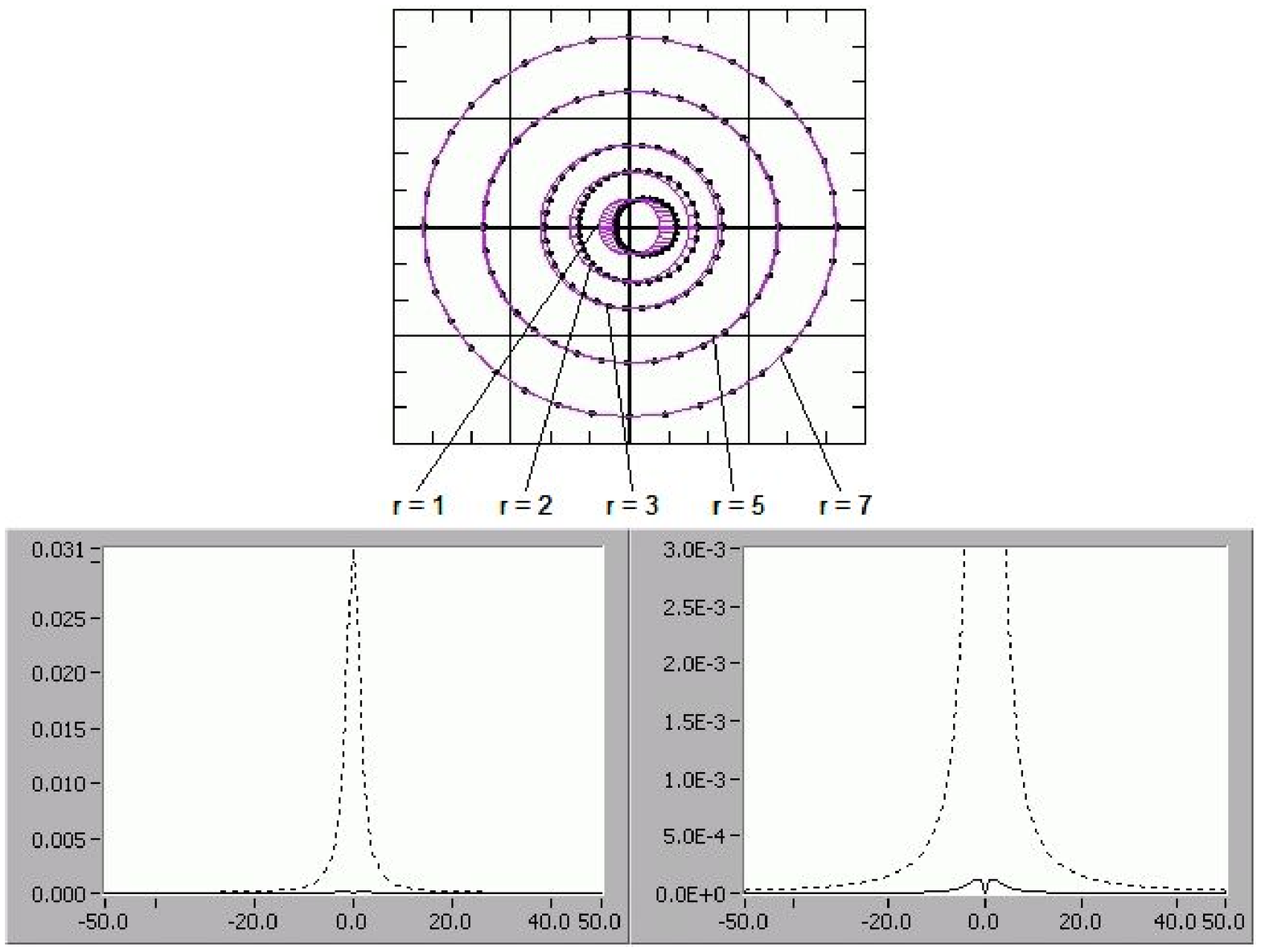}\\
      FIG.5.1.2. n = 1, $F_1$ = 1, $\nu$ = 1
   \end{center} 
 \begin{center}  
     \includegraphics[height=80mm]{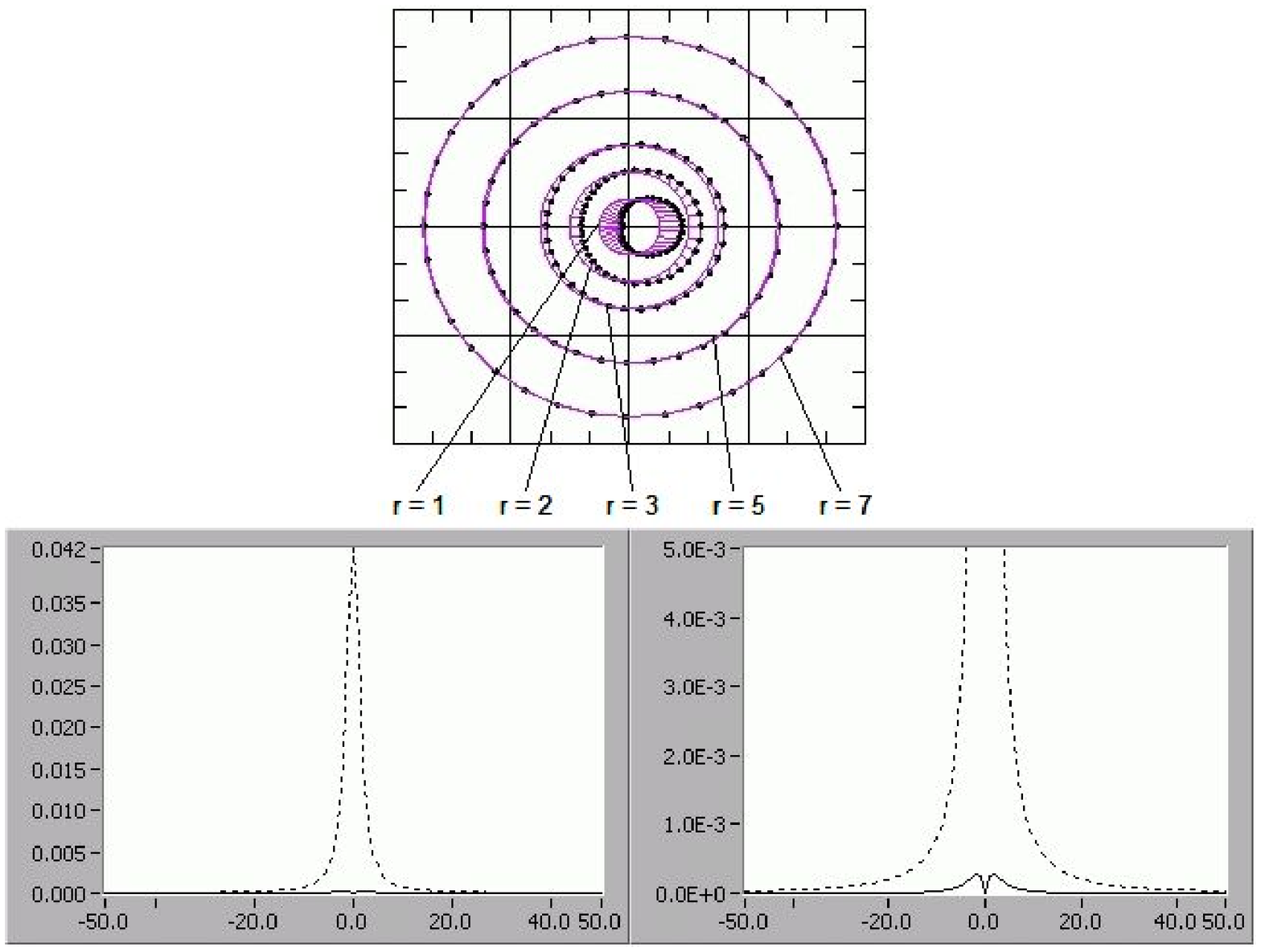}\\
       FIG.5.1.3. n = 1, $F_1$ = 1, $\nu$ = 0.75
   \end{center}        
  \begin{center} 
    \includegraphics[height=80mm]{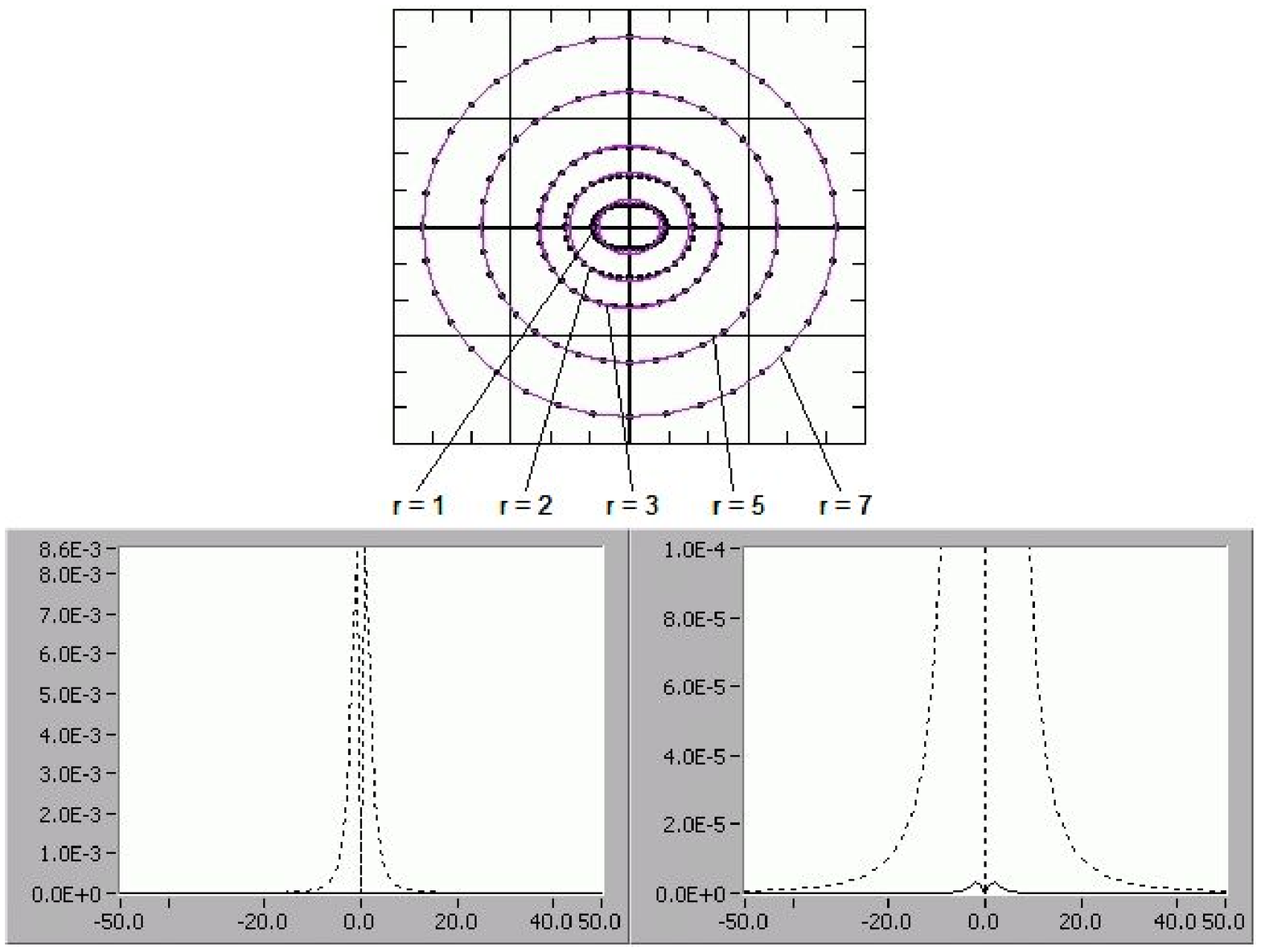}\\
      FIG.5.1.4. n = 2, $F_2$ = 0.5, $\nu$ = 1.5
  \end{center} 
  \begin{center} 
    \includegraphics[height=80mm]{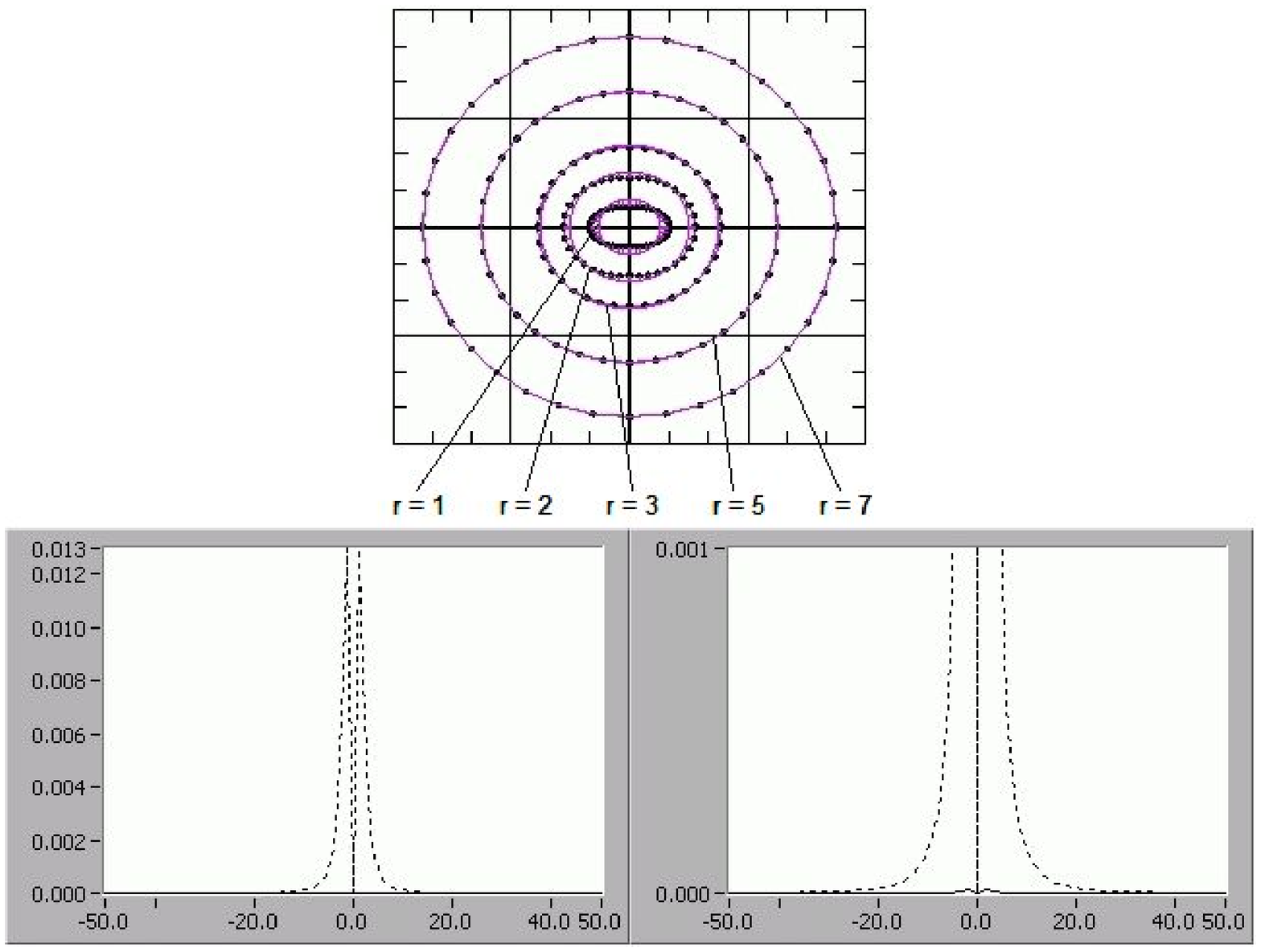}\\
     FIG.5.1.5. n = 2, $F_2$ = 0.5, $\nu$ = 1
  \end{center} 
  \begin{center} 
    \includegraphics[height=80mm]{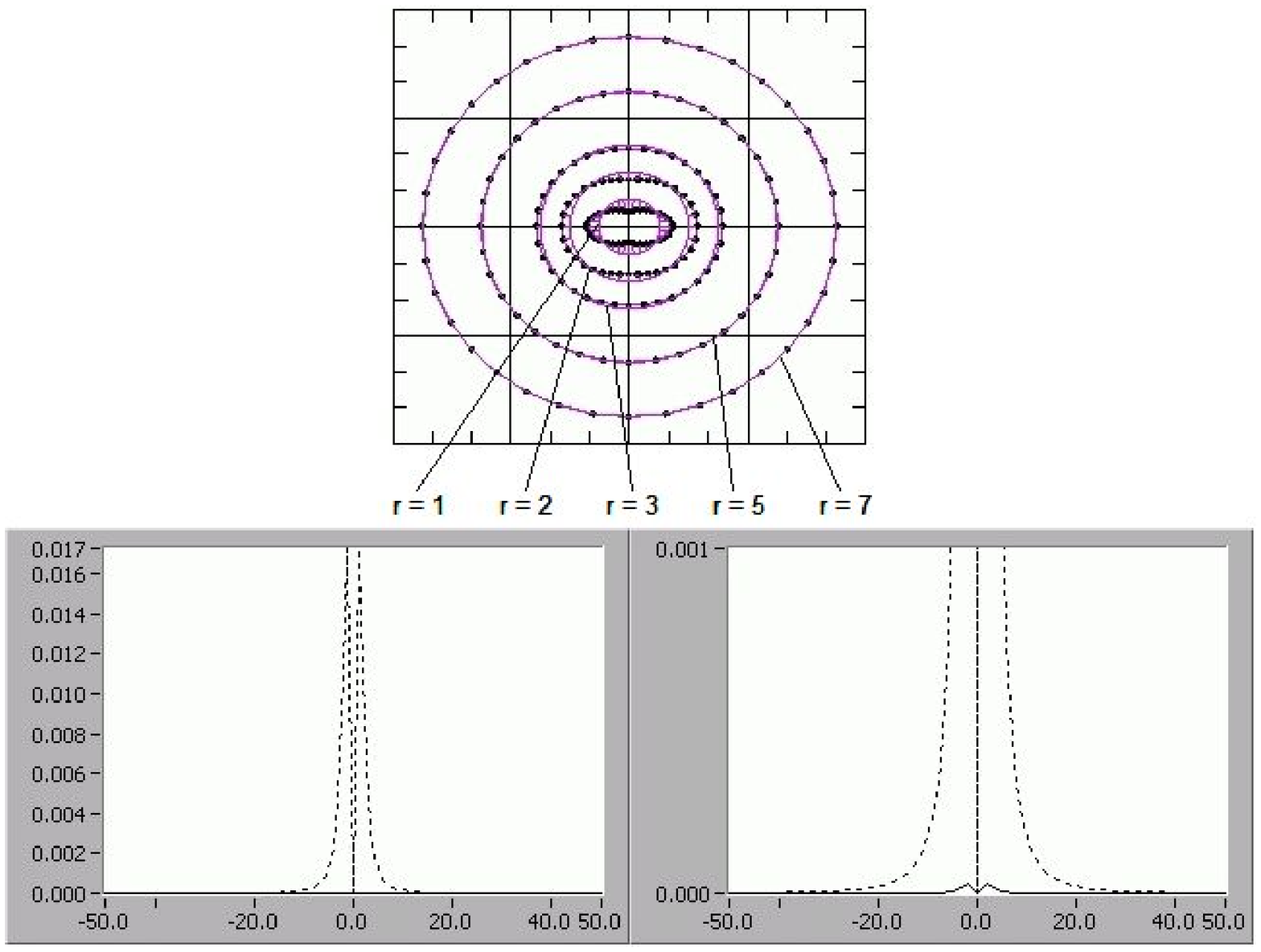}\\
     FIG.5.1.6. n = 2, $F_2$ = 0.5, $\nu$ = 0.75
  \end{center}        
    \begin{center}  
    \includegraphics[height=80mm]{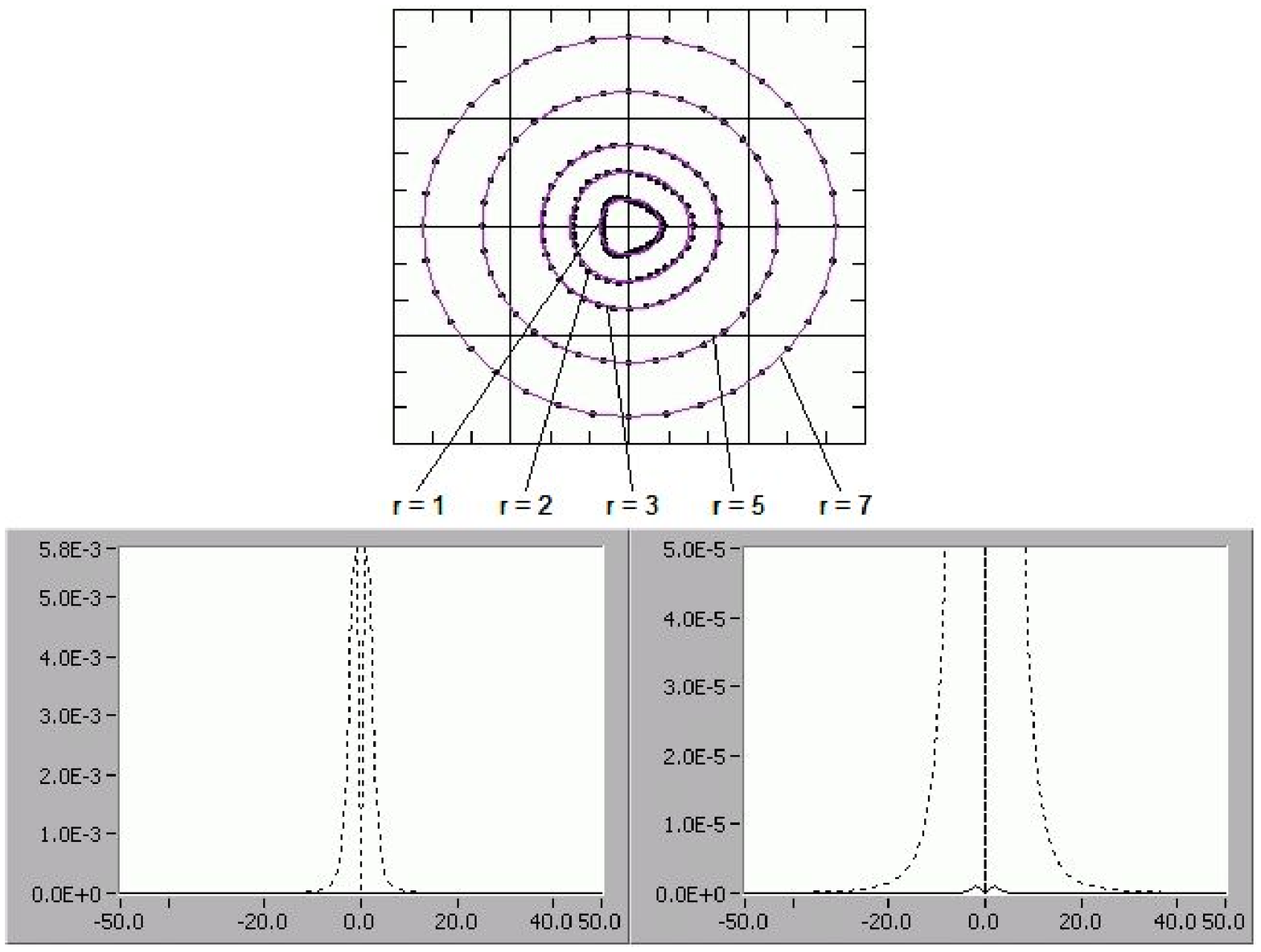}\\
      FIG.5.1.7. n = 3, $F_3$ = 0.33, $\nu$ = 1.5
      \end{center}  
   \begin{center}   
    \includegraphics[height=80mm]{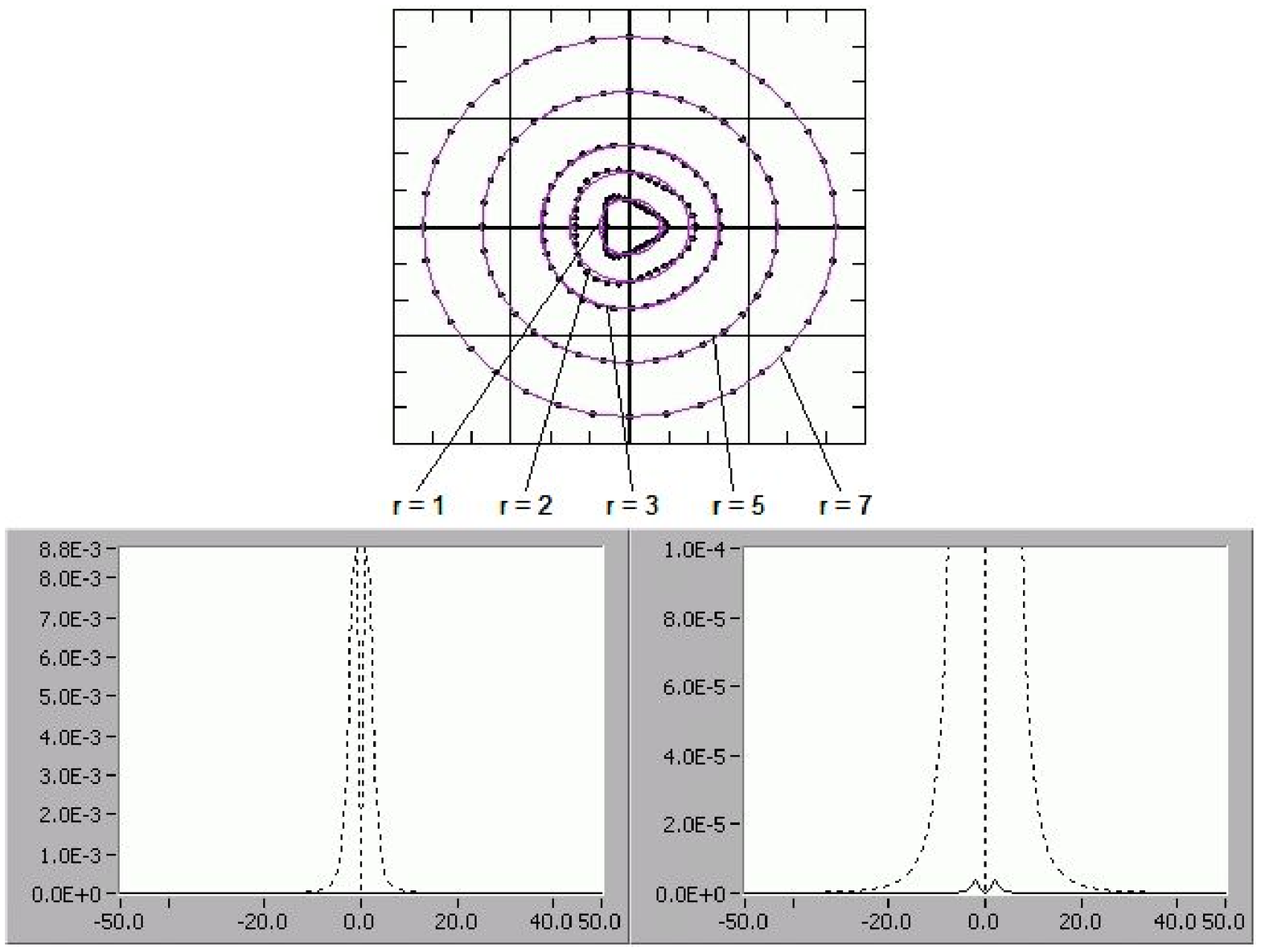}\\
     FIG.5.1.8. n = 3, $F_3$ = 0.33, $\nu$ = 1
      \end{center}
   \begin{center}  
    \includegraphics[height=80mm]{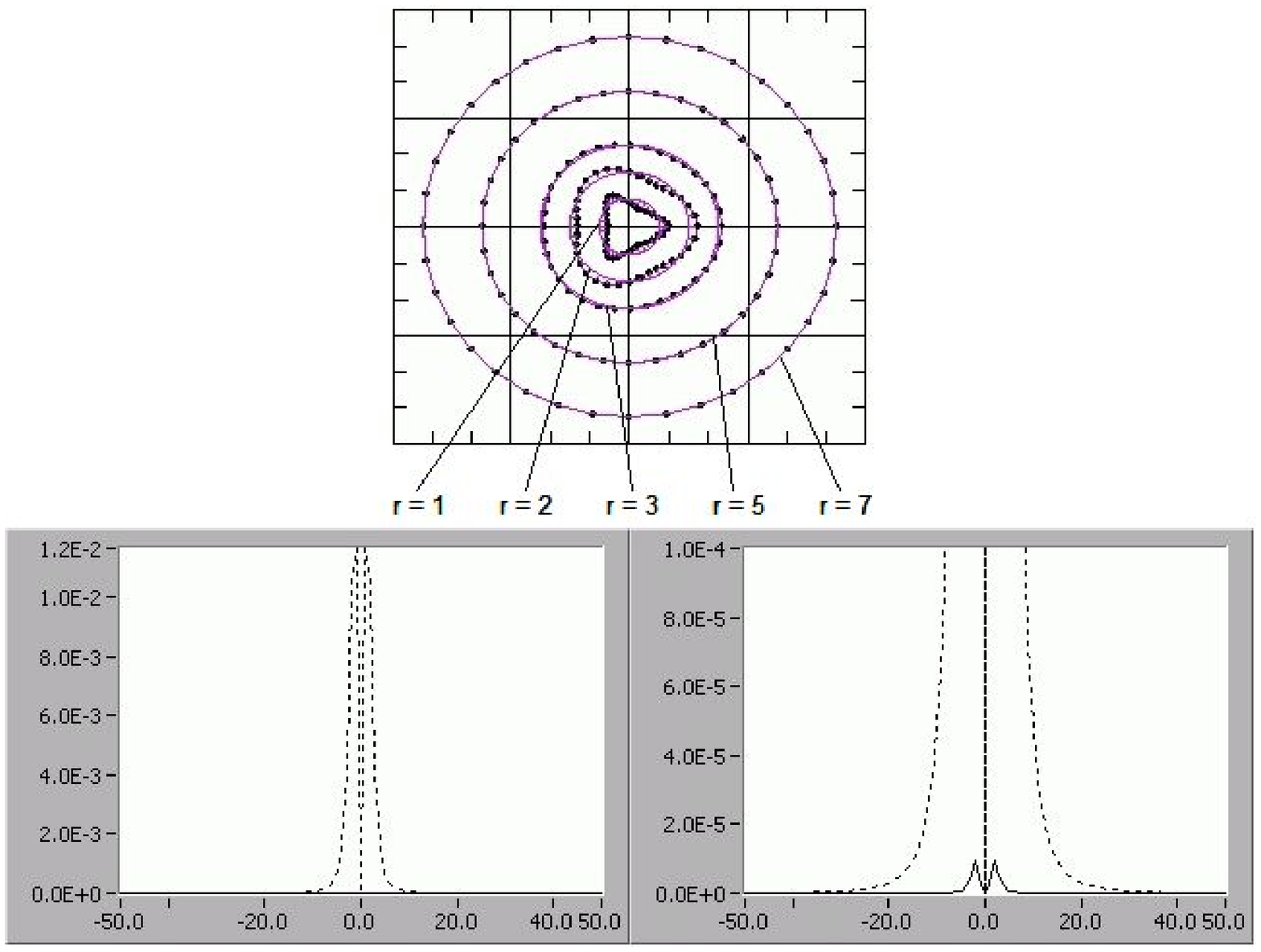}\\
      FIG.5.1.9. n = 3, $F_3$ = 0.33, $\nu$ = 0.75
      \end{center}              
 \begin{center} 
   \includegraphics[height=80mm]{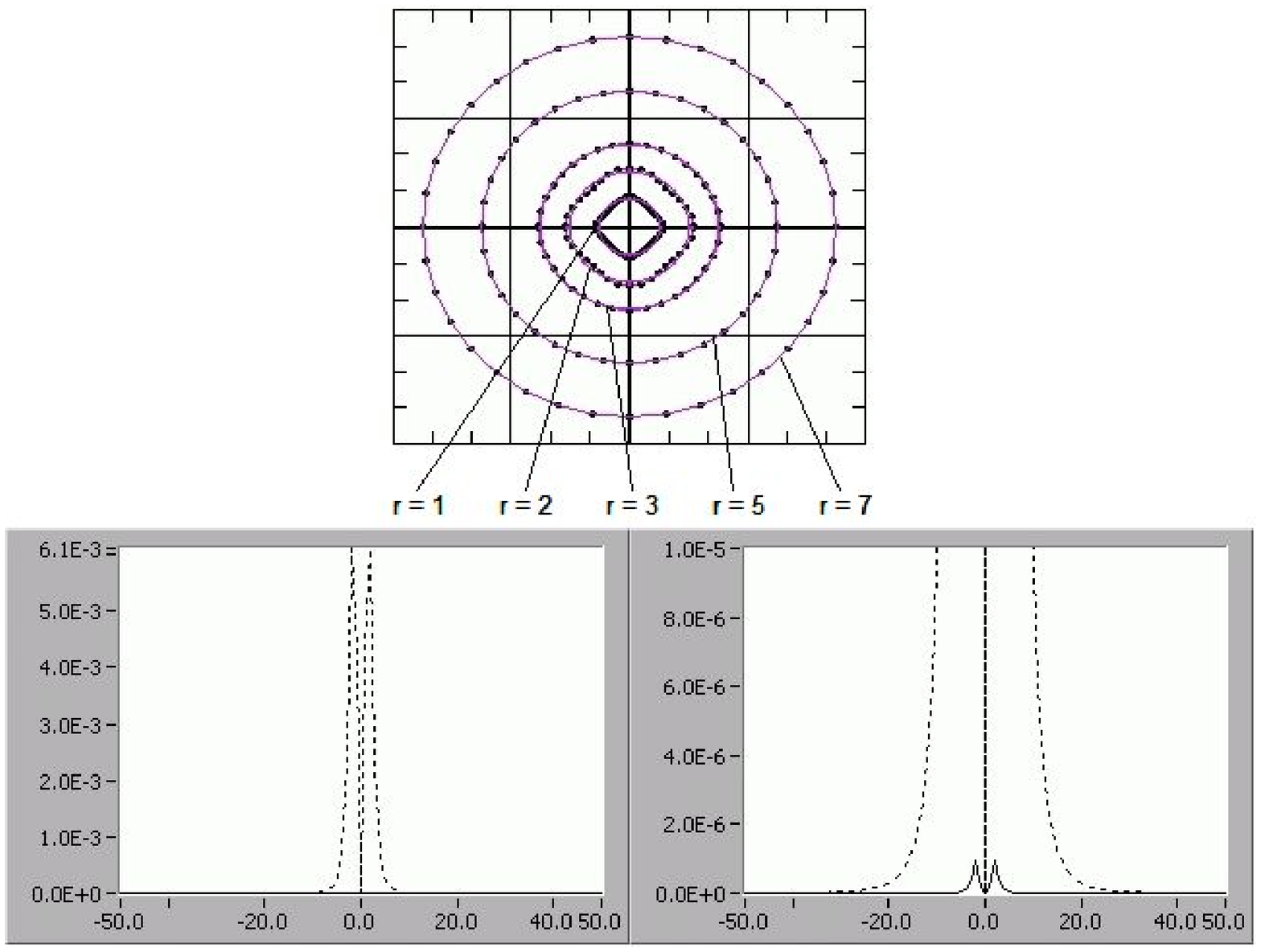}\\
    FIG.5.1.10. n = 4, $F_4$ = 0.25, $\nu$ = 1.5
  \end{center}
 \begin{center} 
   \includegraphics[height=80mm]{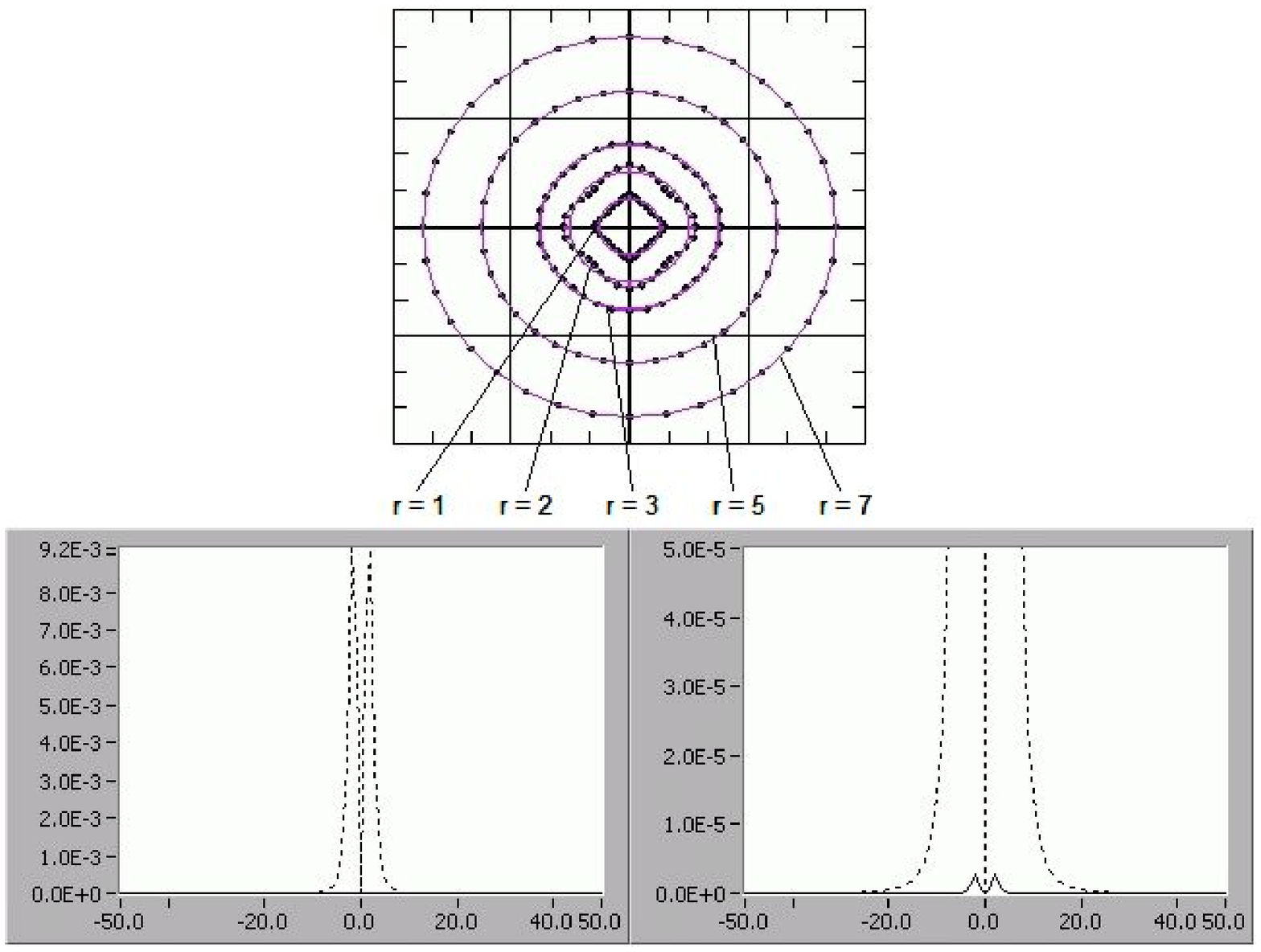}\\
    FIG.5.1.11. n = 4, $F_4$ = 0.25, $\nu$ = 1
  \end{center}    
 \begin{center} 
   \includegraphics[height=80mm]{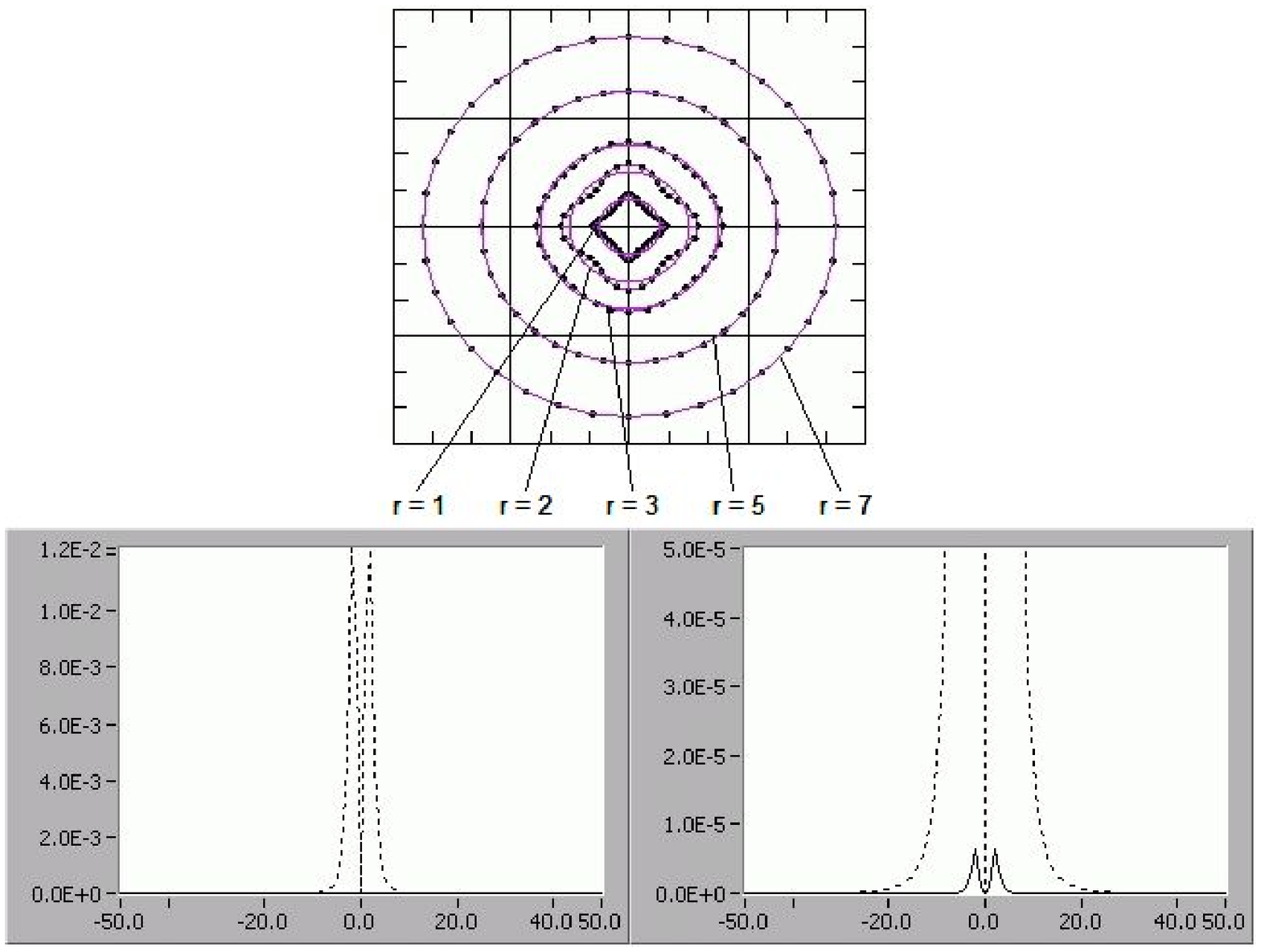}\\
    FIG.5.1.12. n = 4, $F_4$ = 0.25, $\nu$ = 0.75
  \end{center}            
 \begin{center}  
   \includegraphics[height=80mm]{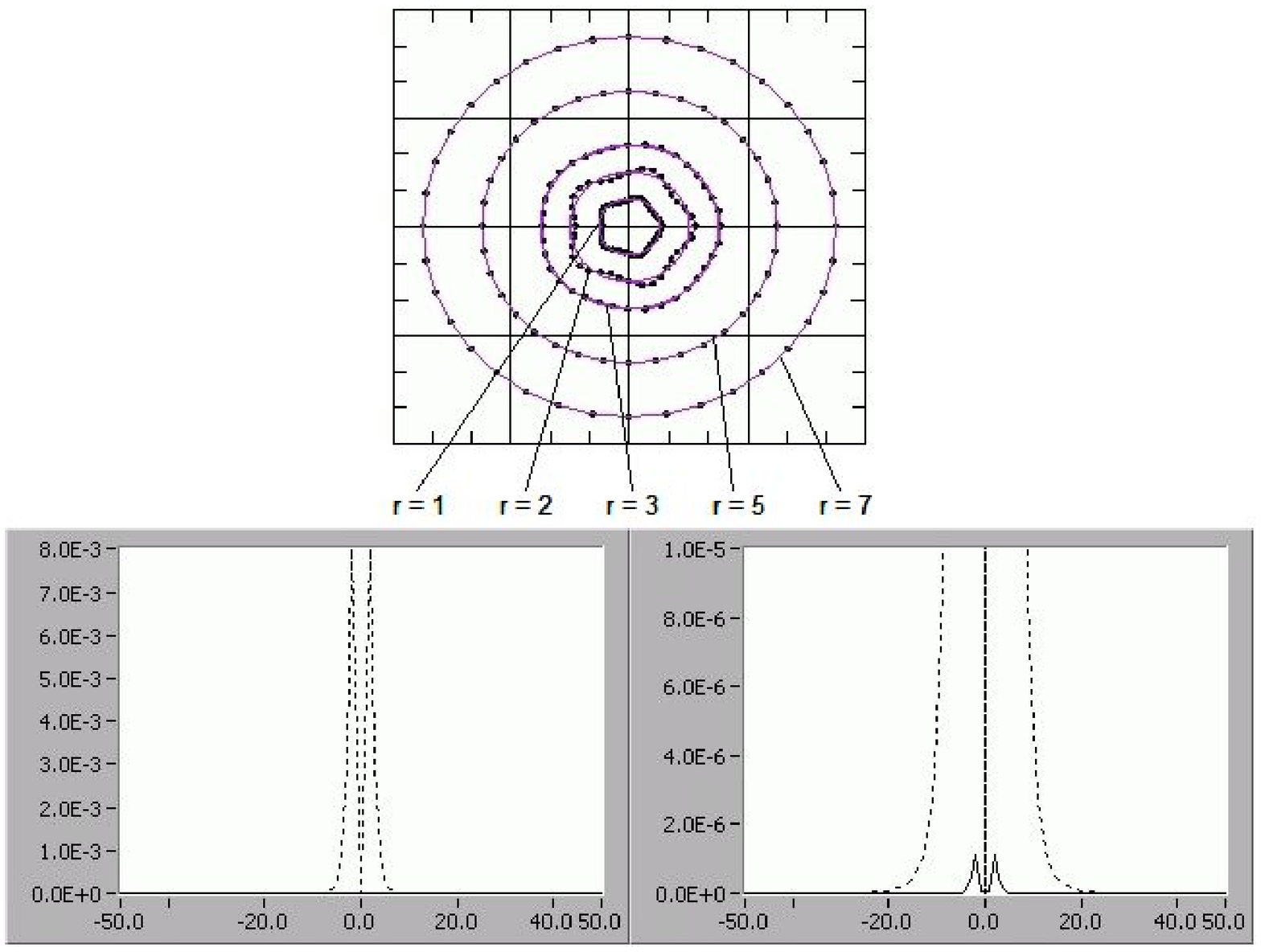}\\ 
   FIG.5.1.13. n = 5, $F_5$ = 0.2, $\nu$ = 1.5
 \end{center}  
 \begin{center} 
   \includegraphics[height=80mm]{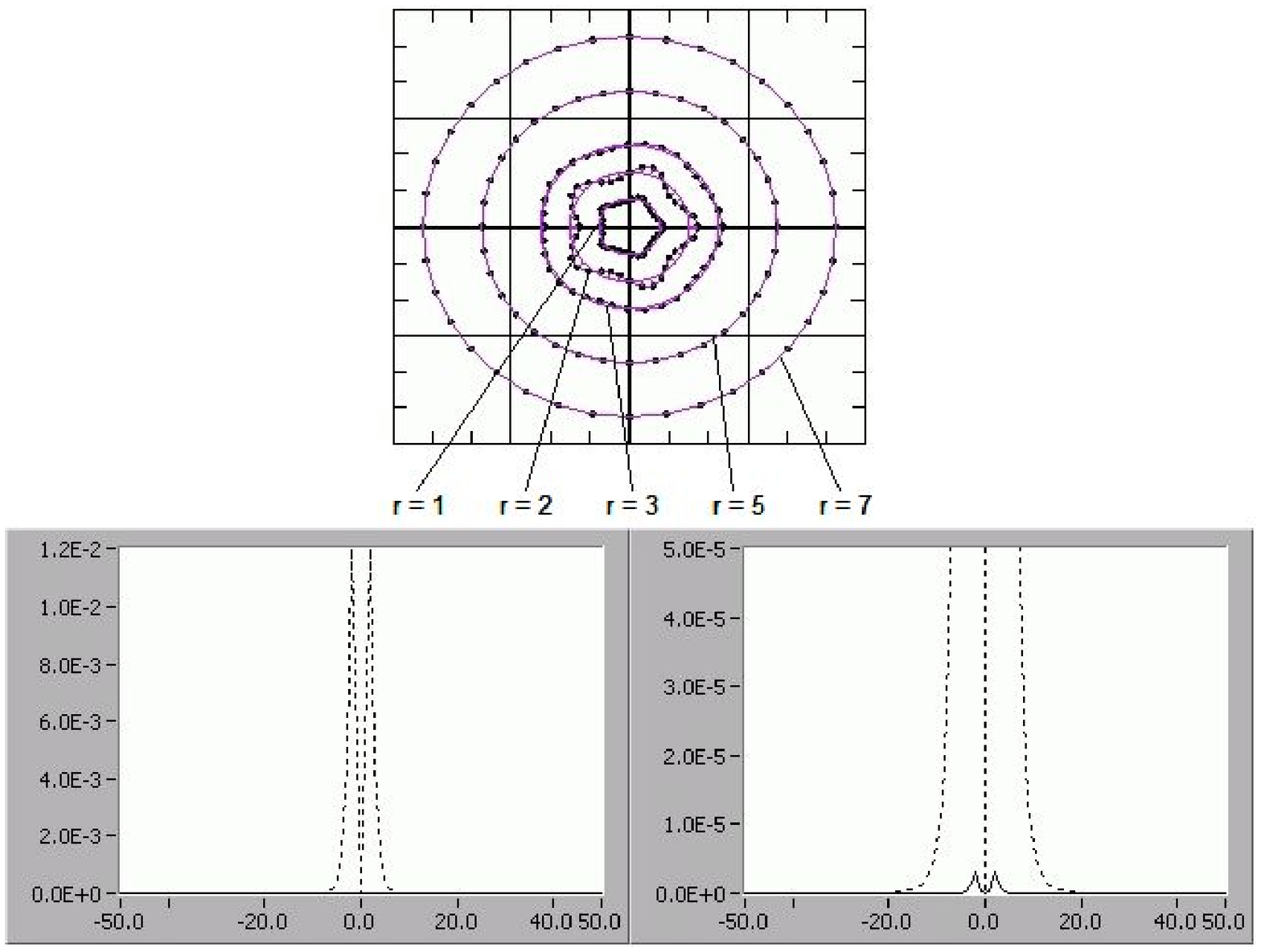}\\
     FIG.5.1.14. n = 5, $F_5$ = 0.2, $\nu$ = 1
 \end{center}  
 \begin{center}  
  \includegraphics[height=80mm]{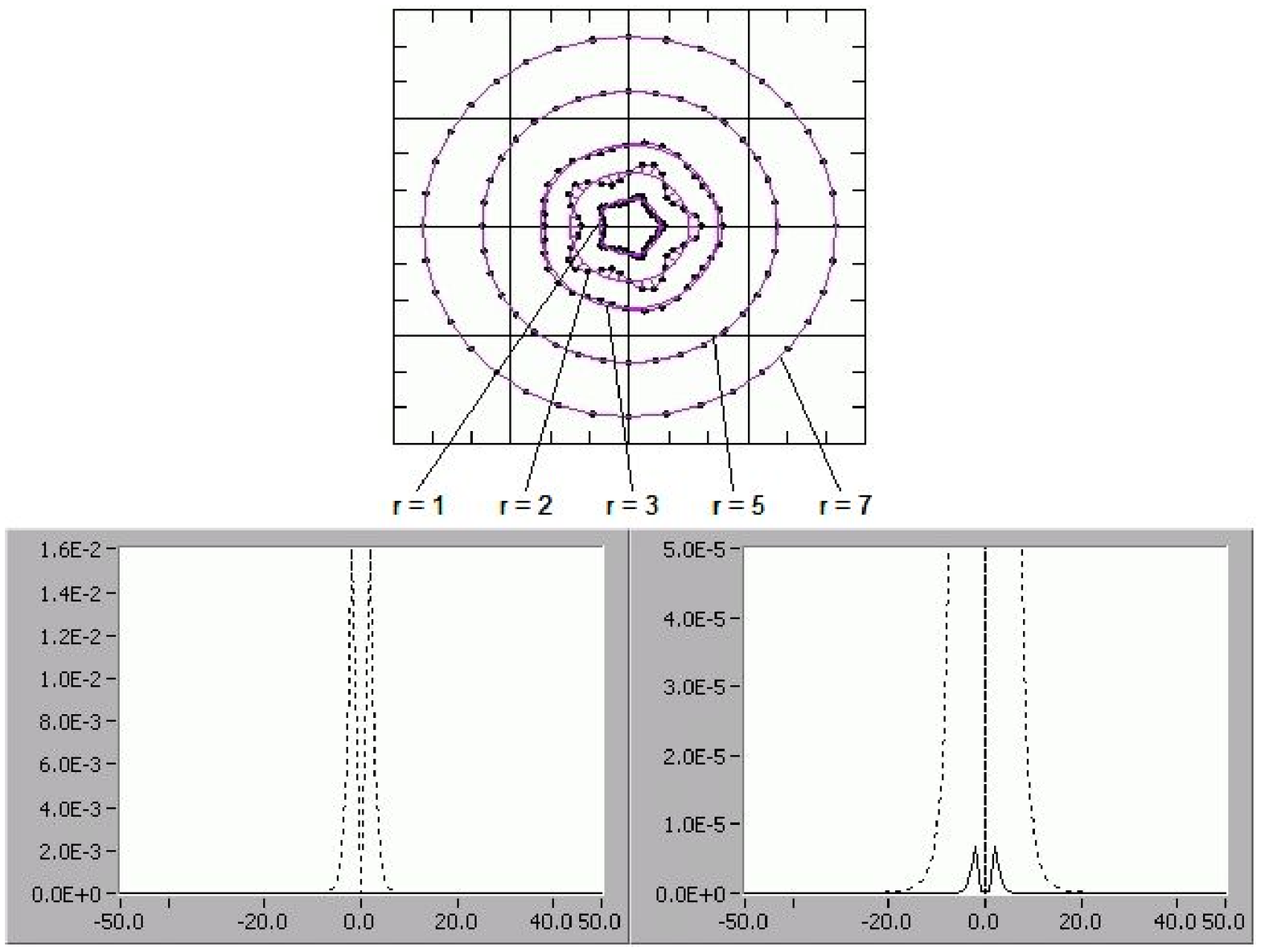}\\
   FIG.5.1.15. n = 5, $F_5$ = 0.2, $\nu$ = 0.75
 \end{center} 
 $\nonumber\\
\nonumber\\$
 $\nonumber\\
\nonumber\\$
 $\nonumber\\
\nonumber\\$
 $\nonumber\\
\nonumber\\$

\appendix\section*{}

The Fourier integral can be stated in the forms: 

$\\N=1$

\begin{equation}\label{A1}
U(\gamma)=F[u(x)]= \frac{1}{(2\pi)^{1/2}} \int_{-\infty}^{\infty}u(x)\, \sz{e}  ^{i\gamma x}dx  \;\;\;\;\;\;\;\;\;\; u(x)= \frac{1}{(2\pi)^{1/2}}\int_{-\infty}^{\infty}U(\gamma) \sz{e}  ^{-i\gamma x}d\gamma
\end{equation}

$\\\\\\N=2$

\begin{eqnarray}\label{A2}
U( \gamma_{1} , \gamma_{2})=F[ u(x_{1} , x_{2})]= \frac{1}{2\pi} \int_{-\infty}^{\infty} \int_{-\infty}^{\infty} u( x_{1} , x_{2})\,\sz{e}  ^{ i( \gamma_{1} x_{1} + \gamma_{2} x_{2}) } dx_{1} dx_{2} 
\nonumber\\
\nonumber\\
u( x_{1} , x_{2})= \frac{1}{2\pi} \int_{-\infty}^{\infty} \int_{-\infty}^{\infty} 
U( \gamma_{1} , \gamma_{2})\,\sz{e}  ^{- i( \gamma_{1} x_{1} + \gamma_{2} x_{2}) } d\gamma_{1} d\gamma_{2}
\end{eqnarray}

$\\N=3$

\begin{eqnarray}\label{A3}
U( \gamma_{1} , \gamma_{2} , \gamma_{3})=F[\, u(x_{1} , x_{2} , x_{3})]= \frac{1}{(2\pi)^{3/2}} \int_{-\infty}^{\infty} \int_{-\infty}^{\infty} \int_{-\infty}^{\infty} u( x_{1} , x_{2} , x_{3})\,\sz{e}  ^{ i( \gamma_{1} x_{1} + \gamma_{2} x_{2} + \gamma_{3} x_{3}) } dx_{1} dx_{2} dx_{3} 
\nonumber\\
\nonumber\\
u( x_{1} , x_{2} , x_{3})= \frac{1}{(2\pi)^{3/2}} \int_{-\infty}^{\infty} \int_{-\infty}^{\infty} \int_{-\infty}^{\infty} U( \gamma_{1} , \gamma_{2},\gamma_{3}  )\, \sz{e}  ^{- i( \gamma_{1} x_{1} + \gamma_{2} x_{2} + \gamma_{3} x_{3}) } d\gamma_{1} d\gamma_{2} d\gamma_{3} 
\end{eqnarray}
$\\\\$
The Laplace integral is usually stated in the following form:

\begin{equation}\label{A4}
U^{\otimes}(\eta)=L[\,u(t)\,]= \int_{0}^{\infty}u(t)\, \sz{e}  ^{-\eta t}dt
\;\;\;\;\; u(t)=\frac{1}{2\pi i}\int_{c- i \infty }^{c + i \infty} U^{\otimes}(\eta) \,\sz{e}  ^{\eta t}d\eta \;\;\;\;\; c > c_{0}
\end{equation}

\begin{equation}\label{A5}
L[\,u^{'}(t)\,]=\eta \,U^{\otimes}(\eta)-u(0)
\end{equation}

\newtheorem{guess}{The convolution theorem}[section]
\begin{guess}
$\\\\$
If integrals
\[ U_{1}^{\otimes}(\eta)= \int_{0}^{\infty}u_{1}(t)\, \sz{e}  ^{-\eta t}d\,t  \;\;\;\;\;\;\;\;\;\; U_{2}^{\otimes}(\eta)= \int_{0}^{\infty}u_{2}(t)\, \sz{e}  ^{-\eta t}d\,t \]

absolutely converge by $Re\, \eta > \sigma_{d}$, then  $U^{\otimes}(\eta)\,= \,U_{1}^{\otimes}(\eta)\, U_{2}^{\otimes}(\eta)$ is Laplace transform of 

\begin{equation}\label{A6}
u(t)=\int_{0}^{t}u_{1}(t-\tau)\,u_{2}(\tau)\,d\,\tau
\\
\end{equation}
\end{guess}

Useful \emph{Laplace integral}:

\begin{equation}\label{A7}
L[\,\sz{e}  ^{\eta_{k}t}\,]\,=\,\int_{0}^{\infty}\sz{e}  ^{-(\eta-\eta_{k})\,t}d\,t
\;=\; \frac{1}{(\eta-\eta_{k})}\;\;\;\;\;\;\;\;\;(Re\,\eta\,>\,\eta_{k})
\end{equation}

De Moivre's formulas:

\begin{equation}\label{A8}
cos \varphi = \frac{1}{2}(\sz{e}^{i\varphi} + \sz{e}^{-i\varphi})\;\;,\;\;sin \varphi = \frac{1}{2i}(\sz{e}^{i\varphi} - \sz{e}^{-i\varphi})
\end{equation}

Bessel function's integral representation:

\begin{equation}\label{A9}
J_{n}(z) = \frac{i^{-n}}{2\pi}\oint_C\sz{e}  ^{iz cos\theta + in\theta} d \theta\;, \;\;\emph{C is the unit circle around the origin.}
\end{equation}

The discontinuous integral of Weber and Schafheitlin:

\begin{equation}\label{A17}
\int_{0}^{\infty}J_{\mu}(at)\cdot J_{\mu-1}(bt)\;dt\;\;=\;\;\Biggl\{
                \begin{array}{lll}   
                b^{\mu-1}/a^{\mu} & (b < a)\\ 
                1/2b & (b = a)\\                            
                 0 &  (b > a)            
                \end{array}
\end{equation}

\begin{equation}\label{A10}
\int_{0}^{\infty}J_{\mu}(\alpha t) \sz{e}^{-\gamma^2 t^2}t^{\mu+1}\;d t = \frac{ \alpha^\mu}{(2\gamma^2)^{\mu+1}}\sz{e}^{-\frac{\alpha^2}{4\gamma^2}},\;\;\;\;\;\;\;\;\;
Re \;\mu > -1, \; Re \;\gamma^2 > 0.
\end{equation} 

\begin{eqnarray}\label{A11}
\int_{0}^{\infty}J_{\mu}(\alpha t) \sz{e}^{-\gamma^2 t^2}t^{\rho-1}\;d t = \frac{\gamma^{-\rho}}{2\cdot\Gamma(\mu+1)} \cdot \Gamma\biggl(\frac{\mu+\rho}{2} \biggr)\cdot \biggl(\frac{\alpha}{2\gamma} \biggr)^{\mu}\cdot \Phi \biggl(\frac{\mu+\rho}{2},\mu+1;-\frac{\alpha^2}{4\gamma^2}\biggr),
\nonumber\\
\nonumber\\
Re \;\gamma^2 > 0, Re (\mu+\rho) > 0.
\end{eqnarray} 

\begin{equation}\label{A12}
\frac{d}{dy}[y^a \cdot \Phi (a,c;-\beta y)] = a \cdot y^{a-1} \cdot \Phi (a+1,c;-\beta y)
\end{equation}

\begin{equation}\label{A13}
\frac{d}{dy}[y^{c-1} \cdot \Phi (a,c;-\beta y)] = (c - 1) \cdot y^{c-2} \cdot \Phi (a,c-1;-\beta y)
\end{equation}

Formula describing connection between the contiguous confluent hypergeometric functions:

\begin{equation}\label{A14}
c \cdot\Phi - c \cdot\Phi(a-) - x\cdot c \cdot\Phi(c+) = 0
\end{equation}
\nonumber\\
\nonumber\\
\nonumber\\
\nonumber\\

\end{document}